\documentclass[reqno]{amsart}
\usepackage[margin = 1.5in]{geometry}
\usepackage{amsmath, amssymb, amsthm, fancyhdr, verbatim, graphicx}
\usepackage{enumerate}
\usepackage[all]{xy}
\usepackage[usenames,dvipsnames]{xcolor}
\usepackage{mathrsfs}
\usepackage{tikz-cd}
\usepackage[T1]{fontenc} 

\usepackage[hyperindex,breaklinks]{hyperref}

\usepackage{framed}
\usepackage[OT2, T1]{fontenc}
\usepackage[titletoc]{appendix}
\usepackage{bbm}

\usepackage{adjustbox}

\usepackage{amssymb} 
\def\acts{\curvearrowright}

\numberwithin{equation}{section}

\DeclareSymbolFont{cyrletters}{OT2}{wncyr}{m}{n}
\DeclareMathSymbol{\Sha}{\mathalpha}{cyrletters}{"58}

\usepackage{color}

\newcommand{\F}{\mathbf{F}}

\newcommand{\CC}{\mathbf{C}}
\newcommand{\G}{\mathbf{G}}

\newcommand{\wt}[1]{\widetilde{#1}}

\newcommand{\Q}{\mathbf{Q}}
\newcommand{\Z}{\mathbf{Z}}

\newcommand{\mf}[1]{\mathfrak{#1}}

\newcommand{\Gal}{\operatorname{Gal}}

\newcommand{\ul}[1]{\underline{#1}}
\newcommand{\ol}[1]{\overline{#1}}
\newcommand{\wh}[1]{\widehat{#1}}

\newcommand{\mbb}[1]{\mathbb{#1}}

\newcommand{\Cal}[1]{\mathcal{#1}}
\newcommand{\A}{\mathbf{A}}

\newcommand{\ft}{{}^{\tau}} 
\newcommand{\co}{\colon}
\newcommand{\mrm}[1]{\mathrm{#1}}
\newcommand{\msf}[1]{\mathsf{#1}}
\newcommand{\bs}{\backslash}
\newcommand{\TT}{\mathbb{T}}

\newcommand{\ld}{{}^L}
\newcommand{\BC}{\mathsf{BC}}
\newcommand{\FF}{\mathbb{F}}
\newcommand{\bbm}[1]{\mathbbm{#1}}
\newcommand{\bu}{\bullet}
\newcommand{\circX}{\stackrel{\circ}X}
\newcommand{\circSht}{\stackrel{\circ}\Sht}
\newcommand{\circHk}{\stackrel{\circ}\Hk}

\newcommand{\inj}{\hookrightarrow}
\newcommand{\surj}{\twoheadrightarrow}
\newcommand{\colim}{\varinjlim}

\newcommand{\OO}{\mbb{O}}

\newcommand{\eq}{B\sigma}
\newcommand{\E}{\mathbf{E}}

\newcommand{\DD}{\mathbb{D}}
\newcommand{\mon}{\G_m\mrm{-mon}}
\newcommand{\fsp}{\mathfrak{sp}}

\newcommand\cA{\mathcal{A}}
\newcommand\cB{\mathcal{B}}
\newcommand\cC{\mathcal{C}}

\newcommand\cE{\mathcal{E}}
\newcommand\cF{\mathcal{F}}

\newcommand\cH{\mathcal{H}}

\newcommand\cK{\mathcal{K}}
\newcommand\cL{\mathcal{L}}

\newcommand\cN{\mathcal{N}}
\newcommand\cO{\mathcal{O}}

\newcommand\cT{\mathcal{T}}

\newcommand\cZ{\mathcal{Z}}

\DeclareMathOperator{\GL}{GL}
\DeclareMathOperator{\SL}{SL}
\DeclareMathOperator{\Frob}{Frob}

\DeclareMathOperator{\N}{\mathbb{N}}

\DeclareMathOperator{\Hom}{Hom}

\DeclareMathOperator{\Ima}{Im\,}

\DeclareMathOperator{\Aut}{Aut}
\DeclareMathOperator{\Rep}{Rep}
\DeclareMathOperator{\Nm}{Nm}
\DeclareMathOperator{\Spec}{Spec\,}

\DeclareMathOperator{\End}{End}

\DeclareMathOperator{\ad}{ad}

\DeclareMathOperator{\Br}{Br}

\DeclareMathOperator{\Res}{Res}

\DeclareMathOperator{\cusp}{cusp}
\DeclareMathOperator{\Stab}{Stab}
\DeclareMathOperator{\Bun}{Bun}
\DeclareMathOperator{\Ext}{Ext}

\DeclareMathOperator{\Id}{Id}
\DeclareMathOperator{\Ad}{Ad}
\DeclareMathOperator{\Gr}{Gr}

\DeclareMathOperator{\pt}{pt}

\DeclareMathOperator{\Sht}{Sht}

\DeclareMathOperator{\Sat}{Sat}

\DeclareMathOperator{\Perf}{Perf}
\DeclareMathOperator{\Psm}{Psm}
\DeclareMathOperator{\Tilt}{Tilt}
\DeclareMathOperator{\Exc}{Exc}

\DeclareMathOperator{\Parity}{Parity}
\DeclareMathOperator{\Perv}{P}
\DeclareMathOperator{\ev}{ev}
\DeclareMathOperator{\Weil}{Weil}

\DeclareMathOperator{\perf}{perf}
\DeclareMathOperator{\Shv}{Shv}
\DeclareMathOperator{\Hk}{Hk}

\DeclareMathOperator{\act}{act}
\DeclareMathOperator{\alg}{alg}
\DeclareMathOperator{\geom}{geom}
\DeclareMathOperator{\Mod}{Mod}
\DeclareMathOperator{\FWeil}{FWeil}
\DeclareMathOperator{\pr}{pr}

\DeclareMathOperator{\TV}{TV}
\DeclareMathOperator{\depth}{depth}

\DeclareMathOperator{\Flat}{Flat}
\DeclareMathOperator{\cts}{cts}
\DeclareMathOperator{\Vect}{Vect}
\DeclareMathOperator{\CT}{CT}

\newcommand{\cHom}{\Cal{H}om}

\newtheorem{thm}{Theorem}[section]
\newtheorem{lemma}[thm]{Lemma}
\newtheorem{prop}[thm]{Proposition}
\newtheorem{cor}[thm]{Corollary}

\newtheorem{conj}[thm]{Conjecture}

\theoremstyle{remark}
\newtheorem{remark}[thm]{Remark} 
\newtheorem{defn}[thm]{Definition}
\newtheorem{const}[thm]{Construction}
\newtheorem{example}[thm]{Example}

\makeatletter
\def\th@remark{%
  \thm@headfont{\bfseries}%
  \normalfont 
  \thm@preskip \thm@preskip 
  \thm@postskip\thm@preskip
}
\def\imod#1{\allowbreak\mkern5mu({\operator@font mod}\,\,#1)}
\makeatother

\numberwithin{equation}{section}

\title[Smith theory and cyclic base change functoriality]{Smith theory and cyclic base change functoriality}

\author{Tony Feng}

\begin{document}

\begin{abstract}
Lafforgue and Genestier-Lafforgue have constructed the global and (semisimplified) local Langlands correspondences for arbitrary reductive groups over function fields. We establish various properties of these correspondences regarding functoriality for cyclic base change: For $\Z/p\Z$-extensions of global function fields, we prove the existence of base change for mod $p$ automorphic forms on arbitrary reductive groups. For $\Z/p\Z$-extensions of local function fields, we construct a base change homomorphism for the mod $p$ Bernstein center of any reductive group. We then use this to prove existence of local base change for mod $p$ irreducible representation along $\Z/p\Z$-extensions, and that Tate cohomology realizes base change descent, verifying a function field version of a conjecture of Treumann-Venkatesh.

The proofs are based on equivariant localization arguments for the moduli spaces of shtukas. They also draw upon new tools from modular representation theory, including parity sheaves and Smith-Treumann theory. In particular, we use these to establish a categorification of the base change homomorphism for mod $p$ spherical Hecke algebras, in a joint appendix with Gus Lonergan. 
\end{abstract}

\maketitle

\tableofcontents

\section{Introduction}

In this paper we prove several results on cyclic base change in the Langlands correspondence over function fields. To set the context for our results, let us recall some history. Global cyclic base change functoriality for reductive groups over number fields, established over many years in increasing generality by work of Saito \cite{S77}, Shintani \cite{Shi79}, Langlands \cite{L80}, Arthur-Clozel \cite{AC89}, Labesse \cite{Lab99}, Harris-Labesse \cite{HL04} and others for cuspidal automorphic representations with characteristic zero coefficients (under some technical assumptions for general groups), is one of the major triumphs of Langlands' program thus far. In addition to its initial applications towards Artin's Conjecture, it plays a crucial role in much subsequent work, such as in automorphy lifting arguments following in the tradition of Wiles, Taylor, etc. 

Our main progress in the present paper is on understanding cyclic base change in the \emph{Local} Langlands correspondence, which was constructed (in a semisimplified form) for all reductive groups over local function fields by Genestier-Lafforgue \cite{GL18}. The proof (over number fields) of global cyclic base change is founded upon the twisted trace formula, a tool which does not seem to apply to our (local and mod $p$) context, and is in any case currently unavailable over function fields due to non-trivial analytic difficulties. We introduce a new strategy, which we use to prove the first general existence results for local base change of all irreducible representations of arbitrary reductive groups over local function fields. Furthermore, we establish descent theorems for cyclic base change that were conjectured by Treumann-Venkatesh; these are new even for specific groups such as $\GL_n$ where the full Local Langlands correspondence (hence in particular the existence of local base change) is already known. These advances involve the construction of a \emph{base change homomorphism for Bernstein centers}, as has been envisaged by Haines in the case of characteristic zero coefficients. En route to the local results we establish new global results as well: we prove the first general existence theorem for cyclic base change of mod $p$ automorphic forms on arbitrary reductive groups over global \emph{function fields}, again without any trace formula arguments. A major novelty of these results is their applicability to completely general groups and representations.  

The proofs assemble a diverse selection of tools ranging from topology (particularly equivariant localization and Tate cohomology) to arithmetic geometry (of moduli stacks of shtukas) to $p$-adic groups (exploiting new constructions with Hecke algebras and Bernstein centers) to modular representation theory (using crucially the recent inventions of parity sheaves and Treumann-Smith theory).

We now proceed to give more precise descriptions of our results. 

\subsection{Local results}

Genestier-Lafforgue have constructed a semi-simplified form of the Local Langlands correspondence over function fields \cite{GL18}. More precisely, let $F_v$ be a local field of positive characteristic \emph{not} equal to $p$ and $W_v$ the Weil group of $F_v$. Let $k$ be an algebraic closure of $\F_p$.\footnote{In this paper our varieties are over fields of characteristic not equal to $p$ while our coefficients are of characteristic $p$. This is to adhere to standard notational conventions for Smith theory; unfortunately, it is at odd with standard notational conventions in arithmetic geometry.} For any reductive group $G$ over $F_v$, \cite{GL18} constructs a map 
\begin{align}\label{eq: LLC}
\left\{ \begin{array}{@{}c@{}}  \text{irreducible admissible representations} \\  \text{$\pi$ of $G(F_v)$ over $k$}\end{array} \right\}/\sim  & \longrightarrow \left\{ \begin{array}{@{}c@{}}  \text{semi-simple $L$-parameters} \\ 
\rho_{\pi} \co W_v \rightarrow \ld G(k)  \end{array} \right\}/\sim.
\end{align}
Here $\ld G$ is Langlands' $L$-group, regarded over $k$. 

Langlands' principle of functoriality predicts that given two reductive groups $H$ and $G$ over $F_v$, and a map of $L$-groups $\phi \co \ld H \rightarrow \ld G$, every $L$-packet of irreducible representations of $H(F_v)$ should admit a ``transfer'' to $G(F_v)$ compatible with $\phi$. In this paper we are concerned with a specific type of functoriality: \emph{base change functoriality}, arising from the case where $H$ is any reductive group over $F_v$, and $G = \Res_{E_v/F_v} (H_{E_v})$ for a cyclic $p$-extension $E_v/F_v$. The relevant map $ \phi_{\BC} \co \ld H \rightarrow \ld G$ is characterized by the property that it is \emph{admissible} and induces the diagonal embedding on their underlying identity connected components (i.e., the respective Langlands dual groups). We emphasize that it is crucial for our results that the degree of the extension coincides with the characteristic of our representations. In this situation, let us say that an irreducible representation $\Pi$ of $G(F_v)$ is a \emph{base change} of an irreducible representation $\pi$ of $H(F_v)$ if $\phi_{\BC}  \circ \rho_{\pi} \cong \rho_{\Pi}$. 


\begin{thm}[Existence of local base change]\label{thm: local existence}
Let $\pi$ be any irreducible representation of $H(F_v)$ over $k$. Then a base change of $\pi$ to $G(F_v)$ exists. 
\end{thm}

For $H = \GL_n$, a full Local Langlands correspondence has been established by Vign\'{e}ras \cite{Vig01}, giving a much more precise result than Theorem \ref{thm: local existence}. It seems reasonable to expect that Vign\'{e}ras' methods could (eventually) be extended to some classical groups, after the stabilization of the twisted trace formula for automorphic forms over function fields is achieved. The novelty of Theorem \ref{thm: local existence} is that it applies uniformly to all reductive groups, and all irreducible representations. A motivation for Theorem \ref{thm: local existence} is a program of the author to compute explicitly the Genestier-Lafforgue parameters of explicit supercuspidal representations, such as those arising in \cite{CO21}; this depends on the residue field being ``large enough'', so is convenient to be able to make an unramified base change.

We also prove a descent result for the above base change situation which was conjectured by Treumann-Venkatesh, and is new even for $H = \GL_n$ whenever $n > 1$. Let $\sigma$ be a generator\footnote{The choice of generator is made for convenience of notation; all constructions involving it will be manifestly independent of the choice.} of $\Gal(E_v/F_v)$; it acts on $G$ and its induced action on $G(F_v) = H(E_v)$ is the Galois action. It is expected that if the isomorphism class of a $k$-representation $\Pi$ of $G(F_v)$ is preserved by $\sigma$, then it should come from base change. For any irreducible admissible representation $\Pi$ of $G(F_v)$ whose isomorphism class is fixed by $\sigma$, there is a unique $\sigma$-action on $\Pi$ compatible with the $G(F_v)$-action (Lemma \ref{lem: 6.1}). Hence we can form the \emph{Tate cohomology} groups $T^0(\Pi)$, $T^1(\Pi)$ with respect to the $\sigma$-action (cf. \S \ref{ssec: tate cohomology}), which retain actions of $H(F_v)  = G(F_v)^{\sigma}$, and are conjecturally admissible $H(F_v)$-representations. We prove: 

\begin{thm}[Tate cohomology realizes cyclic base change]\label{thm: intro 2}
Assume $p$ is an odd good prime\footnote{Explicitly, this means that we require $p>2$ if $\wh{G}$ has simple factors of type $A,B,C$ or $D$; $p>3$ if $\wh{G}$ has simple factors of type $G_2, F_4, E_6, E_7$; and $p>5$ if $\wh{G}$ has simple factors of type $E_8$.} for $\wh{G}$. Let $\Pi$ be an irreducible representation of $G(F_v)$ whose isomorphism class is fixed by $\sigma$, and $\Pi^{(p)}: = \Pi \otimes_{k, \Frob} k$ the Frobenius twist of $\Pi$. Let $\pi$ be any irreducible admissible subquotient of $T^0(\Pi)$ or $T^1(\Pi)$ as an $H(F_v)$-representation and $\rho_{\pi} \co W_v \rightarrow \ld H(k)$ be the corresponding $L$-parameter constructed by Genestier-Lafforgue. Then $\phi_{\BC}  \circ \rho_{\pi} \cong \rho_{\Pi^{(p)}}$. 
\end{thm}

This verifies, for the Genestier-Lafforgue construction of the semi-simplified Local Langlands correspondence, a Conjecture of Treumann-Venkatesh \cite[Conjecture 6.3]{TV} that ``Tate cohomology realizes functoriality''. It had previously been proved for certain depth-zero supercuspidal representations of $\GL_n(F_v)$ by Ronchetti \cite{Ron16}, by direct calculation of the Tate cohomology and comparison to Vign\'{e}ras' work. The difficulty of the calculations, even in those special cases, made them inaccessible to generalization. By contrast, our proof applies uniformly for all groups and all representations under only a very mild condition on $p$, and is completely conceptual; in particular, it avoids any computations with specific models of representations, for example as compact inductions of Deligne-Lusztig representations.

\begin{remark}In \cite{BFHKT}, we will compute Tate cohomology for an interesting class of supercuspidal representations (of arbitrary depth) studied by Chan-Oi \cite{CO21}, which provides many examples where Theorem \ref{thm: intro 2} can be made very concrete. 
\end{remark}

We now proceed to describe our third main local result. Recall that the \emph{Bernstein center (with coefficients in $k$) of $G(F_v)$}, denoted $\mf{Z}(G)$, is the ring of endomorphisms of the identity functor on the category of smooth $G(F_v)$-representations (on $k$-vector spaces). Informally speaking, an element of $\mf{Z}(G)$ is represented by a system of compatible endomorphisms of all smooth $G(F_v)$-representations (commuting with the $G(F_v)$-action). In particular, $\mf{Z}(G)$ acts on any irreducible smooth $G(F_v)$-representation $\Pi$ through a character $\chi_{\Pi} \co \mf{Z}(G) \rightarrow k$. Furthermore, the correspondence \eqref{eq: LLC} turns out to assign isomorphic $L$-parameters to irreducible representations inducing the same character of $\mf{Z}(G)$. The ideas used to establish the preceding theorems also allow us to construct a \emph{base change homomorphism} between the Bernstein centers of $G(F_v)$ and $H(F_v)$ with the property detailed in the following Theorem. 

\begin{thm}[Base change homomorphism for Bernstein centers]\label{thm: intro bernstein center}
Assume $p$ is an odd good prime for $\wh{G}$. Then there is a homomorphism 
\[
\mf{Z}(G) \xrightarrow{\mf{Z}_{\TV}} \mf{Z}(H)
\] 
such that for each irreducible $H(F_v)$-representation $\pi$, the character $\chi_{\pi} \circ \mf{Z}_{\TV} \co \mf{Z}(G) \xrightarrow{\mf{Z}_{\TV}} \mf{Z}(H) \xrightarrow{\chi_{\pi}} k$ has the property that for any irreducible $G(F_v)$-representation $\Pi$ on which $\mf{Z}(G)$ acts through $\chi_{\pi} \circ \mf{Z}_{\TV}$, there is an isomorphism of semi-simple $L$-parameters $\rho_{\Pi} \cong \phi_{\BC} \circ \rho_{\pi}$.
\end{thm}

A base change homomorphism for Bernstein centers, with characteristic zero coefficients, has been sought by Haines \cite{Hai14}, and was constructed in some low-depth cases \cite{Hai09, Hai12} (cf. also \cite{Fe20} for the function field case). Haines also constructed a base change homomorphism for the \emph{stable} Bernstein center of general groups, which in the case of $\GL_n$ coincides with the Bernstein center; since the stable Bernstein center is defined directly in terms of Galois representations, this is rather more direct. Our Theorem \ref{thm: intro bernstein center} is somewhat different since it concerns characteristic $p$, but it is the first such construction that applies to general groups and depth. Its generality and provable connection to the Local Langlands correspondence make it rather new and compelling evidence for Haines' vision. 

\begin{remark}
The construction of the map $\mf{Z}_{\TV}$ applies equally well for local fields of characteristic $0$ having residue characteristic distinct from $p$. However, our argument for proving that it has the ``correct'' effect in terms of the Local Langlands correspondence only works for function fields. The future work \cite{F23} aims to prove analogous results with respect to Fargues-Scholze's construction \cite{FS} of the (semisimplified) local Langlands correspondence for arbitrary local fields. 
\end{remark}

\subsection{Global results} Although our most striking progress is on the local Langlands correspondence, we also obtain new results in the global Langlands correspondence. In fact, the local results mentioned above are themselves deduced from analysis of Lafforgue's machine for constructing the global Langlands correspondence. 

 Now let $G$ be a reductive group over a global function field $F$, of characteristic \emph{not} equal to $p$. Vincent Lafforgue has constructed in \cite[\S 13]{Laff18} a global ``mod $p$'' Langlands correspondence, decomposing the space of cuspidal automorphic functions $C_{\cusp}^{\infty}(G(F) \bs G(\A_F), k)$ into summands indexed by semi-simple \emph{$L$-parameters}, which are certain $\wh{G}(k)$-conjugacy classes of continuous homomorphisms $\rho \co  \Gal(F^s/F) \rightarrow \ld G(k)$. Work of Cong Xue \cite{Xue20a, Xue20, Xue21} extends Lafforgue's theory to the space of all compactly supported automorphic functions, $C_c^{\infty}(G(F) \bs G(\A_F), k)$. See \S \ref{sssec: Xue} for a more precise discussion.  Let us call an $L$-parameter $\rho$ \emph{automorphic} if it arises from Lafforgue(-Xue)'s construction. 

Let $H$ be a reductive group over $F$, and $G = \Res_{E/F} (H_E)$ for a cyclic $p$-extension $E/F$. The relevant map $\phi_{\BC} \co \ld H \rightarrow \ld G$ is the diagonal on the identity connected components. 

\begin{thm}[Existence of global base change]\label{thm: intro 1}
Assume $p$ is an odd good prime for $\wh{G}$. If $\rho \co \Gal(F^s/F) \rightarrow \ld H(k)$ is automorphic, then $\phi_{\BC} \circ \rho \co \Gal(F^s/F) \rightarrow \ld G(k)$ is automorphic. 
\end{thm}

We comment on the relation of Theorem \ref{thm: intro 1} to other base change theorems known in global contexts. To appreciate this it is important to highlight the distinction between ``weak base change'', which is determined Hecke eigensystems at almost all places of the global function, and ``strong base change'' as provided by Theorem \ref{thm: intro 1}, which concerns the entire $L$-parameter. These notions are equivalent for $H = \GL_n$, but for general groups ``strong base change'' is a strictly stronger notion. Indeed, Lafforgue's correspondence can assign different Langlands parameters to Hecke eigenfunctions with the same unramified eigensystem; in fact, it can even assign different parameters to different automorphic forms generating \emph{isomorphic} automorphic representations\footnote{This implies that there is a difference between strong base change and an intermediate notion of base change which demands compatibility with local base change at \emph{every} place (not just the unramified ones).}, with examples occurring already for $\SL_n$ when $n \geq 3$ \cite{Bla94, Lap99}. The reason for this is the failure of local conjugacy to imply global conjugacy; see \cite[\S 0.7]{Laff18} for more discussion of this phenomenon. 

\begin{remark}The distinction between weak and strong base change can be quite important in applications. For example, recent work of Sawin-Templier \cite{ST21} shows that the Ramanujan Conjecture for cuspidal automorphic forms satisfying appropriate local conditions is implied by a strong form of cyclic base change, but weak base change does not suffice for their argument. 
\end{remark}

Our proof of Theorem \ref{thm: intro 1} is inspired by work of Treumann-Venkatesh \cite{TV}, which establishes existence of ``weak base change'' for the cohomology of locally symmetric spaces. The analogue of \cite{TV} in the function field context would guarantee the existence of a ``weak base change'' for mod $p$ automorphic forms. The work of Treumann-Venkatesh is about Hecke operators, but in the function field context it is possible to go beyond Hecke operators to Lafforgue's \emph{excursion operators}, and this is necessary to obtain ``strong base change''; it is also what provides our handhold on the Local Langlands correspondence. 

Over number fields, weak base change results with \emph{characteristic zero coefficients} are known using the twisted trace formula, for all cuspidal automorphic representations of $\GL_n$ \cite{AC89} or, on more general groups, cuspidal automorphic representations satisfying certain local conditions \cite{Lab99}. Over function fields the analogous results are known for $H=\GL_n$ because the full global Langlands correspondence is already known in that case, again using the trace formula. But there are analytic difficulties in the theory of the twisted trace formula over function fields, which prevent parallel results from being known more generally. Instead, forthcoming work \cite{BFHKT} will combine Theorem \ref{thm: intro 1} with automorphy lifting theorems, generalizing those of \cite{BHKT}, in order to obtain existence of cyclic order $p$ base change for automorphic forms on split semisimple groups with characteristic $0$ coefficients, for sufficiently large $p$ and under a ``large image'' assumption (the latter is needed to make the notions of weak and strong base change coincide). 

\subsection{Remarks on the proofs}\label{ssec: intro proofs} We emphasize at the outset that our arguments make no use of the traditional tool for analyzing cyclic base change, namely the twisted trace formula (which is in any case unavailable in our situation). Any serious discussion of the proofs of our main results would require an explanation of the construction of Lafforgue's and Genestier-Lafforgue's correspondences, in addition to a number of other ideas and definitions. To prevent this introduction from becoming overly technical, we confine ourselves to vague hints here. 

The Genestier-Lafforgue correspondence is characterized by local-global compatibility, so the main input to the local results comes from an analysis of the global situation. The Global Langlands parametrization is extracted from the \emph{cohomology of moduli stacks of shtukas}. This idea goes back to Drinfeld \cite{Drin87}, who introduced it to establish the global Langlands correspondence for $\GL_2$ over function fields, later extended to $\GL_n$ by Laurent Lafforgue \cite{Laff02}. For general groups, the role that Langlands' $L$-group $\ld G$ should play presented a puzzle that was definitively resolved by Vincent Lafforgue: via the \emph{Geometric Satake equivalence}, the category $\Rep(\ld G)$ naturally indexes perverse sheaves that lives on the moduli stacks of $G$-shtukas, called $\Sht_G$. 

Summarizing roughly, the global Langlands correspondence involves two major inputs:  
\begin{enumerate}
\item A ``topological'' input, wherein the $p$-adic cohomology of spaces $\Sht_G$ supplies interesting $\Gal(F^s/F)$-representations.   
\item A representation-theoretic input, wherein $L$-parameters into $\ld G$ are extracted using that the coefficient sheaves for these cohomology groups are indexed functorially by $\Rep(\ld G)$. 
\end{enumerate}
This will be explained more in \S \ref{sec: global base change}. For now it is enough to appreciate that in order to produce a \emph{functorial transfer} from $H$ to $G$, we then need to address both of these aspects of Lafforgue's construction. More precisely, we need to: 
\begin{enumerate}
\item Show that cohomology classes on $\Sht_H$ can be ``transferred'' to cohomology classes on $\Sht_G$. 
\item Give a geometric interpretation of the restriction functor $\Rep(\ld G) \rightarrow \Rep(\ld H)$ at the level of perverse sheaves. 
\end{enumerate}
The immediate difficulty of (1) is that in general there is not so much as a non-trivial map relating $\Sht_H$ and $\Sht_G$. In the base change situation there is a natural map, but it is not even Hecke-equivariant, nor is it clear a priori that the map is not too destructive to cohomology groups. Ultimately, we solve (1) in our situation by looking at \emph{Tate cohomology} instead of cohomology, and using a form of \emph{equivariant localization} that relates the Tate cohomology of a space and its fixed points under a $\Z/p\Z$-action. Here we were inspired by work of Treumann-Venkatesh \cite{TV}, where it was shown that such equivariant localization for locally symmetric spaces realized functoriality in that context.

For (2), the obvious difficulty in general is again that we are seeking to transport sheaves between two spaces that are not connected by any visible non-trivial geometric maps. In the base change situation there is a map, but the obvious functors it induces on sheaves do not come close to having the desired effect. In some sense, the problem is a categorified and local version of the problem in the previous paragraph. Our solution to this problem passes through certain ``exotic'' localizations of categories of sheaves called \emph{Tate categories}, which can be seen as a categorification of Tate cohomology. The point is, vaguely speaking, that the desired relations of functoriality are satisfied in the relevant Tate categories. However, this does not interface well with Lafforgue's construction because localization to the Tate category does not interact well with the theory of perverse sheaves; our second main idea here is that this can be fixed by reworking the theory in terms of \emph{parity sheaves} invented by Juteau-Mautner-Williamson \cite{JMW14}. Here we were inspired by work of Leslie-Lonergan \cite{LL}, which used these tools to give a geometric interpretation of the Frobenius contraction functor in modular representation theory. (The key idea that parity sheaves play well with localization to the Tate category is also at the heart of recent work of Riche-Williamson \cite{RW}.) Ultimately, we are able to construct a ``base change functor'' that \emph{categorifies} the base change homomorphism for spherical Hecke algebras, and which is suitable for input into the global setup. The construction is completed in the joint Appendix with Gus Lonergan. 

To complete the proofs of the local results, we also need to exploit some new constructions with local Hecke algebras, in particular the base change homomorphism $\mf{Z}_{\TV}$ for Bernstein centers. A key insight in \cite{TV} is that the base change homomorphism for spherical Hecke algebras admits a more ``geometric'' description when the field extension is cyclic of order $p$ and the coefficients also have characteristic $p$. We generalize this observation to the centers of higher depth Hecke algebras, and then to the Bernstein center, by an analysis of Hecke algebras with respect to the subgroups coming from the Moy-Prasad filtration at a special vertex of the Bruhat-Tits building of $G$.

\subsection{Organization of the paper}

The outline of this paper is as follows.

 In \S \ref{sec: excursion}, we define \emph{excursion algebras} and recall their relation to Langlands parameters. We explain functoriality from the perspective of excursion algebras. 

 In \S \ref{sec: generalities}, we generalize the basic framework of \emph{sheaf-theoretic Smith theory} from \cite{Tr19, RW}, which worked for topological spaces and finite type schemes respectively, to locally finite type schemes. This is needed because our spaces of interest are not of finite type. More specifically, we introduce the notion of Tate categories, the Smith functor $\Psm$ and its properties, Tate cohomology, and explain the relation to classical equivariant localization theorems for $\Z/p\Z$-actions. 

 In \S \ref{sec: parity}, we recall the fundamentals of parity sheaves due to Juteau-Mautner-Williamson, and the analogous notion of ``Tate-parity sheaves'' due to Leslie-Lonergan. We explain how to combine these with the functor $\Psm$ to construct a \emph{base change functor} for parity objects in the Satake category. This functor plays the categorified role of the base change homomorphism for Hecke algebras.

 In \S \ref{sec: global base change}, we prove a collection of global results, including Theorem \ref{thm: intro 1}. First we recall background on moduli spaces of shtukas and Lafforgue's global Langlands correspondence in terms of actions of the excursion algebra on the cohomology of shtukas. Then we establish certain equivariant localization isomorphisms for the Tate cohomology of shtukas in the setting of $\Z/p\Z$-base change, which gives relations between excursion operators in the context of functoriality. These are used later in the local applications, and Theorem \ref{thm: intro 1} is also deduced as an application. 

 In \S \ref{sec: local base change} we prove our local results. We review the relevant aspects of the Genestier-Lafforgue correspondence. After analyzing the Brauer homomorphism for Hecke algebras with respect to subgroups arising from the Moy-Prasad filtration, we are able to construct the map $\mf{Z}_{\TV}$ from Theorem \ref{thm: intro bernstein center}, which we then establish using the global theory and local-global compatibility. Finally, we deduce Theorem \ref{thm: local existence} and Theorem \ref{thm: intro 2}.

\subsection{Acknowledgments}

We thank Jean-Fran\c{c}ois Dat, David Helm, Gus Lonergan, Simon Riche, Gordan Savin, David Treumann, Marie-France Vign\'{e}ras, Geordie Williamson, Zhiwei Yun, and Xinwen Zhu for helpful correspondence related to this work. We thank Laurent Clozel, Jesper Grodal, Tom Haines, and Michael Harris for comments on a draft. We especially thank Michael Harris for many stimulating questions, which led us to discover several new results and applications after the initial version of this paper. The paper benefited immensely from many comments and corrections by the incisive referee, who in particular suggested a big simplification of the proof of Proposition \ref{prop: inj}. During the writing of this paper, the author was supported by an NSF Postdoctoral Fellowship under grant No. 1902927, as well as the Friends of the Institute for Advanced Study. 

\subsection{Notation}\label{ssec: notation}

\begin{itemize}
\item (Coefficients) We let $k$ be an algebraic closure of $\F_p$ (considered with the discrete topology). 

In general we will consider geometric objects over fields of characteristic $\neq p$, and \'{e}tale sheaves over $p$-adically complete coefficients. 

\item ($\sigma$-actions) Throughout the paper, $\sigma$ denotes a generator of a group isomorphic to $\Z/p\Z$. When we say that a widget has a ``$\sigma$-action'', what we mean is that the widget has an action of a cyclic group of order $p$ with chosen generator $\sigma$. 

 Let $N := 1 + \sigma + \ldots + \sigma^{p-1} \in \Z[\sigma]$. We will also denote by $N$ the induced operation on any $\Z[\sigma]$-module.\footnote{This is to be contrasted with the operation $\Nm$, which will mean $\Nm(a) = a * \sigma(a) * \ldots * \sigma^{p-1}(a)$ in the context where there is a monoidal operation $*$.}

If $A$ is a ring or module for $\Z[\sigma]$, then $A^{\sigma}$ denotes the $\sigma$-invariants in $A$. 

\item (Reductive groups) For us, reductive groups are connected by definition. The Langlands dual group $\wh{G}$ is considered as a split reductive group over $k$. For our conventions on the $L$-group, see \S \ref{ssec: L-groups}. 

For any group, $\bbm{1}$ denotes the trivial representation (with the group made clear by context). 

\item (Derived categories of sheaves) If $Y$ is a locally finite type stack and $\Lambda$ is a coefficient ring in which the characteristic of $Y$ is invertible, we let $D^b_c(Y;\Lambda)$ denote the bounded constructible derived category of \'{e}tale sheaves over $\Lambda$; by this we mean complexes whose restriction to any quasi-compact open substack $U \subset Y$ lie in $D^b_c(U; \Lambda)$.

We shall also have occasional to consider larger categories of sheaves, where the constructibility condition is weakened. We let $D^b(Y; \Lambda)$ denote the bounded derived category of \'{e}tale sheaves over $\Lambda$ that are ind-constructible. In other words, it is the full subcategory of the (co-complete) category $D(Y; \Lambda)$, of ind-constructible \'{e}tale sheaves over $\Lambda$, spanned by the bounded objects. 

If $S = \{Y_\lambda\}$ is a stratification of $Y$, then we denote by $D_S^b(Y; \Lambda)$ the full subcategory of $D^b(Y; \Lambda)$ consisting of complexes constructible with respect to the stratification $S$. 

\item (Equivariant derived categories) If a (pro-)algebraic group $\Sigma$ acts on $Y$, then we denote by $D_{c,\Sigma}^b(Y; \Lambda)$ or $D^b_c(X;\Lambda)^{B\Sigma}$ the $\Sigma$-equivariant bounded derived category of constructible sheaves with coefficients in $\Lambda$. We denote by $D_{\Sigma}^b(Y; \Lambda)$ or $D^b(X;\Lambda)^{B\Sigma}$ the analogous categories with the constructibility condition replaced by ind-constructibility, as above. 

When $\Lambda = k$ we may suppress it from the notation, writing instead $D^b_c(Y ) := D^b_c(Y; k)$, etc.

\item Functors between derived categories, e.g. $f_!, f_*, f^!, f^*$, will always denote the derived functors.

\end{itemize}

 \section{Functoriality and the excursion algebra}\label{sec: excursion}

In this section we formalize the \emph{abstract excursion algebra} $\Exc(\Gamma, \ld G)$, a device used to decomposable a space into pieces indexed by Langlands parameters. This notion appears implicitly in \cite{Laff18}, but there it is the image\footnote{This image is denoted $\Cal{B}$ in \cite{Laff18}.} of the abstract excursion algebra in a certain endomorphism algebra which is emphasized. 	

Since we work with non-split groups, we first clarify in \S \ref{ssec: L-groups} our conventions regarding $L$-groups. This is a bit subtle, as one finds (at least) two natural versions of the $L$-group in the literature: the ``algebraic $L$-group'' $\ld G^{\alg}$, following Langlands, and the ``geometric $L$-group'' $\ld G^{\geom}$, derived from the Geometric Satake equivalence. The difference between them is parallel to the difference between $L$-algebraicity and $C$-algebraicity emphasized in \cite{BG14}.

We emphasize that the unadorned notation $\ld G$ denotes the algebraic $L$-group, to be consistent with \cite{Laff18}, although the geometric $L$-group is really what appears more naturally in our arguments. 

We introduce two explicit presentations for the excursion algebra in \S \ref{ssec: excursion presentation 1} and \S \ref{ssec: excursion presentation 2}. The first presentation is more natural for making the connection to Langlands parameters, which we recall in \ref{subsec: excursion gal reps}. The second presentation is more amenable to constructing actions of the excursion algebra, which makes it more convenient for our purposes, and it is the only one that will be used in the sequel. 

Finally in \S \ref{ssec: functoriality for excursion} we explain how functoriality is interpreted in terms of excursion algebras.

 \subsection{Conventions on $L$-groups and Langlands parameters}\label{ssec: L-groups}
 
For a reductive group $G$ over a field $\F$ with separable closure $\F^s$, we regard its Langlands dual group $\wh{G}$ as a split reductive group over $k$. The $L$-group is a certain semi-direct product $\ld G = \wh{G} \rtimes \Gal(\F^s/\F)$. Actually, in the case where $\F$ is a local field we shall instead work with the ``Weil form'' $\wh{G} \rtimes \Weil(\F^s/\F)$. (This is just for consistency with \cite{GL18}; because we consider representations over $k$, in our case it would make no difference to work with the Galois form.) 

\subsubsection{Algebraic $L$-group}\label{sssec: algebraic $L$-group} In fact there are at least two conventions for the definition of the $L$-group. The one which is more traditionally used in the literature is what we shall call the \emph{algebraic $L$-group}, denoted $\ld G^{\alg}$, defined as in \cite[\S 2.5]{TV}. The root datum $\Psi(G)$ of $G_{\F^s}$ determines a pinning for $\wh{G}$, which in turns gives a splitting $\mrm{Out}(\wh{G}) \rightarrow \Aut(\wh{G})$ and an identification $\Aut(\Psi(G)) \cong \mrm{Out}(\wh{G})$. The $\Gal(\F^s/\F)$-action on $\Psi(G)$ transports to an action $\act^{\alg}$ of $\Gal(\F^s/\F)$ on $\wh{G}$, and we define $\ld G^{\alg}$ to be the semidirect product 
\[
\ld G^{\alg} := \wh{G} \rtimes_{\act^{\alg}} \Gal(\F^s/\F).
\]
Since the action $\act^{\alg}$ factors through a finite quotient, we may regard $\ld G^{\alg}$ as a pro-algebraic group over $k$. 

\subsubsection{Geometric $L$-group}\label{sssec: geometric L-group} We now make a different construction of the $L$-group, using the Tannakian theory, following \cite[Appendix A]{Zhu15} and \cite[\S 5.5]{Zhu17}. We begin with the Geometric Satake equivalence, 
\[
\Perv_{L^+G_{\F^s}}(\Gr_{G,\F^s};k)  \cong \Rep_k(\wh{G}).
\]
The Galois group $\Gal(\F^s/\F)$ acts on $\Gr_{G,\F^s}$, inducing an action on the neutralized Tannakian category $(\Perv_{L^+G_{\F^s}}(\Gr_{G,\F^s};k), \underbrace{H^*(-)}_{\text{fiber functor}})$. By \cite[Lemma A.1]{Zhu15} this in turn induces an action $\act^{\geom}$ of 
$\Gal(\F^s/\F)$ on $\wh{G}_k$. We define 
\[
\ld G^{\geom} := \wh{G}_k \rtimes_{\act^{\geom}} \Gal(\F^s/\F).
\]
In the case at hand we shall see that $\act^{\geom}$ also factors through a finite quotient of $\Gal(\F^s/\F)$, so we may also regard $\ld G^{\geom}$ as a pro-algebraic group. 

\subsubsection{Relation between the two $L$-groups} The relation between these two actions is as follows. We let $\rho$ be the half sum of positive coroots of $\wh{G}$, and we denote by $\rho \co \G_m \rightarrow \wh{G}_{\ad}$ the corresponding cocharacter. With $\mrm{cyc}_p \co \Gal(\F^s/\F) \rightarrow \F_p^{\times} $ denoting the mod $p$ cyclotomic character, let $\chi$ denote the composite 
\[
\Gal(\F^s/\F) \xrightarrow{\mrm{cyc}_p} \F_p^{\times} \hookrightarrow  k^{\times} \xrightarrow{\rho} \wh{G}_{\ad}(k).
\] 
This induces a homomorphism $\Ad_{\chi} \co \Gal(\F^s/\F)  \rightarrow \Aut(\wh{G})$. 

\begin{prop}
We have $\act^{\geom} = \act^{\alg} \circ \Ad_{\chi}$.
\end{prop}

\begin{proof}
When $\wh{G}$ is over $\Q_{p}$, this is \cite[Proposition 1.6]{Zhu15}. More generally, it is established in \cite[\S VI.11]{FS} over any $p$-adic ring.
\end{proof}

Given a choice of lift $\wt{\chi} \co \Gal(\F^s/\F)  \rightarrow \wh{G}(k)$ of $\chi$, which could for example come from a square root of the mod $p$ cyclotomic character, we get an isomorphism $\ld G^{\alg} \xrightarrow{\sim} \ld G^{\geom}$ by 
\begin{equation}\label{eq: L-groups isom}
(g, \gamma) \mapsto (g \wt{\chi}(\gamma^{-1}), \gamma).
\end{equation}
By \cite[Remark 5.5.8]{Zhu17}, a square root of the cyclotomic character exists whenever $\mrm{char}(\F) >0$. (However, in general it can happen that $\ld G^{\alg}$ and $\ld G^{\geom}	$ are not isomorphic; for an example see \cite[Example 5.5.9]{Zhu17}.)

At different points we will want to consider both versions of $L$-groups. If we write $\ld G$ without a superscript, then by default we mean the algebraic $L$-group $\ld G^{\alg}$. 
 
 \subsubsection{Representation categories}\label{sssec: algebraic reps}
 
 For any Galois extension $\F'/\F$ such that $G_{\F'}$ is split, the analogous construction to \S \ref{sssec: algebraic $L$-group} gives a ``finite form'' algebraic $L$-group $\wh{G} \rtimes_{\act^{\alg}} \Gal(\F'/\F)$. We define the category of ($k$-linear) algebraic representations of $\ld G^{\alg}$ to be 	
 \[
 \Rep_k(\ld G^{\alg}) := \varinjlim_{\F'} \Rep_k(\wh{G} \rtimes_{\act^{\alg}}  \Gal(\F'/\F)).
 \]
 Let $ \Rep_k(\ld G^{\geom})  := \Rep_k(\wh{G})^{B\Gal(\F^s/\F), \geom}$ denote the category of continuously $\Gal(\F^s/\F)$-equivariant objects in $\Rep_k(\wh{G})$ with respect to the geometric action. The Geometric Satake equivalence induces by descent an equivalence 
 \begin{equation}\label{eq: arithmetic Geometric satake}
 \Perv_{L^+G}(\Gr_G;k) \cong \Rep_k(\wh{G})^{B\Gal(\F^s/\F), \geom}
 \end{equation}
where the action of $\Gal(\F^s/\F)$ on $\Rep_k(\wh{G})$ on the RHS is via $\act^{\geom}$, and on the LHS, $\Gr_G$ is considered over $\F$. By definition, on the right side we take are taking objects on which $\Gal(\F^s/\F)$ acts \emph{continuously} with its Krull topology. Since $k$ is algebraic over $\F_p$, \emph{in this case} $\Rep_k(\wh{G})^{B\Gal(\F^s/\F), \geom}$ can be identified with $\varinjlim_{\F'/\F} \Rep_k(\wh{G})^{B\Gal(\F'/\F), \geom}$ where the limit runs over finite Galois extensions $\F'/\F$ over which the geometric action factors.

An isomorphism \eqref{eq: L-groups isom} gives an embedding $\Rep_k(\ld G^{\alg}) \hookrightarrow  \Rep_k(\wh{G})^{\Gal(\F^s/\F), \geom}$, which as just remarked is an equivalence for our choice of $k$. See \cite[Proposition A.10]{Zhu15} for a description of the essential image in general.  


\subsubsection{$L$-parameters}\label{sssec: LP}

\begin{defn}Let $\Gamma$ be a topological group and $\ol{\Gamma}$ be a quotient of $\Gamma$ acting on $\wh{G}$. An \emph{$L$-parameter from $\Gamma$ to $\wh{G}(k) \rtimes \ol{\Gamma}$} is a $\wh{G}(k)$-conjugacy class of continuous homomorphisms $\rho \co \Gamma \rightarrow \wh{G}(k) \rtimes \ol{\Gamma}$, which has the property that the composite map $\Gamma \rightarrow \wh{G} \rtimes \Gamma \rightarrow \ol{\Gamma}$ is the given quotient $\Gamma \surj \ol{\Gamma}$.

Equivalently, we may view $\rho$ as an element of the continuous cohomology group $H^1_{\cts}(\Gamma, \wh{G}(k))$, where the action of $\Gamma$ on $\wh{G}(k)$ is the given one (via $\Gamma \rightarrow \ol{\Gamma}$) in the semi-direct product. 

We will consider $L$-parameters with $\wh{G}(k) \rtimes \Gamma$ being either $\ld \wh{G}^{\alg}(k)$ or $\ld \wh{G}^{\geom}(k)$, and $\Gamma$ being either $\Gal(F^s/F)$ for a global field $F$ or $\Weil(F_v^s/F_v)$ for a local field $F_v$.

Note that the algebraic $\Gamma$-action on $\wh{G}(k)$ factors through a finite quotient $\Gamma \surj \Gal(\F'/\F)$. It is clear that $L$-parameters into $\ld G^{\alg}(k)$ are in bijection (under restriction) with $L$-parameters into $\wh{G}(k) \rtimes \Gal(\F'/\F)$ for any such $\F'$; indeed, the set of all such $L$-parameters is identified with $H^1_{\mrm{cts}}(\Gamma, \wh{G}(k))$. 

 We say that a homomorphism $\rho \co \Gamma \rightarrow  \ld G^{\alg}(k)$ is \emph{semisimple}\footnote{Also called ``completely reducible'' in \cite{BHKT}.} if whenever it factors though a parabolic $\ld P^{\alg}(k) \subset \ld G^{\alg}(k)$, it also factors through a Levi $\ld M^{\alg}(k) \subset \ld P^{\alg}(k)$ (see \cite[\S 3]{Bor79} for the notion of parabolic and Levi subgroups of an $L$-group). 
\end{defn}


\subsection{Presentation of the excursion algebra}\label{ssec: excursion presentation 1} 
 Let $\Gamma$ be a group, which is either $\Gal(F^s/F)$ for a global field $F$ or $\Weil(F^s/F)$ for a local field $F$. Let $G$ be a reductive group over $F$ and $\ld G^{\alg}$ the algebraic $L$-group as defined in \S \ref{sssec: algebraic $L$-group}. 
 
 We will define the \emph{excursion algebra} $\Exc(\Gamma, \ld G^{\alg})$ to be the commutative algebra over $k$ presented by explicit generators and relations given below. (The topology on $\Gamma$ will not be relevant for the definition of $\Exc(\Gamma, \ld G^{\alg})$.) For a more conceptual perspective see \cite[\S 2]{Zhu20}, wherein the excursion algebra is denoted $k[\Cal{R}_{\Gamma, \ld G^{\alg}//\wh{G}}]$. 

\subsubsection{Generators}\label{sssec: generators 1}
We define $\Cal{O}(\ld G_k^{\alg})  := \varinjlim_{F'/F} \Cal{O}(\wh{G}_k \rtimes \Gal(F'/F))$ where the limit runs over finite extensions $F'/F$ over which the $\Gamma$-action on $\wh{G}_k$ factors. 

Generators of $\Exc(\Gamma, \ld G^{\alg})$ will be denoted $S_{I, f, (\gamma_i)_{i \in I}}$, where the indexing set $(I, f, (\gamma_i)_{i \in I})$ consists of: 
\begin{enumerate}[(i)]
\item $I$ is a finite (possibly empty) set,
\item $f \in \Cal{O}(\wh{G}_k \bs (\ld G_k^{\alg})^I / \wh{G}_k) := \Cal{O}((\ld G_k^{\alg})^I )^{\wh{G}_k  \times \wh{G}_k }$, where the quotient is for the actions of $\wh{G}_k$ by diagonal left and right translation, respectively, and
\item $\gamma_i \in \Gamma$ for each $i \in I$.
\end{enumerate}


\subsubsection{Relations}\label{sssec: relations 1} 

Next we describe the relations. (Compare \cite[\S 10]{Laff18}.)
\begin{enumerate}[(i)]
\item $S_{\emptyset, f, *} = f(1_G)$, an element of $k \subset \Exc(\Gamma, \ld G^{\alg})$. 
\item The map $f \mapsto S_{I, f, (\gamma_i)_{i \in I}}$ is a $k$-algebra homomorphism in $f$, i.e. 
\begin{align*}
S_{I, f+f', (\gamma_i)_{i \in I}} &= S_{I, f, (\gamma_i)_{i \in I}} + S_{I, f', (\gamma_i)_{i \in I}}, \\
S_{I, ff', (\gamma_i)_{i \in I}} &= S_{I, f, (\gamma_i)_{i \in I}} \cdot S_{I, f', (\gamma_i)_{i \in I}},
\end{align*}
and 
\[
S_{I, \lambda f, (\gamma_i)_{i \in I}} = \lambda S_{I, f, (\gamma_i)_{i \in I}} \text{ for all $\lambda \in k$}.
\]
\item For all maps of finite sets $\zeta \co I \rightarrow J$, all $f \in \Cal{O}(\wh{G}_k \bs (\ld G_k^{\alg})^I / \wh{G}_k)$, and all $(\gamma_j)_{j \in J} \in \Gamma^J$, we have 
\[
S_{J, f^{\zeta}, (\gamma_j)_{j \in J}} = S_{I, f, (\gamma_{\zeta(i)})_{i \in I}}
\]
where $f^{\zeta} \in \Cal{O}(\wh{G}_k \bs (\ld G_k^{\alg})^J / \wh{G}_k)$ is defined by $f^{\zeta} ((g_j)_{j \in J}) := f((g_{\zeta(i)})_{i \in I})$. 

\item For all $f \in  \Cal{O}(\wh{G}_k \bs (\ld G_k^{\alg})^I / \wh{G}_k)$ and $(\gamma_i)_{i \in I}, (\gamma_i')_{i \in I}, (\gamma_i'')_{i \in I} \in \Gamma^I$, we have 
\[
S_{I \sqcup I \sqcup I, \wt{f}, (\gamma_i)_{i \in I} \times  (\gamma_i')_{i \in I} \times (\gamma_i'')_{i \in I} } = S_{I, f, (\gamma_i (\gamma_i')^{-1} \gamma_i'')_{i \in I} },
\]
where $\wt{f} \in \Cal{O}(\wh{G}_k \bs (\ld G_k^{\alg})^{I \sqcup I \sqcup I}/ \wh{G}_k)$ is defined by 
\[
\wt{f}((g_i)_{i \in I} \times (g_i')_{i \in I} \times (g_i'')_{i \in I}) = f((g_i (g_i')^{-1} g_i'')_{i \in I}).
\]
\item If $f$ is inflated from a function on $\Gamma^I$, then $S_{I , f, (\gamma_i)_{i \in I}}$ equals the scalar $f((\gamma_i)_{i \in I})$. More generally, if $J$ is a subset of $I$ and $f$ is inflated from a function on $(\wh{G}_k \bs (\ld G_k^{\alg})^J / \wh{G}_k) \times \Gamma^{I \setminus J}$, then we have
\[
S_{I, f, (\gamma_i)_{i \in I}} = S_{J, \check{f}, (\gamma_j)_{j \in J}}
\]
where $\check{f}((g_j)_{j \in J}) := f((g_j)_{j \in J}, (\gamma_i)_{i \in I \setminus J})$. (Compare \cite[p. 164]{Laff18}.) 
\end{enumerate}

\begin{defn}
The \emph{excursion algebra} $\Exc(\Gamma, \ld G^{\alg})$ is the $k$-algebra with generators and relations specified as above. 
\end{defn}

\subsection{Constructing Galois representations}\label{subsec: excursion gal reps}

The following result of Lafforgue (generalized to modular coefficients by B\"{o}ckle-Harris-Khare-Thorne) explains how to obtain Langlands parameters from characters of $\Exc(\Gamma, \ld G^{\alg})$.

\begin{prop}[{\cite[Theorem 4.5]{BHKT}, \cite[\S 13]{Laff18}}]\label{prop: galois rep}
For any character $\nu \co \Exc(\Gamma, \ld G^{\alg}) \rightarrow k$, there is a semisimple $L$-parameter $\rho_{\nu} \co \Gamma \rightarrow \ld G^{\alg}(k)$ (for the discrete topology on $\Gamma$), unique up to conjugation by $\wh{G}(k)$, which is characterized by the following condition:

For all $n \in \N$, $f \in \Cal{O}(\wh{G}_k \bs (\ld G_k^{\alg})^{n+1}/ \wh{G}_k)$, and $(\gamma_0, \ldots, \gamma_n)\in \Gamma^{n+1}$, we have 
\begin{align}\label{eq: compatibility of galois rep} 
\nu(S_{\{0, \ldots, n\}, f, (\gamma_0, \gamma_1, \ldots, \gamma_n)}	) = f((\rho_\nu(\gamma_0\gamma_n), \rho_\nu( \gamma_1\gamma_n), \ldots, \rho_\nu(\gamma_{n-1} \gamma_n), \rho_\nu(\gamma_n))).
\end{align}
\end{prop}

\begin{remark} See also \cite[Corollary VIII.4.3]{FS} for more perspectives on, and generalizations of, this statement. 
\end{remark}

\begin{remark}
In Proposition \ref{prop: galois rep}, the datum of $\rho_{\nu}$ up to conjugation is equivalent to that of a cohomology class $[\rho_{\nu}] \in H^1(\Gamma, \wh{G}(k))$ where $\Gamma$ is given the \emph{discrete} topology.
\end{remark}

\subsection{Another presentation for the excursion algebra}\label{ssec: excursion presentation 2}
We will now describe a second presentation of $\Exc(\Gamma, \ld G^{\alg})$, following \cite[Lemma 0.31]{Laff18}, which is more useful for constructing actions of $\Exc(\Gamma, \ld G^{\alg})$ in practice. 

\subsubsection{Generators}\label{sssec: excursion generators 2} We take a set of generators indexed by tuples of data of the form $(I, W, x, \xi, (\gamma_i)_{i \in I})$, where: 
\begin{enumerate}[(i)]
\item $I$ is a finite set, 
\item $W \in \Rep_k((\ld G^{\alg})^I)$ (cf. \S \ref{sssec: algebraic reps}),
\item $x \in W$ is a vector invariant under the diagonal $\wh{G}_k$-action, 
\item $\xi \in W^*$ is a functional invariant under the diagonal $\wh{G}_k$-action,
\item $\gamma_i \in \Gamma$ for each $i$.
\end{enumerate}
The corresponding generator of $\Exc(\Gamma, \ld G^{\alg})$ will be denoted by $S_{I, W, x, \xi, (\gamma_i)_{i \in I}} \in \Exc(\Gamma, \ld G^{\alg})$.

\subsubsection{Relations}\label{sssec: relations 2} Next we describe the relations. 
\begin{enumerate}[(i)]
\item $S_{\emptyset, x, \xi, \emptyset} = \langle x, \xi \rangle$, an element of $k \subset \Exc(\Gamma, \ld G^{\alg})$. 
\item For any morphism of $(\ld G_k^{\alg})^I$-representations $u \co W \rightarrow W'$ and functional $\xi' \in (W')^*$ invariant under the diagonal $\wh{G}_k$-action, we have
\begin{equation}\label{eq: relation a-1}
S_{I, W, x, u^* (\xi'), (\gamma_i)_{i \in I} } = S_{I, W', u(x), \xi', (\gamma_i)_{i \in I}},
\end{equation}
where $u^* \co (W')^* \rightarrow W^*$ denotes the dual to $u$.

\item For two tuples $(I_1, W_1, x_1, \xi_1, (\gamma_i^1)_{i \in I_1})$ and $(I_2, W_2, x_2, \xi_2, (\gamma_i^2)_{i \in I_2})$ as in \S \ref{sssec: excursion generators 2}, we have 
\begin{equation}\label{eq: relation a-2}
S_{I_1 \sqcup I_2, W_1 \boxtimes W_2, x_1 \boxtimes x_2, \xi_1 \boxtimes \xi_2, (\gamma_i^1)_{i \in I_1} \times (\gamma_i^2)_{i \in I_2}} = S_{I_1, W_1, x_1, \xi_1, (\gamma_i^1)_{i \in I_1}} \circ S_{I_2, W_2, x_2, \xi_2, (\gamma_i^2)_{i \in I_2}}.
\end{equation}
Also,
\begin{equation}\label{eq: relation a-2b}
S_{I_1 \sqcup I_2, W_1 \oplus W_2, (x_1, x_2) ,  \xi_1 \oplus \xi_2, (\gamma_i^1)_{i \in I_1} \times (\gamma_i^2)_{i \in I_2}} = S_{I_1, W_1, x_1, \xi_1, (\gamma_i^1)_{i \in I_1}} + S_{I_2, W_2, x_2, \xi_2, (\gamma_i^2)_{i \in I_2}}.
\end{equation} 
Furthermore, the assignment $(I,W, x, \xi, (\gamma_i)_{i \in I}) \mapsto S_{I, W, x, \xi, (\gamma_i)_{i \in I}} \in \Exc(\Gamma, \ld G^{\alg})$ is $k$-linear in $x$ and $\xi$. 

\item Let $\zeta \co I \rightarrow J$ be a map of finite sets. Suppose $W \in \Rep((\ld G)^I)$, $x  \in  W^{\Delta(\wh{G})}$, $\xi \co (W^*)^{\Delta(\wh{G})}$, and $(\gamma_j)_{j \in J} \in \Gamma^J$. Letting $W^{\zeta}$ be the restriction of $W$ under the functor $\Rep((\ld G)^I) \rightarrow \Rep((\ld G)^J)$ induced by $\zeta$, we have 
\begin{equation}\label{eq: relation a-1.5}
S_{J, W^{\zeta}, x, \xi, (\gamma_j)_{j \in J}}= S_{I, W, x, \xi, (\gamma_{\zeta(i)})_{i \in I}}.
\end{equation}

\item Let $\delta_W \co \bbm{1} \rightarrow W \otimes W^*$ and $\ev_W \co W^* \otimes W \rightarrow \bbm{1}$ be the natural counit and unit. Conflating $x$ with a $\wh{G}$-invariant map $\bbm{1} \rightarrow W|_{\Delta(\wh{G})}$ and similarly for $\xi$, we have 
\begin{equation}\label{eq: relation a-3}
S_{I, W, x, \xi, (\gamma_i (\gamma'_i)^{-1} \gamma_i'')_{i \in I}} = S_{I \sqcup I \sqcup I, W \boxtimes W^* \boxtimes W, \delta_W \boxtimes x, \xi \boxtimes \ev_W, (\gamma_i)_{i \in I} \times (\gamma_i')_{i \in I} \times (\gamma_i'')_{i \in I}}.
\end{equation}

\item For $J \subset I$, if $W$ is inflated from a representation of $(\ld G^{\alg})^J \times \Gamma^{I \setminus J}$, then we have
\[
S_{I, W, x, \xi, (\gamma_i)_{i \in I}} = S_{J, W|_{(\ld G^{\alg})^J },  ((1_j)_{j \in J}, (\gamma_i)_{i \in I \setminus J}) \cdot x , \xi, (\gamma_j)_{j \in J}}.
\]
\end{enumerate}


\subsubsection{Relation between the presentations}\label{sssec: relation between presentations}

The two presentations in \S \ref{ssec: excursion presentation 1} and \S \ref{ssec: excursion presentation 2} are related as follows. The generator $S_{I, W, x, \xi, (\gamma_i)_{i \in I}}$ corresponds to $S_{I, f_{x, \xi}, (\gamma_i)_{i \in I}}$ where $f_{x, \xi}$ is the function on $(\ld G_k)^I$ given by $(g_i)_{i \in I} \mapsto \langle \xi, (g_i)_{i \in I} \cdot x \rangle$. The assumptions on $\xi$ and $x$ imply that $f_{x, \xi}$ is invariant under the left and right diagonal $\wh{G}_k$-actions. The relations in \S \ref{sssec: relations 2} imply that $S_{I, W, x, \xi, (\gamma_i)_{i\in I}}$ depends only on $f_{x, \xi}$ (and not on the choice of $x, \xi$) by \cite[Lemme 10.6]{Laff18}.


\subsection{Functoriality for excursion algebras}\label{ssec: functoriality for excursion}

A homomorphism of $L$-groups $\phi \co \ld H^{\alg} \rightarrow \ld G^{\alg}$ is \emph{admissible} if it lies over the identity map on $\Gamma$, i.e. the diagram below commutes. 
\[
\begin{tikzcd}
 \ld H^{\alg} \ar[r,"\phi"] \ar[d] &  \ld G^{\alg} \ar[d] \\
 \Gamma \ar[r, equals, "\Id"] & \Gamma
\end{tikzcd}
\]

\begin{lemma}\label{lem: functoriality homomorphism}
Let $\phi \co \ld H^{\alg} \rightarrow \ld G^{\alg}$ be an admissible homomorphism. Then there is a homomorphism $\phi^* \co \Exc(\Gamma, \ld G^{\alg}) \rightarrow \Exc(\Gamma, \ld H^{\alg})$ which in terms of the description of $k$-points of their corresponding spectra given in Proposition \ref{prop: galois rep}, sends $\rho \in H^1(\Gamma, \wh{H}(k))$ to $\phi \circ \rho \in H^1(\Gamma, \wh{G}(k))$. 
\end{lemma}

\begin{proof}
The map $\phi$ induces $\Res_{\phi} \co \Rep_k(\ld G^{\alg}) \rightarrow \Rep_k(\ld H^{\alg})$. At the level of generators, the map $\phi^*$ sends 
\[
S_{V, x, \xi, \{\gamma\}_{i \in I}} \mapsto S_{\Res_{\phi}(V), \Res_{\phi}(x), \Res_{\phi}(\xi), \{\gamma_i\}_{i \in I}}.
\]
We verify by inspection that this map sends relations to relations. To see that this indeed induces composition with $\phi$ at the level of Langlands parameters, use \eqref{eq: compatibility of galois rep}.
\end{proof}

\begin{defn}[Base change]\label{defn: phi_{BC}}
 In the base change situation, where $H$ is a reductive group over $F$ and $G = \Res_{E/F} (H_E)$, the relevant morphism of $L$-groups $\phi_{\BC} \co \ld H^{\alg} \rightarrow \ld G^{\alg}$ is defined by the formula $(h, \gamma) \mapsto (\Delta(h), \gamma)$. In fact this same formula also defines the corrresponding map of geometric $L$-groups $\phi_{\BC}^{\geom} \co \ld H^{\geom} \rightarrow	 \ld G^{\geom}$, so $\phi_{\BC}^{\geom}$ and $\phi_{\BC}$ are compatible with \eqref{eq: L-groups isom} if we use the same choice of square root of the cyclotomic character in the latter to define isomorphisms $\ld H^{\alg} \simeq \ld H^{\geom}$ and $\ld G^{\alg} \simeq \ld G^{\geom}$. We denote
\[
\phi_{\BC}^* \co \Exc(\Gamma, \ld G^{\alg}) \rightarrow \Exc(\Gamma, \ld H^{\alg})
\]
the induced map of excursion algebras. 
\end{defn}

\section{Smith theory in locally finite type}\label{sec: generalities}

Classical Smith theory concerns a type of equivariant localization that relates the mod $p$ cohomology of a topological space with the mod $p$ cohomology of its fixed points under a $\Z/p\Z$-action. Treumann proposed in \cite{Tr19} that this could be understood in terms of a ``sheaf-theoretic Smith theory'' formalism, which he developed at least in the context of complex algebraic varieties in the analytic topology. An algebraic version of this theory was built in \cite{RW} for $p$-adic \'{e}tale sheaves on finite type schemes (over fields where $p$ is invertible). We will need generalizations of this theory from finite type to \emph{locally finite type}. This is because we will want to apply the theory to the moduli spaces of shtukas, which are of not of finite type, but are locally of finite type.

Let us comment on some of the technical issues that arise in doing so. Because the cohomology of locally finite type schemes is not necessarily finite-dimensional, already the basic formalism of constructible sheaves and perfect complexes from \cite{Tr19, RW} does not apply. For example, we will have to enlarge the notion of ``Tate category'' to encompass the objects of interest. 

We do not strive for the maximum possible generality, but our theory at least encompass all examples of interest that will show up in this paper. In particular, we will use tricks to avoid discussing Smith theory for \emph{stacks}, which presents an interesting problem that could potentially refine our applications. For steps that are very similar to the case of finite type schemes as treated already in \cite{RW}, we will only sketch the proofs. 



\subsection{The Tate category}
Let $\Lambda$ be a $p$-adic coefficient ring; we will be interested in the cases where $\Lambda = k$ or $W(k)$. We will denote by $\Lambda[\sigma]$ the group ring of $\langle \sigma \rangle$ with coefficients in $\Lambda$. Our geometric objects will be over a field of characteristic $\neq p$ and we will consider $\Lambda$-adic sheaves. 

Let $Y$ be a separated, locally finite type scheme over a field. We let $\Perf(Y; \Lambda[\sigma]) \subset D^b_c(Y; \Lambda[\sigma])$ be the full subcategory consisting of complexes whose stalks at all geometric points of $Y$ are perfect over $\Lambda[\sigma]$. 

\begin{defn} We define $\Flat^b(Y; \Lambda[\sigma]) \subset D^b(Y; \Lambda[\sigma])$ to be the full subcategory consisting of bounded complexes whose stalks at all geometric points of $Y$ are represented by bounded complexes of flat (but not necessarily finite) $\Lambda[\sigma]$-modules. 
\end{defn}

The following Lemma will not be used essentially in the rest of the paper, but it may help to clarify the nature of $\Flat^b(Y; k[\sigma]) $. We thank Jesper Grodal for pointing out a gap in the original argument and also for suggesting its fix. 

\begin{lemma}\label{lem: flat =  finite tor dimension}The subcategory $\Flat^b(Y; k[\sigma]) \subset D^b(Y; k[\sigma])$ coincides with the full subcategory of objects which locally have finite tor-amplitude over $k[\sigma]$. 
\end{lemma}

\begin{proof}
For any commutative ring $A$, a complex of $A$-modules has finite tor-amplitude if and only if it is represented by a bounded complex of flat $A$-modules \cite[Tag 08G1]{stacks-project}. This shows that any object of $D^b(Y; k[\sigma])$ with locally finite tor-amplitude over $k[\sigma]$ lies in $\Flat^b(Y; k[\sigma]) $.

For the converse direction, note that a complex has tor-amplitude in $[a,b]$ if and only if all its stalks at geometric points have tor-amplitude in $[a,b]$ by \cite[Tag 0DJJ]{stacks-project}. Hence it suffices to show that if $\cK \in D^b(Y; k[\sigma])$ has finite tor-amplitude at all stalks, then its tor-amplitude is uniformly bounded. The key observation is: if $\cK \in D^b(Y; k[\sigma])$ is represented by a global complex concentrated in degrees $[a,b]$ and all its stalks at geometric points have finite tor-amplitude, then in fact all of its stalks at geometric points have tor-amplitude in $[a,b]$. Given this observation, we may conclude by using \cite[Tag 0DJJ]{stacks-project}.

Next we prove the observation. It suffices to show that if $\cK$ is a complex over $k[\sigma]$ supported in degrees $[a,b]$ which has finite tor-amplitude, then $\cK$ has tor-amplitude $[a,b]$. Since $k$ is Artinian, $k[\sigma]$ is also Artinian. For a module over an Artinian local ring, the properties of being flat and projective coincide \cite[Tag 051E]{stacks-project}. Therefore $\cK$ also has finite projective dimension \cite[Tag 0A5M]{stacks-project}. Furthermore, since $k[\sigma]$ has finitistic dimension zero, the projective amplitude of $\cK$ lies in $[a,b]$. Hence $\cK$ is represented by a complex of projective $k[\sigma]$-modules supported in degrees $[a,b]$, and therefore has tor-amplitude in $[a,b]$. 	
\end{proof}


\begin{defn} The \emph{(constructible) Tate category of $Y$} (with respect to $\Lambda$) is the Verdier quotient category $D^b_c(Y; \Lambda[\sigma]) / \Perf(Y; \Lambda[\sigma])$. 


This is the construction considered under the name ``Tate category'' in \cite{Tr19}, at least when $Y$ is a complex-analytic variety. According to \cite[Remark 4.1]{Tr19}, the category $D^b_c(Y; \Lambda[\sigma]) / \Perf(Y; \Lambda[\sigma])$ can be regarded as a derived category of perfect complexes over a certain ``$E_{\infty}$-ring spectrum'' $\Cal{T}_{\Lambda}$. So we will denote the corresponding Tate categories by $\Perf(Y; \Cal{T}_\Lambda)$. For our purposes $\Cal{T}_{\Lambda}$ can be thought of as just a notational device.  

We will require the following enlargement of the constructible Tate category. We define the \emph{(bounded ind-constructible) Tate category of $Y$} (with respect to $\Lambda$) to be the Verdier quotient category 
\[
\Shv(Y; \cT_{\Lambda}) := D^b(Y; \Lambda[\sigma]) / \Flat^b(Y; \Lambda[\sigma]).
\]

We denote the tautological projection maps from $D^b_c(Y; \Lambda[\sigma])$ to $\Perf(Y; \Cal{T}_\Lambda)$, and from $D^b(Y; \Lambda[\sigma])$ to $\Shv(Y; \cT_{\Lambda})$ by 
\[
\TT^* \co D^b(Y; \Lambda[\sigma]) \rightarrow \Shv(Y; \cT_{\Lambda}), \quad \text{ and } \quad \TT^*  \co D^b_c(Y; \Lambda[\sigma]) \rightarrow \Perf(Y; \Cal{T}_\Lambda).
\]
\end{defn}

Note that the fully faithful embedding $D^b_c(Y; \Lambda[\sigma])  \rightarrow D^b(Y; \Lambda[\sigma])$ carries $\Perf(Y; \Lambda[\sigma])$ into $\Flat^b(Y; \Lambda[\sigma])$ and so induces a functor 
\begin{equation}\label{eq: perf to shv}
\Perf(Y; \Cal{T}_\Lambda) \rightarrow \Shv(Y; \cT_{\Lambda}),
\end{equation}
which is conservative (e.g., because $\Perf(Y; \Lambda[\sigma]) \subset D^b_c(Y; \Lambda[\sigma])$ can also be characterized as the full subcategory of objects locally having finite tor-amplitude over $\Lambda[\sigma]$, according to Lemma \ref{lem: flat =  finite tor dimension}).

\begin{example}[{\cite[Proposition 4.2]{Tr19}}]\label{ex: Tate on point} 
The (bounded ind-constructible) Tate category over a point (meaning the spectrum of a separably closed field) is  $D^b(\Lambda[\sigma])/\Flat^b(\Lambda[\sigma])$. In this category the shift-by-2 functor is isomorphic to the identity functor, as one sees by considering the nullhomotopic complex 
\[
0 \rightarrow V \rightarrow V \otimes \Lambda[\sigma]  \xrightarrow{1-\sigma} V \otimes \Lambda[\sigma] \rightarrow V \rightarrow 0
\]
whose middle two terms project to $0$ in the Tate category. \end{example}

\subsection{The Smith operation}
Let $Y$ be a separated, locally finite type scheme with a $\sigma$-action that is \emph{admissible} in the sense of \cite[Expos\'{e} 5, D\'{e}finition 1.7]{SGA1}. By \cite[Remark 2.2]{RW}, this is automatic if $Y$ is exhausted by quasi-projective schemes over a field.

There is an equivariant bounded derived category $D_\sigma^b(Y; \Lambda)$. We distinguish this from the equivariant bounded constructible derived category $D_{c, \sigma}^b(Y; \Lambda) $, a full subcategory of $D_\sigma^b(Y; \Lambda)$ that will also be of interest to us. Note that since $\sigma$ acts trivially on the $\sigma$-fixed subscheme $Y^\sigma \subset Y$, we have an equivalence of derived categories 
\begin{equation}\label{eq: derived cat for trivial action}
D_{\sigma}^b(Y^\sigma;\Lambda) \cong D^b(Y^\sigma; \Lambda[\sigma]), \quad \text{and} \quad 
D_{c,\sigma}^b(Y^\sigma;\Lambda) \cong D^b_c(Y^\sigma; \Lambda[\sigma]).
\end{equation}

Then the ``Smith operation'' (cf. \cite[Definition 4.2]{Tr19}) is the functor 
\begin{equation}\label{eq: Psm for ft schemes}
\Psm  := \TT^* \circ i^* \co D_{c,\sigma}^b(Y; \Lambda) \rightarrow \Perf(Y^{\sigma}; \Cal{T}_{\Lambda})
\end{equation}
defined as the composition of $i^* \co D_{c,\sigma}^b(Y; \Lambda) \rightarrow D_{c,\sigma}^b(Y^{\sigma}; \Lambda) \stackrel{\eqref{eq: derived cat for trivial action}}\cong D^b_c(Y^{\sigma}; \Lambda[\sigma])$ with the projection $\TT^*$ to $\Perf(Y^{\sigma}; \Cal{T}_{\Lambda})$.

We extend this definition to bounded ind-constructible Tate categories in the analogous manner, defining 
\begin{equation}\label{eq: Psm for lft schemes}
\Psm  := \TT^* \circ i^* \co D_{\sigma}^b(Y; \Lambda) \rightarrow \Shv(Y^{\sigma}; \Cal{T}_{\Lambda}).
\end{equation}

\begin{remark}For $\cF \in D_{c,\sigma}^b(Y; \Lambda)$, there is potential confusion about whether ``$\Psm(\cF)$'' denotes the result of applying \eqref{eq: Psm for ft schemes} or \eqref{eq: Psm for lft schemes}. But there is a natural isomorphism between the functors 
\[
D_{c,\sigma}^b(Y; \Lambda) \xrightarrow{\Psm} \Perf(Y^{\sigma}; \Cal{T}_{\Lambda}) \xrightarrow{\eqref{eq: perf to shv}} \Shv(Y^{\sigma}; \Cal{T}_{\Lambda})
\]
and 
\[
D_{c,\sigma}^b(Y; \Lambda) \xrightarrow{\eqref{eq: perf to shv}} D_{\sigma}^b(Y; \Lambda) \xrightarrow{\Psm} \Shv(Y^{\sigma}; \Cal{T}_{\Lambda}),
\]
so the meaning is unambiguous once the ambient category is specified. When the distinction is important, we will take care to specify the ambient category. 
\end{remark}

The following properties are used to prove that our extended version of $\Psm$ retains the good behavior enjoyed by the constructible version. 

\begin{lemma}\label{lem: big Tate 1}
Retain the notation and assumptions above. Assume that the $\sigma$-action on $Y$ is free. Let $q \co Y \rightarrow Y / \sigma$ denote the quotient (which exists as a map of schemes by admissibility of the $\sigma$-action on $Y$). Then for any $\cF \in D^b(Y; \Lambda[\sigma])$, we have $q_* \cF \in \Flat^b(Y/\sigma; \Lambda[\sigma])  \subset D^b(Y/\sigma; \Lambda[\sigma])$. 
\end{lemma}

\begin{proof}
The same argument as \cite[Lemma 2.3]{RW} works here. To summarize it: for any geometric point $\ol{y} \rightarrow Y/\sigma$, and $\ol{x} \rightarrow Y$ lifting it, we have 
\[
(q_* \cF)_{\ol{y}} \cong \cF_{\ol{x}} \otimes_{\Lambda} \Lambda[\sigma],
\]
which is visibly in $\Flat^b(\ol{y}; \Lambda[\sigma])$. 
\end{proof}

\begin{lemma}\label{lem: big Tate 2}
Retain the notation and assumptions above. Let $U := Y \setminus Y^{\sigma}$ be the open complement of the $\sigma$-fixed locus of $Y$, and $j \co U \inj Y$ be its inclusion into $Y$. Then for any $\cF \in D^b(U; \Lambda[\sigma])$, and any geometric point $\ol{y}$ of $Y^{\sigma}$, the stalk $(Rj_* \cF)_{\ol{y}}$ lies in $\Flat^b(\ol{y}; \Lambda[\sigma]) \subset D^b(\ol{y}; \Lambda[\sigma])$. 
\end{lemma}

\begin{proof}
A similar argument as in \cite[Proposition 2.5]{RW} works here. Since the map $q \co Y \rightarrow Y/\sigma$ is totally ramified at $\ol{y}$, we have a $\sigma$-equivariant identification $(Rj_* \cF)_{\ol{y}} \cong (q_* Rj_* \cF)_{q(\ol{y})}$. Then by the commutativity of the diagram 
\[
\begin{tikzcd}
U \ar[r, hook, "j"] \ar[d, "q_U"] & Y \ar[d, "q"] \\
U/\sigma \ar[r, hook, "\ol{j}"] & Y/\sigma
\end{tikzcd}
\]
we have $(q_* Rj_* \cF)_{q(\ol{y})} \cong (R\ol{j}_* q_{U*} \cF)_{q(\ol{y})}$. Now Lemma \ref{lem: big Tate 1} implies that $q_{U*} \cF$ has finite tor-amplitude, and combining \cite[Tag 0F10]{stacks-project} with \cite[Expos\'{e} XVII, Th\'{e}or\`{e}me 5.2.11]{SGA4-3} implies that $R\ol{j}_*$ preserves finiteness of tor-amplitude, so their composition has locally finite tor-amplitude. 
\end{proof}

The good properties of $\Psm$ come from the following Lemma, 
which was proved for finite type schemes in \cite[Lemma 3.5]{RW} (following \cite[Theorem 4.7]{Tr19} in the topological situation). 

\begin{lemma}\label{lem: * vs !}
Retain the notation and assumptions above. Let $i \co Y^\sigma \hookrightarrow Y$. Then for any $\cF \in D^b_{\sigma}(Y; \Lambda)$, the cone of $i^! \cF \rightarrow i^* \cF$ belongs to $\Flat^b(Y^\sigma; \Lambda[\sigma])$. 
\end{lemma}

\begin{proof}Let $j \co Y \setminus Y^\sigma \inj Y$. Consider the exact triangle $i_* i^! \cF \rightarrow \cF \rightarrow j_* j^* \cF$ on $Y$. Applying $i^*$ to it yields the exact triangle in $D^b(Y^{\sigma}; \Lambda[\sigma])$: 
\[
i^! \cF \rightarrow i^* \cF \rightarrow i^* Rj_* j^* \cF.
\]
By Lemma \ref{lem: big Tate 2}, $i^* Rj_* j^* \cF \in \Flat^b(Y^{\sigma}; \Lambda[\sigma])$. 
\end{proof}


\begin{lemma}\label{lem: pushforward preserves perfect complexes}
Suppose $f \co Y \rightarrow S$ is a locally finite type and separated $\sigma$-equivariant morphism between locally finite type schemes, of bounded dimension. Then $Rf_! \co D^b(Y; \Lambda[\sigma]) \rightarrow D^b(S; \Lambda[\sigma])$ carries $\Flat^b(Y; \Lambda[\sigma])$ to $\Flat^b(S; \Lambda[\sigma])$. 
\end{lemma}

\begin{proof}
We may write $Y$ as a filtered colimit of open subschemes $Y_{\alpha}$ of finite type. Then for $\cF \in D^b(Y; \Lambda[\sigma])$, we have an identification of $Rf_! \cF$ with the colimit over $Rf_! (\cF|_{Y_{\alpha}})$. Since filtered colimits preserve flatness, we are reduced to the same statement in the finite type situation (where one can also replace ``$\Flat^b$'' by $\Perf$), which is obtained by combining \cite[Tag 0F10]{stacks-project} and \cite[Expos\'{e} XVII, Th\'{e}or\`{e}me 5.2.10]{SGA4-3}.
\end{proof}

\begin{remark}\label{remark: need big}
Note that Lemma \ref{lem: pushforward preserves perfect complexes} would \emph{not} have been true with ``$\Flat^b$'' replaced by ``$\Perf$''. This is why we need to consider ind-constructible sheaves when not in a finite type situation.
\end{remark}

\begin{cor}\label{cor: pushforward from free is flat}
Suppose $f \co Y \rightarrow S$ is a locally finite type and separated morphism between locally finite type schemes, of bounded dimension. Suppose $\sigma$ acts trivially on $S$ and freely on $Y$, and $f$ is $\sigma$-equivariant. Then $Rf_! \co D^b(Y; \Lambda[\sigma]) \rightarrow D^b(S; \Lambda[\sigma])$ lands in $\Flat^b(Y; \Lambda[\sigma])$.
\end{cor}

\begin{proof}
By the hypotheses, we may factor $f$ as the composition 
\[
Y \xrightarrow{q} Y/\sigma \xrightarrow{\ol{f}} S.
\]
Then apply Lemma \ref{lem: big Tate 1} to $q_!$ and Lemma \ref{lem: pushforward preserves perfect complexes} to $\ol{f}_!$. 
\end{proof}

\subsection{Functors on Tate categories}\label{ssec: functors on Tate cat}
Throughout this subsection we let $f \co Y \rightarrow S$ denote a $\sigma$-equivariant locally finite type morphism of locally finite type schemes with admissible $\sigma$-action, of bounded dimension. 

\subsubsection{Pullback} Since $f^* \co D^b_{\sigma}(S^{\sigma}; k) \rightarrow D^b_{\sigma}(Y^{\sigma}; k)$ preserves stalks, it preserves flat and perfect complexes, and so descends to the Tate category to induce $f^* \co \Shv(S^{\sigma}; \Cal{T}_{\Lambda}) \rightarrow \Shv(Y^{\sigma}; \Cal{T}_{\Lambda})$ and $f^* \co \Perf(S^{\sigma}; \Cal{T}_{\Lambda}) \rightarrow \Perf(Y^{\sigma}; \Cal{T}_{\Lambda})$.

\subsubsection{Proper pushforward} By Lemma \ref{lem: pushforward preserves perfect complexes}, $Rf_!\co D^b(Y^{\sigma}; \Lambda[\sigma]) \rightarrow D^b(S^{\sigma}; \Lambda[\sigma])$ descends to 
\[
Rf_! \co \Shv(Y^{\sigma}; \cT_{\Lambda}) \rightarrow \Shv(S^{\sigma}; \cT_{\Lambda}).
\]

\begin{prop}\label{prop: equivariant localization} The following diagram commutes: 
\[
\begin{tikzcd}
D_\sigma^b(Y;\Lambda) \ar[d, "\Psm"] \ar[r, "Rf_!"]  &  D^b_\sigma(S;\Lambda)   \ar[d, "\Psm"]  \\
\Shv(Y^{\sigma}; \Cal{T}_\Lambda)  \ar[r, "Rf_!"]  & \Shv(S^{\sigma}; \Cal{T}_\Lambda)
\end{tikzcd} 
\]
\end{prop}

\begin{proof}
We may as well replace $S$ by $S^{\sigma}$ and thus assume that the $\sigma$-action on $S$ is trivial. Let $\cF \in D_\sigma^b(Y;\Lambda)$. Denoting $i \co Y^{\sigma} \inj Y$ and $j$ the inclusion of the open complement, we have a distinguished triangle in $D^b_{\sigma}(Y; \Lambda)$:
\[
j_! j^* \cF  \rightarrow  \cF \rightarrow  i_* i^* \cF.
\] 
Abbreviate $f^{\sigma} := f \circ i \co Y^{\sigma} \rightarrow S$. By definition $\sigma$ acts freely on $U$, which implies that $Rf_! \circ  (j_! j^* \cF) \in \Flat^b(S; \Lambda[\sigma])$ by Corollary \ref{cor: pushforward from free is flat}. Hence the cone of $Rf_!  \cF \rightarrow Rf^{\sigma}_! (i^*  \cF)$ lies in $\Flat^b(S; \Lambda[\sigma])$, and therefore becomes $0$ in $\Shv(S; \cT_{\Lambda})$. Hence we have 
\[
\TT^*(Rf_! \cF) \cong  \TT^* ( Rf^{\sigma}_! (i^* \cF)) \cong Rf_! \Psm(\cF) \in \Shv(S; \cT_{\Lambda}),
\]
which exactly expresses the desired commutativity. 
\end{proof}

\subsection{Tate cohomology}\label{ssec: tate cohomology}
For a $\Lambda[\sigma]$-module $M$, its \emph{Tate cohomology groups} are 
\[
T^0(M) := \frac{M^{\sigma}}{N \cdot M}, \quad 
T^1(M) := \frac{\ker(N \co M \rightarrow M)}{(1-\sigma) \cdot M}, 
\]
(Recall that $N := 1 + \sigma + \ldots + \sigma^{p-1}$.) We will generalize this to complexes and then sheaves. 

\subsubsection{Tate cohomology of complexes} 

The exact sequence of $\Lambda[\sigma]$-modules
\[
0 \rightarrow \Lambda \rightarrow \Lambda[\sigma] \xrightarrow{1-\sigma} \Lambda[\sigma] \rightarrow \Lambda \rightarrow 0
\]
induces a morphism 
\begin{equation}\label{eq: lambda shift by 2}
\Lambda \rightarrow \Lambda[2] \in D^b(\Lambda[\sigma]).
\end{equation}

Given a bounded-below complex of $\Lambda[\sigma]$-modules $C^{\bu}$, we define its \emph{Tate cohomology} as 
\begin{align*}
T^i(C^{\bu}) &= \colim_{n \rightarrow \infty} \Hom_{D(\Lambda[\sigma])}(\Lambda, C^{\bu}[i+2n]) 
\end{align*}
where the transition maps are those induced by \eqref{eq: lambda shift by 2}. 

Evidently $T^i(C^{\bu})$ is 2-periodic in $i$. It is clear that this construction descends to the derived category, so we can view Tate cohomology as a collection of functors 
\[
T^i \co D^b(\Lambda[\sigma]) \rightarrow \mrm{Mod}_{/T^0(\Lambda)}.
\]
Now we specialize to the case where $\Lambda = k $. Note that by \cite[Tag 051E]{stacks-project}, a module over $k[\sigma]$ is flat if and only if it is free. Since Tate cohomology of free $k[\sigma]$-complexes vanishes (by inspection), this construction further factors through the Tate category, inducing 
\[
T^i \co \Shv(\pt; \cT_{k}) \rightarrow  \mrm{Vect}_{/k}.
\]

\begin{remark}
If $C^{\bu}$ is a bounded $k[\sigma]$-module, which will always be the case for us in practice, then we may regard $C^{\bu} \in D(\pt; k[\sigma])$ and the argument of \cite[Proposition 4.5.1]{LL} gives a natural isomorphism 
\[
T^i(C^{\bu}) \cong \Hom_{\Shv(\pt; \cT_k)}(k, \TT^* C^{\bu}[i]).
\]
\end{remark}

The following is obvious from the definition but important enough to record.

\begin{lemma}\label{lem: tate cohomology trivial action}
Suppose $C^{\bullet} \in D^b(k[\sigma])$ is inflated from $D^b(k)$, i.e. $\sigma$ acts trivially on $C^\bullet$. Then $T^* C^\bullet \cong H^*(C^{\bullet}) \otimes T^*(k)$, where $k$ is equipped with the trivial $\sigma$-action in the formation of $T^*(k)$. 
\end{lemma}

\subsubsection{Tate cohomology sheaves}\label{sssec: tate cohomology sheaves}

If $S$ has the trivial $\sigma$-action, then $\Shv(S; \Cal{T}_\Lambda)$ is defined. Given $\Cal{F} \in D^b(S; \Lambda[\sigma])$, we define \emph{Tate cohomology sheaves} 
\[
T^i \Cal{F} := \colim_{n \rightarrow \infty} \cHom_{D(S; \Lambda[\sigma])}(\Lambda, \cF[i+2n]) 
\]
where the transition maps are induced by the map $\Lambda \rightarrow \Lambda[2] \in D(S; \Lambda[\sigma])$ pulled back from \eqref{eq: lambda shift by 2}. The $T^i\cF$ are \'{e}tale sheaves of $T^0\Lambda$-modules, where $T^0\Lambda$ is the 0th Tate cohomology of $\Lambda$ viewed as a trivial $\sigma$-module.

For $\Lambda = k$, we also have the description
\[
T^i \cF \cong \cHom_{\Shv(S; \cT_k)}(k, \TT^* \cF[i])
\]
on $S$, which is an \'{e}tale sheaf of $T^0(k) = k$-modules. 

\subsubsection{Tate cohomology for a morphism} For $\Cal{F} \in D^b_{\sigma}(Y; k)$, we have $Rf_! \cF \in D^b_{\sigma}(S; k)$. If $S$ has the trivial $\sigma$-action, then we can form $T^i Rf_! \cF$, which we call the ``relative Tate cohomology of $\cF$''.

 If $S$ is the spectrum of a separably closed field equipped with the trivial $\sigma$-action, then we will abbreviate $T^i(Y;\Cal{F}) := T^i (Rf_! \cF)$, and call it the ``Tate cohomology of $Y$ with coefficients in $\Cal{F}$''. 

\begin{remark}\label{rem: tate cohomology factor over tate category}
Note that if $\sigma$ acts trivially on $Y$ and on $S$, then the construction $\cF \mapsto Rf_! \cF$ factors over $\Shv(Y; \cT_k)$ by Lemma \ref{lem: pushforward preserves perfect complexes}. In this situation we will also regard $T^i Rf_!$ as a functor on $\Shv(Y; \cT_k)$. 
\end{remark}

\subsubsection{The long exact sequence for Tate cohomology}\label{sssec: LES} 
Given a distinguished triangle $\Cal{F}' \rightarrow \Cal{F} \rightarrow \Cal{F}'' \in D^b_\sigma(Y;k)$, we have a long exact sequence 
\[
\begin{tikzcd}[row sep = tiny]
& & \ldots \ar[r] &T^{-1} Rf_!( \Cal{F}'')  \\
\ar[r] &  T^0 Rf_!(\Cal{F}') \ar[r] &  T^0 Rf_!(\Cal{F}) \ar[r] &  T^0 Rf_!( \Cal{F}'') \\
\ar[r] &  T^1 Rf_!( \Cal{F}') \ar[r] &  T^1 Rf_!( \Cal{F} ) \ar[r] &  T^1Rf_!(\Cal{F}'') \\
\ar[r] &  T^2 Rf_!(\Cal{F}' ) \ar[r]  & \ldots 
\end{tikzcd}
\]



\subsubsection{Equivariant localization}\label{ssec: general localization}

We explain that Proposition \ref{prop: equivariant localization} encompasses the classical equivariant localization theorems of ``Smith theory'', e.g., \cite[Theorem 4.2]{Qui71}. Proposition \ref{prop: equivariant localization} says that for $\Cal{F} \in D_{\sigma}^b(Y; k)$ we have 
\[
\Psm (Rf_! \Cal{F}) \cong (Rf|_{Y^{\sigma}})_!  \Psm(\Cal{F}) \in \Shv(S^\sigma;\Cal{T}_{k}).
\]
In particular, if $L$ is the spectrum of a separably closed with the trivial action of $\sigma$, we get  
\begin{equation}\label{eq: equivariant localization}
\TT^* R\Gamma_c(Y; \Cal{F}) \cong R\Gamma_c (Y^{\sigma}; \Psm(\Cal{F})).
\end{equation}

\section{Parity sheaves and the base change functor}\label{sec: parity}
We begin by indicating where this section is headed. 

The Geometric Satake equivalence $\Perv_{L^+G}(\Gr_G;k) \cong \Rep_k(\wh{G})$ provides the link between $G$ and its Langlands dual group. In the situation of functoriality, we have a map $\wh{H} \rightarrow \wh{G}$ and we would like to describe the induced restriction operation $\Rep_k(\wh{G}) \rightarrow \Rep_k(\wh{H})$ on the other side of the Geometric Satake equivalence, as a geometric operation on perverse sheaves. 

In the context of base change it is even the case that there is an embedding $\Gr_H \hookrightarrow \Gr_G$, and when seeking to describe functoriality it is natural to look to the Smith operation. (One motivation is that the papers \cite{Tr19, TV} verify that the \emph{function-theoretic} Smith operation is indeed related to functoriality for Hecke algebras.) However, the Smith operation lands in a Tate category, and in Example \ref{ex: Tate on point} we saw that in the Tate category, the shift-by-$2$ functor is isomorphic to the identity functor. This makes it seem unlikely that one can capture the notion of ``perverse sheaf'' in the Tate category. 

Juteau-Mautner-Williamson invented the theory of \emph{parity sheaves}, which have seen significance in modular representation theory. Parity sheaves are cut out in the derived category by constraints on the parity of cohomological degrees, and can therefore make sense in a context where cohomological degrees are only defined modulo $2$. The notion of \emph{Tate-parity sheaves} was introduced in \cite{LL} as an analog of parity sheaves for the Tate category, and was found to enjoy analogous properties. 

After briefly reviewing the notions of parity and Tate-parity sheaves in \S \ref{ssec: parity} and \S \ref{ssec: Tate-parity}, we will establish that the Smith operation respects the parity property, at least under certain conditions satisfied in our application of interest. Using ``coefficient lifting'' properties of parity sheaves, this will allow us to ultimately define a functor $\BC$ from parity sheaves on $\Gr_G$ to parity sheaves on $\Gr_H$, which realizes base change functoriality on the geometric side. We note that in this section, we will only need the ``constructible'' version of Smith theory for schemes, and not the generalizations developed in \S \ref{sec: generalities}.

\subsection{Parity sheaves}\label{ssec: parity}

We begin with a quick review of the theory of parity sheaves. We will take coefficients in a ring $\Lambda$ which is a complete local PID, i.e., a field or complete DVR; in our applications of interest will be either $k$ or $\OO := W(k)$. 

 Let $Y$ be a stratified variety over a separably closed field of characteristic $\neq p$, with stratification $S = \{Y_{\lambda}\}$. For the theory of parity sheaves to work, we need to assume that the (induced) stratification on $Y$ is JMW, meaning:
\begin{itemize}
\item for any two finite $\Lambda$-free local systems $\Cal{L}, \Cal{L}'$ on a stratum $Y_{\lambda}$, we have $\Ext^i(\Cal{L}, \Cal{L}')$ is free over $\Lambda$ for all $i$, and vanishes when $i$ is odd.
\end{itemize}
This holds for Kac-Moody flag varieties over separably closed fields, and in particular for affine flag varieties over separably closed fields \cite[\S 4.1]{JMW14}.

Fix a \emph{pariversity} $\dagger \co S \rightarrow \Z/2\Z$. In this paper we will always take the \emph{dimension pariversity} $\dagger(Y_\lambda) := \dim Y_{\lambda} \mod{2}$, so we will sometimes omit the pariversity from the discussion. Recall that \cite{JMW14} define \emph{even} complexes (with respect to the pariversity $\dagger$) to be those $\Cal{F} \in D_S^b(Y; \Lambda)$ such that for all $i_{\lambda} \co Y_{\lambda} \hookrightarrow Y$, for $Y_\lambda \in S$, $i_{\lambda}^* \Cal{F}$ and $i_{\lambda}^!\Cal{F}$ have cohomology sheaves which are $\Lambda$-free and supported in degrees congruent to $\dagger(Y_\lambda)$ modulo $2$, and \emph{odd} complexes analogously. They define \emph{parity complexes} to be direct sums of even and odd complexes. The full subcategory of ($S$-constructible) parity complexes (with coefficients in $\Lambda$) is denoted $\Parity_S(Y; \Lambda)$.

\begin{thm}[{\cite[Theorem 2.12]{JMW14}}]\label{thm: 2.12 JMW} Let $\Cal{F}$ be an indecomposable parity complex. Then:
\begin{itemize}
\item $\Cal{F}$ has irreducible support, which is therefore of the form $\ol{Y}_{\lambda}$ for some $Y_\lambda \in S$,
\item $i_{\lambda}^* \Cal{F}$ is a shifted $\Lambda$-free local system $\Cal{L}[m]$, and 
\item Any indecomposable parity complex supported on $\ol{Y}_{\lambda}$ and extending $\Cal{L}[m]$ is isomorphic to $\Cal{F}$. 
\end{itemize}
\end{thm}

A \emph{parity} sheaf (with respect to $\dagger$) is an indecomposable parity complex (with respect to $\dagger$) with $Y_{\lambda}$ the dense stratum in its support and extending $\Cal{L}[\dim Y_{\lambda}]$. Given $\Cal{L}[\dim Y_{\lambda}]$, it is not clear in general that a parity sheaf extending it extends exists. If it does exist, then Theorem \ref{thm: 2.12 JMW} guarantees its uniqueness, and we denote it by $\Cal{E}(\lambda, \Cal{L})$. The existence is guaranteed for $\Gr_G$ with the usual stratification by $L^+G$-orbits; $\Cal{E}(\lambda, \Cal{L})$ can moreover be promoted to a $L^+G$-equivariant complex if $p$ is not a torsion prime for $G$ \cite[Theorem 1.4]{JMW16}. If $\Cal{E}(\lambda, \Cal{L})$ exists for all $\lambda$ and $\Cal{L}$, we will say that ``all parity sheaves exist''.

 \subsection{Tate-parity sheaves}\label{ssec: Tate-parity}

As we have seen, the cohomological grading in the Tate category is only well-defined modulo $2$, so it does not seem to make sense to talk about perverse sheaves in the Tate category. However, elements of the Tate category have Tate cohomology sheaves (\S \ref{sssec: tate cohomology sheaves}), which are indexed by $\Z/2\Z$, so it \emph{could} make sense to talk about an analog of parity sheaves in the Tate category. As Leslie-Lonergan \cite{LL} observed, for this to work we must take coefficients in the \emph{integral} version of the Tate category, meaning $\Lambda = \OO = W(k)$, because then we have (say by \cite[Proposition 4.6.1]{LL})
\begin{equation}\label{eq: homotopy groups}
\Ext_{\Perf(\Cal{T}_{\OO})}^*(\TT^* (\OO),  \TT^* (\OO)) = \bigoplus_{i \in \Z} k[2i]
\end{equation}
is supported in even degrees. This is necessary for the assumption of non-vanishing odd Exts in the definition of the JMW stratification. 

For a stratification $S$ on $Y$, we define $\Perf_S(Y; \Cal{T}_{\OO}) \subset \Perf(Y; \Cal{T}_{\OO})$ to be the full subcategory generated by images of objects in $D^b_S(Y; \OO[\sigma])$. Letting $\Perf_S(Y; \OO[\sigma]) \subset \Perf(Y; \OO[\sigma])$ be the full thick subcategory of $S$-constructible objects, we have by \cite[Corollary 4.5.2]{LL} that 
\[
D^b_S(Y; \OO[\sigma]) / \Perf_S(Y; \OO[\sigma]) \xrightarrow{\sim} \Perf_S(Y; \Cal{T}_{\OO}) .
\]

\begin{defn}[{\cite[Definition 5.3.1]{LL}}] Let $\Cal{F} \in \Perf_S(Y; \Cal{T}_{\OO})$. Fix a pariversity $\dagger \co S \rightarrow \Z/2\Z$.
\begin{enumerate}
\item For $? \in \{*,!\}$, we say $\Cal{F}$ is $?$-\emph{Tate-even} (with respect to $\dagger$) if for each $Y_\lambda \in S$, we have 
\[
T^{\dagger(Y_\lambda)+1}(i_{\lambda}^?\Cal{F}) = 0.
\]
\item For $? \in \{*,!\}$, we say $\Cal{F}$ is $?$-\emph{Tate-odd} (with respect to $\dagger$) if $\Cal{F}[1]$ is ?-Tate-even. 
\item We say $\Cal{F}$ is \emph{Tate-even} (resp. \emph{Tate-odd}) if $\Cal{F}$ is both $*$-Tate even (resp. odd) and $!$-Tate even (resp. odd).
\item We say $\Cal{F}$ is \emph{Tate-parity complex} (with respect to $\dagger$) if it is isomorphic within $\Perf_S(Y; \Cal{T}_{\OO})$ to the direct sum of a Tate-even complex and a Tate-odd complex.\footnote{This is to be distinguished from the (upcoming) notion of \emph{Tate-parity sheaf}, which is more restrictive.} 
\end{enumerate}
The full subcategory of ($S$-constructible) Tate-parity complexes (with coefficients in $\Cal{T}_{\OO}$) is denoted $\Parity_S(Y; \Cal{T}_{\OO})$. If $S$ arises from the orbits of a group $G$, then the corresponding stratification is denoted with a subscript $(G)$. 
\end{defn}

Parallel to Theorem \ref{thm: 2.12 JMW}, we have the following result in this context: 

\begin{prop}[{\cite[Theorem 5.5.6]{LL}}] Let $\Cal{F}$ be an indecomposable Tate-parity complex. 
\begin{enumerate}
\item The support of $\Cal{F}$ is of the form $\ol{Y}_{\lambda}$ for a unique stratum $Y_{\lambda}$. 
\item Suppose $\Cal{G}$ and $\Cal{F}$ are two indecomposable Tate-parity complexes such that $\mrm{supp}(\Cal{G}) = \mrm{supp}(\Cal{F})$. Letting $j_{\lambda} \co Y_{\lambda} \hookrightarrow Y$ be the inclusion of the unique stratum open in this support, if $j_{\lambda}^* \Cal{G} \cong j_{\lambda}^* \Cal{F}$ then $\Cal{G} \cong \Cal{F}$. 
\end{enumerate}

\end{prop}

\begin{proof}
The same argument as in \cite[Theorem 2.12]{JMW14} works.  
\end{proof}

We define $\epsilon^* \co D_c^b(Y;\OO) \rightarrow D_c^b(Y;\OO[\sigma])$ for the inflation through the augmentation $\epsilon \co \OO[\sigma] \surj \OO$. Recall that $\TT^* \co D_c^b(Y;\OO[\sigma]) \rightarrow \Perf(Y; \Cal{T}_{\OO})$ denotes projection to the Tate category. We are interested in Tate complexes that come from the composite functor
\[
\TT^* \epsilon^* \co D^b_S(Y; \OO) \xrightarrow{\epsilon^*} D^b_S(Y; \OO[\sigma])  \xrightarrow{\TT^*} \Perf_S(Y; \Cal{T}_{\OO}).
\]

\begin{defn}
A \emph{Tate-parity sheaf} $\Cal{F} \in \Perf_S(Y; \Cal{T}_{\OO})$ is an indecomposable Tate-parity complex with the property that its restriction to the unique stratum $Y_{\lambda}$ which is dense in its support is of the form $\TT^* \epsilon^* \Cal{L}[\dim Y_{\lambda}]$ for an indecomposable $\OO$-free local system $\Cal{L}$ on $Y_{\lambda}$. If such an $\Cal{F}$ exists then it is unique, and we denote it by $\Cal{E}_{\Cal{T}}(\lambda, \Cal{L})$. 
\end{defn}

If $\Cal{E}_{\Cal{T}}(\lambda, \Cal{L})$ exists for all $Y_\lambda \in S$ and all $\Cal{L}$, we will say that ``all Tate-parity sheaves exist'' (for $Y,S$).

\subsection{Modular reduction}

We now explain that the functor $\TT^*$ has good properties that one would expect from ``base change of coefficients'' functors for categories of sheaves in classical rings. We will suppress mention of the pariversity $\dagger$.

\begin{prop}[{\cite[Proposition 5.6.3, Theorem 5.6.4]{LL}}]\label{prop: 5.16-17} \hfill 

(1) If $\Cal{F} \in D_S^b(X; \OO)$ is even/odd, then $\TT^* \epsilon^* \Cal{F} \in \Perf_S(X; \Cal{T}_{\OO})$ is Tate-even/odd. 

(2) If the parity sheaf $\Cal{E} = \Cal{E}(\lambda, \Cal{L})$ exists and satisfies $\Hom_{D^b_S(Y; \OO)} (\Cal{E}, \Cal{E}[n]) = 0$ for all $n<0$ (this holds for example if $\Cal{E}$ is perverse\footnote{In fact this is both necessary and sufficient by \cite[Lemma 6.6]{MR18}, which we thank Simon Riche for pointing out to us.}) then $\Cal{E}_{\Cal{T}}(\lambda, \Cal{L})$ exists and we have
\[
\TT^* \epsilon^* \Cal{E}(\lambda, \Cal{L}) \cong \Cal{E}_{\Cal{T}}(\lambda, \Cal{L}).
\]
\end{prop}

\begin{remark}\label{rem: modular reduction JMW} The Proposition (and its proof) are analogous to the following results of parity sheaves \cite[\S 2.5]{JMW14}. Let $\FF$ denote the base change functor 
\[
\FF = k \stackrel{L}\otimes_{\OO} (-)  \co D_S^b(Y; \OO) \rightarrow D_S^b(Y; k).
\]
The functor $\FF$ enjoys following properties. 
\begin{enumerate}
\item $\Cal{F} \in D^b_S(X; \OO)$ is a parity sheaf if and only if $\FF(\Cal{E}) \in D^b_S(X; k)$ is a parity sheaf. 
\item If $\Cal{E}(\lambda, \Cal{L}) $ exists, then $\Cal{E}(\lambda, \FF \Cal{L})$ exists and we have 
\[
\FF\Cal{E}(\lambda, \Cal{L}) \cong \Cal{E}(\lambda, \FF \Cal{L}).
\]
\end{enumerate}
\end{remark}

\begin{proof}[Proof of Proposition \ref{prop: 5.16-17}]
We reproduce the proof from \cite{LL} because it brings up certain ideas that will be needed later. The operation $\TT^* \epsilon^* $ is compatible with formation of $i_{\lambda}^*$ or $i_{\lambda}^!$. Hence to prove (1) we reduce to examining $T^i \epsilon^* \Cal{L}$ for a local system $\Cal{L}$ of free $\OO$-modules, with the trivial $\sigma$-action. This reduces to the fact that the Tate cohomology of $\OO$ is supported in even degrees, which is \eqref{eq: homotopy groups}.

For (2), we just need to check that $\TT^* \epsilon^* \Cal{E}(\lambda, \Cal{L})$ is indecomposable. Since $\Parity_S(Y; \Cal{T}_{\OO})$ is Krull-Remak-Schmidt by \cite[Proposition 5.5.2]{LL}, it suffices to check that 
the endomorphism ring of $\TT^* \epsilon^* \Cal{E}(\lambda, \Cal{L})$ is local. According to \cite[Proposition 4.6.1]{LL}, for $\Cal{F}, \Cal{G} \in D^b_S(Y; \OO)$ we have 
\begin{equation}\label{eq: morphisms between Tate-parity}
\Hom_{\Perf(Y; \Cal{T}_{\OO})}(\TT^*\epsilon^* \Cal{F}, \TT^* \epsilon^* \Cal{G}) \cong \bigoplus_{i \in \Z}  \Hom_{D^b_S(Y; k)} (\FF \Cal{F}, \FF \Cal{G}[2i]). 
\end{equation}
We apply this to $\Cal{F} = \Cal{G} = \Cal{E}(\lambda, \Cal{L})$. Since $\Cal{E}(\lambda, \Cal{L})$ is parity, \cite[(2.13)]{JMW14} applies to show that
\[
\Hom_{D^b_S(Y; k)} (\FF \Cal{E}(\lambda, \Cal{L}), \FF \Cal{E}(\lambda, \Cal{L})) = \FF \otimes \Hom_{D^b_S(Y; \OO)} (\cE(\lambda, \cL), \cE(\lambda, \cL)).
\]
By indecomposability of $\cE(\lambda, \cL)$, the ring $\Hom_{D^b_S(Y; \OO)} (\cE(\lambda, \cL), \cE(\lambda, \cL))$ is local, so $ \FF \otimes \Hom_{D^b_S(Y; \OO)} (\cE(\lambda, \cL), \cE(\lambda, \cL))$ is also local. This shows that the subalgebra on the RHS of \eqref{eq: morphisms between Tate-parity} indexed by $i=0$ is local, and the assumption implies that the summands of \eqref{eq: morphisms between Tate-parity} indexed by negative $i$ vanish. This implies the desired locality of the graded algebra \eqref{eq: morphisms between Tate-parity}.
\end{proof}

What we have seen can be summarized by the slogan: 
\begin{quote}\label{quote}
If all parity sheaves exist and have vanishing negative self-Exts, then all Tate-parity sheaves exist and $\TT^* \circ \epsilon^*$ induces a bijection between parity sheaves and Tate-parity sheaves. 
\end{quote}

\subsection{The lifting functor}\label{ssec: lifting functor} We will now define a functor lifting Tate-parity sheaves to parity sheaves. In fact the preceding slogan already tells us what to do about objects, so we just need to specify what happens on morphisms.

\begin{defn}
A \emph{normalized} (Tate-)parity complex is a direct sum of (Tate-)parity sheaves \emph{with no shifts}. Hence, under our assumptions, an indecomposable (Tate)-parity complex is normalized if and only if its restriction to the dense open stratum in its support $Y_{\lambda}$ is isomorphic to $\Cal{L}[\dim Y_{\lambda}]$ (resp. $\TT^* \epsilon^* \Cal{L}[\dim Y_{\lambda}]$) for an indecomposable local system $\cL$. We denote the full subcategory of normalized (Tate)-parity complexes by $\Parity^0_S(Y; \OO) \subset \Parity_S(Y; \OO)$ (resp. $\Parity^0_S(Y; \Cal{T}_{\OO}) \subset \Parity_S(Y; \Cal{T}_{\OO})$), and called them the \emph{categories of normalized (Tate)-parity sheaves}. 
\end{defn}

 Under the assumption that all parity sheaves exist and have vanishing negative self-Exts, Proposition \ref{prop: 5.16-17} implies that $\cE(\lambda, \cL) \mapsto \TT^* \epsilon^* \cE(\lambda, \cL) \cong \cE_{\cT}(\lambda, \cL)$ induces a bijection between normalized parity sheaves and normalized Tate-parity sheaves. We then have a ``lifting functor'' \cite[Theorem 5.6.6]{LL}
\[
L \co \Parity_S^0(Y; \Cal{T}_{\OO}) \rightarrow \Parity^0_S(Y; k)
\]
sending $\Cal{E}_{\Cal{T}}(\lambda, \Cal{L})$ to $\Cal{E}(\lambda, \Cal{L} \otimes_{\OO} k)$ on objects, and on morphisms inducing projection to the summand indexed by $i=0$ under identification \eqref{eq: morphisms between Tate-parity}. It can be thought of as an ``intermediate'' reduction between $\OO$ and $k$ in the sense that the following diagram commutes:
\begin{equation}\label{eq: lifting triangle}
\begin{tikzcd}
\Parity^0_S(Y; \OO) \ar[r, "\TT^* \epsilon^*"] \ar[dr, "\FF"'] &  \Parity_S^0(Y; \Cal{T}_{\OO})\ar[d, "L"] \\
&  \Parity^0_S(Y; k)
\end{tikzcd}
\end{equation}
 
 \subsection{Parity sheaves on the affine Grassmannian and tilting modules}\label{ssec: tilting}
 We now consider the preceding theory in the context of the affine Grassmannian $\Gr_G$ over a separably closed field $\F$, with the stratification by $L^+G$-orbits. Since this is a special case of a Kac-Moody flag variety, the stratification is JMW by \cite[\S 4.1]{JMW14}.

If $p$ is a good prime for $\wh{G}$, \cite[Corollary 1.6]{MR18} implies that all parity sheaves exist, and that all normalized parity sheaves are perverse. Therefore, the category of normalized parity sheaves corresponds under the Geometric Satake equivalence to some subcategory of $\Rep_k(\wh{G})$, and it is natural to ask what this is. The answer is given in terms of \emph{tilting modules} for $\wh{G}$ (recall that these are the objects of $\Rep_k(\wh{G})$ having both a filtration by standard objects, and a filtration by costandard objects). The tilting property is preserved by direct sum and tensor products (the latter assertion is a non-trivial theorem). Let $\Tilt_k(\wh{G}) \subset \Rep_k(\wh{G})$ denote the full subcategory of tilting modules.

\begin{thm}[{\cite[Corollary 1.6]{MR18}}]\label{thm: parity = tilting} If $p$ is good for $G$, then the Geometric Satake equivalence restricts to an equivalence\footnote{
Strictly speaking, the cited references employ the trivial pariversity instead of the dimension pariversity. Since dimensions of Schubert strata in $\Gr_G$ have constant parity on connected components, the trivial pariversity and dimension pariversity lead to the same notion of parity complexes in this case, so the only difference is in the notion of ``normalization''. We follow \cite{LL} in the use of the dimension pariversity so that perverse sheaves are $\dagger$-even.}
\[
\Parity^0_{L^+G}(\Gr_G;k) \cong \Tilt_k(\wh{G}).
\]
\end{thm}

\begin{proof}
The proof in \cite{MR18} is written for the affine Grassmannian over the complex numbers but adapts to our situation with some small modifications. First, one takes $\Gr_G$ over $\F$ instead of over $\CC$ as in \cite{MR18}. The proof of Theorem \ref{thm: parity = tilting} follows formally as in \cite[\S 6.5]{MR18} from an equivalence of categories, between the category of Iwahori-equivariant parity sheaves with coefficients in $k$ on $\Gr_G$, and the category of tilting objects in the heart of Bezrukavnikov's exotic t-structure on $\check{G} \times \G_m$-equivariant tilting objects on the Springer resolution $\wt{\cN}$. This equivalence is in turn proved by a Soergel bimodule argument. The analysis of the ``coherent side'' in \cite[\S 4,5]{MR18} is literally the same as in our situation. The analysis of the ``constructible side'' in \cite[\S 3]{MR18} applies verbatim to $\Gr_G$ over $\F$ except at one point: in \cite[Proof of Lemma 3.6, p.22]{MR18} the property that ``the map $\cO(\mf{t}^*/W \times \mf{t}^*/W) \rightarrow H^*_{L^+G}(\Gr_G; \mf{R})$ factors through $\cO(\Delta)$'' is proved using the ``loop group presentation'' of the complex affine Grassmannian; an alternate argument for this fact, that works in arbitrary characteristic, is provided in \cite[Lemma 5.2.4]{Zhu17}. 
\end{proof}

\begin{remark}
A much shorter argument for Theorem \ref{thm: parity = tilting}, but with a slightly worse bound on $p$, is given in \cite[Theorem 1.8]{JMW16}. 
\end{remark}

We need a few facts about the representation theory of tilting modules. For our arithmetic applications, the key point is that there are ``enough'' tilting modules to generate the derived category of $\Rep_k(\wh{G})$, as articulated by the statement below (which in fact applies to general highest weight categories).

\begin{prop}[{\cite[Proposition 7.17]{Riche}}]\label{prop: tiltings generate}The natural projection from the bounded homotopy category $K^b(\Tilt_k(\wh{G}))$ to $D^b(\Rep_k(\wh{G}))$ is an equivalence. 
\end{prop}

\subsection{Base change functoriality for the Satake category}\label{ssec: smith for BC}

We now consider a specific geometric situation relevant to Langlands functoriality for $p$-cyclic base change. Let $\F$ be a field of characteristic $\neq p$. We will consider reductive groups, and their affine Grassmannians, over $\F$.


	\subsubsection{The base change setup}\label{sssec: base change setup} We now specialize the situation a bit further: $H$ is any reductive group over $\F$ and $G = H^p$. We let $\sigma$ act on $G$ by cyclic rotation, sending the $i$th factor to the $(i+1)$st (mod $p$) factor. Then it is clear that the stratification on $\Gr_G$ by $L^+G$-orbits induces by restriction the stratification on $\Gr_H$ by $L^+H$-orbits.

Evidently the ``diagonal'' embedding $H \hookrightarrow G$ realizes $H$ as the fixed points of $G$ under the automorphism $\sigma$. This map $H \hookrightarrow G$ also induces a diagonal map $\Gr_H \rightarrow \Gr_G$. 

\begin{lemma}\label{lem: gr fixed points}
The diagonal map induces an isomorphism $\Gr_H \cong \Gr_G^{\sigma}$ as subfunctors of $\Gr_G$. 
\end{lemma}

\begin{proof} We have $\Gr_G \cong (\Gr_H)^p$, with $\sigma$ acting by cyclic rotation of the factors, from which the claim is clear. 
\end{proof}

Henceforth we assume that $p$ is odd and good for $\wh{G}$, so that the results of \S \ref{ssec: tilting} apply. 


We aim to give a ``geometric'' description of the corresponding functor under the Geometric Satake equivalence, $\Parity_{L^+G}(\Gr_G; k) \rightarrow \Parity_{L^+H}(\Gr_H;k)$, in terms of Smith theory. (Of course, one could give an ``ad hoc'' description using that $G = H^p$. The point is to define a functor that does not make reference to this, which will then generalize well, using descent, to the situation where $G = \Res_{\E/\F}(H)$ for a non-trivial field extension $\E/\F$.)


\begin{defn}\label{defn: Nm}
Given $\Cal{F} \in \Perv_{L^+G}(\Gr_G; \Lambda)$, we define 
\[
\Nm(\Cal{F}) := \Cal{F} \star {}^{\sigma} \Cal{F} \star \ldots \star {}^{\sigma^{p-1}} \Cal{F} \in \Perv_{L^+G \rtimes \sigma}(\Gr_G; \Lambda),
\]
equipped with the $\sigma$-equivariant structure coming from the commutativity constraint for $(\Perv_{L^+G}(\Gr_G;\Lambda), \star)$: 
\begin{equation}\label{eq: equivariant structure}
{}^{\sigma}\Nm(\Cal{F}) = {}^{\sigma} \Cal{F} \star \ldots \star  {}^{\sigma^{p-1}} \Cal{F}  \star \Cal{F} \xrightarrow{\sim} \Cal{F} \star {}^{\sigma} \Cal{F} \star \ldots \star {}^{\sigma^{p-1}} \Cal{F} = \Nm(\Cal{F}).
\end{equation}
Using the realization functor $\Perv_{L^+G \rtimes \sigma}(\Gr_G; \Lambda) \rightarrow D_{L^+G \rtimes \sigma}(\Gr_G; \Lambda)$, we view $\Nm(\Cal{F}) \in  D_{L^+G \rtimes \sigma}(\Gr_G;\Lambda)$ (so that we may apply the Smith functor, for example). Equipping a general object of $D_{L^+G}(\Gr_G;\Lambda)$ with a $\sigma$-equivariant structure is much more involved than just specifying isomorphisms \eqref{eq: equivariant structure} (satisfying cocycle conditions), so we emphasize that we construct $\Nm(\Cal{F})$ first as a $\sigma$-equivariant perverse sheaf, and then apply the realization functor to get a $\sigma$-equivariant object of $D_{L^+G}(\Gr_G;\Lambda)$. 
\end{defn}

\begin{remark}
In our applications we will assume that $p$ is large enough so that all parity sheaves are perverse. The properties of being $L^+G$-constructible and $L^+G$-equivariant are equivalent for perverse sheaves on $\Gr_G$. Therefore, we will not need to worry about any extra complications coming from the equivariance. 
\end{remark}

\begin{lemma}\label{lem: LL crit 1}
Let $i \co \Gr_H \cong \Gr_G^{\sigma} \hookrightarrow \Gr_G$. For $\Cal{F} \in \Perv_{L^+G}(\Gr_G; \OO)$, regard $\Nm(\Cal{F}) \in \Perv_{L^+G \rtimes \sigma}^b(\Gr_G; \OO)$ as in Definition \ref{defn: Nm} above. Suppose that all the cohomology sheaves of $\Cal{F}$ have $\OO$-free stalks and costalks. 
\begin{enumerate}
\item[(i)] The stalks of the cohomology sheaves of $i^* \Nm(\Cal{F})$ have an $\OO[\sigma]$-stable filtration with associated graded a direct sum of either trivial or free $\OO[\sigma]$-modules. 
\item[(ii)] The costalks of the cohomology sheaves of $i^! \Nm(\Cal{F})$ have an $\OO[\sigma]$-stable filtration with associated graded a direct sum of either trivial or free $\OO[\sigma]$-modules. 
\end{enumerate}
\end{lemma}

\begin{proof} If $\cF$ has a finite $L^+G$-equivariant filtration whose associated graded satisfies the hypotheses of the Lemma, then the statement of Lemma for $\cF$ can be checked on the associated graded. Since $\Gr_G^{\ul{\lambda}}$ is a product of homogeneous spaces for (a finite type quotient of) $L^+ H$, there is a finite $L^+G$-equivariant filtration of $\cF$ with associated graded of being a direct sum of sheaves of the form $\Cal{F}_1 \boxtimes \ldots \boxtimes \Cal{F}_p$ where each $\cF_i \in \Perv_{L^+H}(\Gr_H;  \OO)$. Hence we reduce to the case where $\cF$ is itself of this form. Then 
\[
\Nm(\Cal{F}) \approx (\Cal{F}_1 \star \Cal{F}_2 \star \ldots \star \Cal{F}_p) \boxtimes (\Cal{F}_2 \star \ldots \star \Cal{F}_p \star \Cal{F}_1) \boxtimes \ldots \boxtimes (\Cal{F}_p \star \Cal{F}_1 \star \ldots \star \Cal{F}_{p-1}),
\]
with $\sigma$ acting by rotating the tensor factors, and the $\sigma$-equivariant structure coming from the commutativity constraint. 

Write $\Cal{F}' := \Cal{F}_1 \star \Cal{F}_2 \star \ldots \star \Cal{F}_p \in \Perv_{L^+H}(\Gr_H; \OO)$. Since $i$ may be identified with the diagonal embedding $\Gr_H \hookrightarrow \Gr_H^p$, we have $i^*(\Nm \Cal{F})  \approx (\Cal{F}')^{\otimes p}$, with $\sigma$-equivariant structure given by cyclic rotation of the tensor factors. In particular, the stalk of $i^* (\Nm \Cal{F})$ at $x \in \Gr_H$ is the tensor-induction of the stalk of $\Cal{F}'_x$ from $\OO$ to $\OO[\sigma]$. Hence it suffices to prove that any cohomology sheaf of such a tensor induction has an $\OO[\sigma]$-equivariant filtration by either trivial or free $\OO[\sigma]$-modules. This is verified by explicit inspection: choosing a basis for $\Cal{H}^j(\Cal{F}'_x)$, the induced basis of $\Cal{H}^j(\Cal{F}'_x)^{\otimes p}$ is grouped into either trivial or free orbits under the $\sigma$-action. 

The argument for (ii) is completely analogous (alternatively, we could deduce it simply by applying Verdier duality to (i)). 
\end{proof}

\subsubsection{Smith theory for parity sheaves}

We return momentarily to the general setup for Smith theory: $Y$ is a variety over $\F$ with an admissible $\sigma$-action and $Z = Y^{\sigma}$. We assume that $\F$ is separably closed and the stratification $S$ on $Y$ satisfies the JMW condition.

\begin{prop}[Variant of {\cite[Theorem 6.1.1]{LL}}]\label{prop: Tate-parity}Assume that each stratum $Y_\lambda$ is smooth. Suppose $\Cal{E} \in D^b_{S, \sigma}(Y; \OO)$ is a parity complex satisfying the condition:
\begin{itemize}\label{eq: stalk condition}
\item[(*)] all $*$ and $!$-stalks of cohomology sheaves of $\Cal{E}$ at fixed points $y \in Y$ have an $\OO[\sigma]$-stable filtration with associated graded being a direct sum of trivial or free $\OO[\sigma]$-modules. 
\end{itemize}
Then $\Psm(\Cal{E}) \in \Perf_S(Z	; \Cal{T}_{\OO})$ is Tate-parity with respect to the induced stratification $Z_{\lambda} = Y_{\lambda} \cap Z$ and the induced pariversity $\dagger_{Z}(\lambda) := \dagger_Y(\lambda)$. 
\end{prop}

\begin{proof} This theorem is closely related to \cite[Theorem 6.1.1]{LL}, but \emph{loc. cit.} imposes the stronger condition that the $\sigma$-action on all stalks is trivial. This is satisfied in their application (to the loop-rotation action), but not in ours, so we need to re-do the argument in the requisite generality. 

Let $Z  = Y^{\sigma}$ and take the induced stratification on $Z$. Let $i \co Z \rightarrow Y$, $i_{\lambda}^Y \co Y_{\lambda} \hookrightarrow Y$, $i_{\lambda}^Z \co Z_{\lambda} \hookrightarrow Z$, $i^{\lambda} \co Z_{\lambda} \hookrightarrow Y_{\lambda}$. Without loss of generality suppose $\Cal{E}$ is an even complex on $Y$. We are given that $(i_{\lambda}^Y)^{?} \Cal{E}$ has $\OO$-free cohomology sheaves supported in degrees congruent to $\dagger_Y(\lambda)$ mod $2$, where $? \in \{*, !\}$; we want to show that $(i_{\lambda}^Z)^? \Psm(\Cal{E})$ has Tate-cohomology sheaves supported in degrees congruent to $\dagger_Z(\lambda)$ mod $2$. Unraveling the definitions, we have 
\begin{align*}
(i_{\lambda}^Z)^* \Psm(\Cal{E}) & =(i_{\lambda}^Z)^* \TT^* i^* \Cal{E}  \cong \TT^* (i_{\lambda}^Z)^* i^* \Cal{E}  \cong  \TT^* (i^{\lambda})^* (i_{\lambda}^Y)^*\Cal{E}.
\end{align*}
Similarly, using Lemma \ref{lem: * vs !} we have 
\begin{equation}\label{eq: !-restrict}
(i_{\lambda}^Z)^! \Psm(\Cal{E}) \cong \TT^* (i^{\lambda})^! (i_{\lambda}^Y)^! \Cal{E}.
\end{equation}
By hypothesis, $(i_{\lambda}^Y)^* \Cal{E}$ has its cohomology sheaves supported in degrees congruent to $\dagger_Y(\lambda) \pmod{2}$. So the stalks of $(i^{\lambda})^* (i_{\lambda}^Y)^* \Cal{E}$ are supported in degrees congruent to $\dagger_Y(\lambda) \pmod{2}$, and we must verify that their Tate cohomology groups are also supported in degrees of a single parity. 

By assumption (*), all the stalks have an $\OO[\sigma]$-stable filtration with associated graded being a direct sum of trivial or free $\OO[\sigma]$-modules. For trivial $\OO[\sigma]$-modules the odd Tate cohomology groups vanish by \eqref{eq: homotopy groups}, while for free $\OO[\sigma]$-modules all the Tate cohomology groups vanish. Hence all odd Tate cohomology groups vanish by the long exact sequence for Tate cohomology (\S \ref{sssec: LES}). This shows that the Tate cohomology sheaves of $(i^{\lambda})^* (i_{\lambda}^Y)^* \Cal{E}$ are supported in degrees congruent to $\dagger_Y(Y_\lambda) \pmod{2}$.

To show that $(i^{\lambda})^! (i_{\lambda}^Y)^! \Cal{E}$ also has Tate cohomology sheaves supported in degrees congruent to $\dagger_Y(\lambda) \pmod{2}$, we make a similar analogous argument using \eqref{eq: !-restrict} instead. This shows that $ \TT^* (i^{\lambda})^* (i_{\lambda}^Y)^! \Cal{E}$ lies in degrees congruent to $\dagger_Y(Y_\lambda) \pmod{2}$, and then we conclude by observing that $(i^{\lambda})^* (i_{\lambda}^Y)^! \Cal{E}$ differs from $(i^{\lambda})^! (i_{\lambda}^Y)^! \Cal{E}$ by an even shift (and twist) by the Gysin isomorphism, which applies because the strata are assumed to be smooth (noting that the smoothness of $Z_\lambda$ follows from the smoothness of $Y_\lambda$ by \cite[Proposition A.8.11]{CGP}, so $Z_\lambda \inj Y_\lambda$ is a regular embedding). 
\end{proof}

For an $\OO$-linear abelian category $\msf{C}$, with all Hom-spaces being free $\OO$-modules, we abbreviate
\[
\msf{C} \otimes_{\OO} k := \msf{C} \otimes_{\OO-\Mod} (k-\Mod).
\]

\begin{lemma}\label{lem: category base changes}
Suppose that all the strata $Y_{\lambda}$ are simply connected and all parity sheaves $\Cal{E}(\lambda, \Cal{L})$ exist, for all $Y_\lambda \in S$. Then we have that 
\[
\mrm{Parity}^0_{S, \sigma}(Y; \OO)  \otimes_{\OO} k \xrightarrow{\sim} \mrm{Parity}^0_{S, \sigma}(Y; k).
\] 
\end{lemma}

\begin{proof}
To see that the functor is well-defined, we note:
\begin{itemize}
\item The Hom-spaces of $\mrm{Parity}^0_{S, \sigma}(Y; \OO) $ are all free $\OO$-modules by \cite[Remark 2.7]{JMW14}, so that the domain is well-defined. 
\item The functor lands in parity sheaves since the modular reduction of a $\OO$-parity sheaf is a $k$-parity sheaf by Remark \ref{rem: modular reduction JMW}.
\end{itemize}
It is essentially surjective because every $k$-parity sheaf lifts to a $\OO$-parity sheaf under our assumption that all parity sheaves exist and all strata are simply connected (which implies that all $k$-local systems on strata lift to $\OO$, since they are trivial). The fact that the functor is fully faithful follows from \cite[(2.39)]{JMW14}. 
\end{proof}

\subsubsection{The base change functor}\label{sssec: base change functor}

We return now to the base change setup of \S \ref{sssec: base change setup}, with $\F$ separably closed. Let $\Cal{F} \in \Parity^0_{L^+G}(\Gr_G; \OO)$. Then $\Cal{F} \in \Perv_{L^+G}(\Gr_G; \OO)$ is perverse since $p$ is good for $\wh{G}$ (this is a part of Theorem \ref{thm: parity = tilting}), and $\Nm(\Cal{F}) \in \Parity^0_{L^+G \rtimes \sigma}(\Gr_G; \OO)$ is a parity sheaf by \cite[Theorem 1.5]{JMW16}. Furthermore, the $\sigma$-equivariant structure on $\Nm(\Cal{F})$ satisfies the assumption (*) of Proposition \ref{prop: Tate-parity} by Lemma \ref{lem: LL crit 1}. 

The Schubert cells of $\Gr_G$ are indexed by tuples $\lambda := (\lambda_1, \ldots, \lambda_p) \in X_*(G)^+$, with each $\lambda_i \in X_*(H)^+$, and we have 
\[
\begin{cases}
\Gr_G^{\lambda} \cap \Gr_H = \Gr_H^{\lambda_1} & \lambda = (\lambda_1, \ldots, \lambda_1), \\
\Gr_G^{\lambda} \cap \Gr_H  = \emptyset & \text{otherwise.}
\end{cases}
\]
We claim that as long as $p>2$, the induced pariversity coincides with the dimension pariversity on $\Gr_H$, i.e., for $\lambda = (\lambda_1, \ldots, \lambda_p) \in X_*(G)^+$, we have 
\[
\dim \Gr_H^{\lambda_1}  \equiv \dim \Gr_G^{\lambda} \pmod{2}.
\]
This will imply that:
\begin{enumerate}
\item We may apply Proposition \ref{prop: Tate-parity} to deduce that $\Psm(\Nm(\Cal{F})) \in \Parity_{(L^+H)}(\Gr_H; \Cal{T}_{\OO})$ is Tate-parity with respect to the dimension pariversity on $\Gr_H$. 
\item $\Psm(\Nm(\Cal{F})) \in \Parity^0_{(L^+H)}(\Gr_H; \Cal{T}_{\OO})$, i.e., is normalized.
\end{enumerate}
To prove the claim, we may focus on the case where $\lambda_1 = \ldots = \lambda_p$ or else the statement is vacuous. By \cite[Proposition 2.1.5]{Zhu17} we have $\dim \Gr_G^{\lambda} =  \langle 2\rho_G, \lambda \rangle$. So we just have to verify that $\langle 2 \rho_G , (\lambda_1, \ldots, \lambda_1) \rangle \equiv \langle 2\rho_H, \lambda_1 \rangle \pmod{2}$. Indeed, $\rho_G = (\rho_H, \ldots, \rho_H)$, so $\langle 2 \rho_G , (\lambda_1, \ldots, \lambda_1) \rangle = p \langle 2 \rho_H, \lambda_1 \rangle$, and $p$ is odd.

Thanks points (1) and (2) above, we can apply the lifting functor $L$ to $\Psm(\Nm(\Cal{F}))$. By Lemma \ref{lem: category base changes}, the composite functor $L \circ \Psm \circ \Nm$ factors uniquely through a functor $\Parity^0_{L^+G}(\Gr_G;k) \rightarrow \mrm{Parity}^0_{L^+H}(\Gr_H; k)$.

\begin{const}[Frobenius twist of categories]\label{const: frob twist C} Let $\Frob$ be the absolute Frobenius of $k$. Given a $k$-linear category $\msf{C}$, there is another $k$-linear category $\msf{C}^{(p)} := \msf{C} \otimes_{k, \Frob} k$. Concretely, it is equivalent to the category which has the same objects as $\msf{C}$, and morphisms 
\[
\Hom_{\msf{C}^{(p)}}(x,y) = \Hom_{\msf{C}}(x,y)^{(p)} := \Hom_{\msf{C}}(x,y) \otimes_{k, \Frob} k.
\]
The tautological map $\Hom_{\msf{C}}(x,y) \rightarrow \Hom_{\msf{C}}(x,y)^{(p)}$ is $\Frob$-semilinear over $k$, and induces an equivalence $\Frob_{\msf{C}} \co \msf{C} \xrightarrow{\sim} \msf{C}^{(p)}$ which is $\Frob$-semilinear. The functor $\Frob_{\msf{C}} \co \msf{C} \rightarrow \msf{C}^{(p)}$ is characterized by the universal property that any $\Frob$-semilinear functor $F \co \msf{C} \rightarrow \msf{D}$ (meaning a functor between $k$-linear categories that is $\Frob$-semilinear over $k$ on morphisms) factors uniquely through a $k$-linear functor $\msf{C}^{(p)} \rightarrow \msf{D}$. 
\[
\begin{tikzcd}
\msf{C} \ar[d, "\Frob_{\msf{C}}"'] \ar[dr, "F"] & \\
\msf{C}^{(p)} \ar[r, dashed] & \msf{D}
\end{tikzcd}
\]

Now, given a presentation 
\begin{equation}\label{eq: Frob structure on C}
F_0 \co \msf{C} \cong \msf{C}_0 \otimes_{\F_p} k := \msf{C}_{0} \otimes_{\Vect_{/\F_p}} \Vect_{/k}
\end{equation}
 for some $\F_p$-linear category $\msf{C}_0$, then there is another, \emph{$k$-linear} equivalence $ \msf{C} \xrightarrow{\sim} \msf{C}^{(p)}$, which with reference to \eqref{eq: Frob structure on C} is the tensor product of $\Id_{\msf{C}_0}$ with the $k$-linear equivalence $\Vect_{/k}^{(p)} \cong \Vect_{/k}$ induced by the $k$-linear isomorphism $k \otimes_{\Frob, k} k \cong k$. Therefore \eqref{eq: Frob structure on C} induces a $\Frob$-semilinear equivalence $\Frob_{F_0} \co \msf{C} \xrightarrow{\sim} \msf{C}$.

 \end{const}

\begin{defn}\label{defn: BC on Gr}  We define
\[
\msf{BC}^{(p)} \co \Parity^0_{L^+G}(\Gr_G; k) \rightarrow \Parity_{L^+H}^0(\Gr_H; k) 
\]
to be the functor unique filling in the commutative diagram
\begin{equation}\label{eq: dashed arrow}
\begin{tikzcd}[column sep = large]
\mrm{Parity}^0_{L^+G}(\Gr_G; \OO)  \ar[r, "\Psm \circ \Nm"] \ar[d, "\FF"] & \mrm{Parity}^0_{(L^+H)}(\Gr_H; \Cal{T}_{\OO}) \ar[d, "L"] \\
\Parity^0_{L^+G}(\Gr_G; k)  \ar[r, dashed, "\BC^{(p)}"] & \mrm{Parity}^0_{L^+H}(\Gr_H; k).
\end{tikzcd}
\end{equation}

One more step is required to define what we call the \emph{base change functor} $\BC$. Note that $\BC^{(p)}$ is $\Frob$-semilinear over $k$; we wish to linearize it. It is evident from the definitions that the equivalence $D(\Gr_G; k) \cong D(\Gr_G; \F_p) \otimes_{\F_p} k$ induces 
\[
F_0 \co \Parity^0_{L^+G}(\Gr_G; k) \cong \Parity^0_{L^+G}(\Gr_G; \F_p) \otimes_{\F_p} k.
\]
Let $\Frob_p := \Frob_{F_0} \co \Parity^0_{L^+G}(\Gr_G; k) \xrightarrow{\sim} \Parity^0_{L^+G}(\Gr_G; k)$ be the $k$-semilinear equivalence induced by $F_0$, as explained in Construction \ref{const: frob twist C}. We define
\[
\BC :=   \BC^{(p)} \circ \Frob_p^{-1}: \Parity^0_{L^+G}(\Gr_G; k)  \rightarrow \mrm{Parity}^0_{L^+H}(\Gr_H; k).
\]
\end{defn}

\begin{remark} The construction of $\BC$ was motivated by a similar functor ``$LL$'' appearing in \cite[\S 6.2]{LL}, which gives a partial geometric description of the Frobenius contraction functor on $\wh{G}$. Another motivation was the ``normalized Brauer homomorphism'' of \cite[\S 4.3]{TV}, which our construction categorifies. \end{remark}

\begin{thm}\label{thm: BC on Gr}
Let $\Res_{\BC} \co \Rep_k(\wh{G}) \rightarrow	 \Rep_k(\wh{H})$ be restriction along the diagonal embedding. We also denote by $\Res_{\BC}$ the same functor restricted to the subcategories of tilting modules.\footnote{Note that it is not obvious that $\Res_{\BC}$ preserves the tilting property, but this follows from the non-trivial theorem (building on work of many authors -- see the discussion around \cite[Theorem 1.2]{JMW16}) that tensor products of tilting modules are tilting.} The following diagram commutes: 
\[
\begin{tikzcd}
 \Parity^0_{L^+G}(\Gr_G; k) \ar[d, "\sim"]  \ar[r, "\BC"] &  \Parity^0_{L^+H}(\Gr_H; k) \ar[d, "\sim"] \\
\mrm{Tilt}_k(\wh{G}) \ar[r, "\Res_{\BC}"] & \mrm{Tilt}_k(\wh{H})
\end{tikzcd}
\]
\end{thm}

The proof is given in Appendix \ref{sec: categorical base change}. 

The triangulated structure on $D_{L^+H}(\Gr_H; k[\sigma])$ equips $\Perf_{L^+H}(\Gr_G; \cT_k)$ with the notion of cone (namely, the image of a cone in $D_{L^+H}(\Gr_H; k[\sigma])$). We say that a sequence $\cA \rightarrow \cB \rightarrow \cC$ in $D_{L^+H}(\Gr_H; \cT_k)$ is \emph{exact} if the induced map $\mrm{Cone}(\cA \rightarrow \cB) \rightarrow \cC $ in $\Perf_{L^+H}(\Gr_H; \cT_k)$ is an isomorphism. 

\begin{lemma}\label{lem: exactness}
The composite functor
\[
\Rep_k(\wh{G}) \xrightarrow{\Frob_p^{-1}} \Rep_k(\wh{G}) \xrightarrow{\Sat} \Perv_{L^+G}(\Gr_G; k) \xrightarrow{\Nm} \Perv_{L^+G \rtimes \sigma}(\Gr_G; k) \xrightarrow{\Psm} \Perf_{(L^+H)}(\Gr_H; \cT_k) 
\]
is exact, i.e., sends exact sequences to exact sequences in the above sense. 
\end{lemma}

The proof requires some notions from Appendix \ref{sec: categorical base change}, and will be postponed to \S \ref{ssec: proof of exactness}.

\subsubsection{Equivariantization and Galois descent}\label{sssec: equivariantization} 
Assuming Theorem \ref{thm: BC on Gr}, let us give a few variants related to descent to a ground field which is not separably closed. Suppose $H$ base changed from some subfield $\F_0 \subset \F$, and $G = \Res_{\E_0/\F_0} (H_{\E_0})$ for some Galois extension $\E_0/\F_0$ with Galois group $\Z/p\Z$. Then $G_{\F} \approx (H_{\F})^p$ and $\Aut(\F/\F_0)$ acts on $H_{\F},G_{\F}$ and therefore also on $\Gr_{H, {\F}}, \Gr_{G, {\F}}$. 

\begin{lemma}[Galois equivariance]\label{lem: bc equivariance} In the situation above, the functor 
\[
\BC \co \Parity^0_{L^+G}(\Gr_{G, \F}; k) \rightarrow \Parity^0_{L^+H}(\Gr_{H, \F}; k) 
\]
is equivariant with respect to the action of $\Aut(\F/\F_0)$. 
\end{lemma}

\begin{proof}
The constituent functors $\Nm$, $i^*$, $\TT^*$, and $L$ are all $\Aut(\F/\F_0)$-equivariant, as is $\Frob_p^{-1}$. It remains only to see that the dashed arrow in \eqref{eq: dashed arrow} is $\Aut(\F/\F_0)$-equivariant. This follows because $L \circ \Psm \circ \Nm$ and $\FF$ both have this property, and $\FF$ is essentially surjective and full.  
\end{proof}

We refer to \cite{DGNO} for the theory of ``equivariantization'' and ``de-equivariantization'' of categories. Given a group $\Gamma$ acting on categories $\msf{C}, \msf{D}$ and a $\Gamma$-equivariant functor $F \co \msf{C} \rightarrow \msf{D}$, the \emph{$\Gamma$-equivariantization} of $F$ is the functor $F^{B \Gamma} \co \msf{C}^{B \Gamma} \rightarrow \msf{D}^{B \Gamma}$. If $\msf{C}$ and $\msf{D}$ are derived categories of sheaves and $F$ is induced by geometric operations that are $\Gamma$-equivariant, then the equivariantization construction exists for equivariant derived categories. (We make this remark because if the $\Gamma$-equivariantization of a derived category is not the same as the $\Gamma$-equivariant derived category, and it is the latter that we want to consider.) 

Thanks to Lemma \ref{lem: bc equivariance}, the equivariantization of $\BC$ induces 
\[
\BC^{B\Aut(\F/\F_0)} \co \Parity^0_{L^+G}(\Gr_{G, \F}; k)^{B\Aut(\F/\F_0)}  \rightarrow \Parity^0_{L^+H}(\Gr_{H, \F}; k)^{B\Aut(\F/\F_0)}.
\]
We \emph{define} $\Parity^0_{L^+G}(\Gr_{G, \F_0};k) := \Parity^0_{L^+G}(\Gr_{G, \F}; k)^{B\Aut(\F/\F_0)} $ and similarly for $H$ (note that in \S \ref{ssec: parity} parity sheaves were only defined for varieties over separably closed fields, since the axioms of a JMW stratification would not otherwise be satisfied). We define $\Tilt_k(\ld G)$ to be the subcategory of $\Rep_k(\ld G)$ consisting of representations whose restriction to $\wh{G}$ is tilting, and $\mrm{Tilt}_k(\ld G^{\geom})$ to be the full subcategory  of $\Rep_k(\ld G^{\geom})$ consisting of representations whose restriction to $\wh{G}$ is tilting; then $\mrm{Tilt}_k(\ld G^{\geom}) \cong \Tilt_k(\wh{G})^{\Aut(\F/\F_0), \geom}$ and similarly for $\ld H$. 
 
  Then applying $\Aut(\F/\F_0)$-equivariantization to Theorem \ref{thm: BC on Gr} yields: 

\begin{cor}
The following diagram is commutative. 
\[
\begin{tikzcd}[column sep = huge]
 \Parity^0_{L^+G}(\Gr_{G, \F_0}; k) \ar[d, "\sim"]  \ar[r, "\BC^{B\Aut(\F/\F_0)}"] &  \Parity^0_{L^+H}(\Gr_{H, \F_0}; k) \ar[d, "\sim"] \\
\mrm{Tilt}_k(\ld G^{\geom}) \ar[r, "\Res_{\BC}"] & \mrm{Tilt}_k(\ld H^{\geom})
\end{tikzcd}
\]
\end{cor}

\subsubsection{} Let $H/\F_0$ and $G/\F_0$ be as before. The following compatibility statement will be needed later.

\begin{lemma}\label{lem: integral cube} The cube 
\begin{equation}\label{eq: integral cube}
\begin{tikzcd}[column sep = tiny]
& & \Parity^0_{L^+G \rtimes \sigma }(\Gr_G; k) \ar[rrr] \ar[ddd, hook] & & &  \ar[ddd, "\TT^* \epsilon^*"] \Parity^0_{L^+H}(\Gr_H; k) \\
\\ 
\Parity^0_{L^+G\rtimes \sigma }(\Gr_G; \OO)  \ar[rrr, "\Psm"] \ar[ddd, hook] \ar[uurr, "\FF"]  & & & \Parity^0_{(L^+H)}(\Gr_H; \cT_{\OO})  \ar[ddd, hook] \ar[uurr, "L"] \\
& & D^b_{c, L^+G \rtimes \sigma}(\Gr_G;k) \ar[rrr, "\Psm"] &  & & \Perf_{(L^+H)}(\Gr_H; \cT_k)   \\ 
\\
D^b_{c, L^+G \rtimes \sigma }(\Gr_G;\OO)  \ar[rrr, "\Psm"] \ar[uurr, "\FF"]  & & & \Perf_{(L^+ H)}(\Gr_H; \cT_{\OO}) \ar[uurr, "\FF"] 
\end{tikzcd}
\end{equation}
commutes, where $L$ is lifting functor $L \co \Parity^0_{L^+H}(\Gr_H; \cT_{\OO})  \rightarrow \Parity^0_{L^+H}(\Gr_H; k)$ from \S \ref{ssec: lifting functor} and the unlabeled arrow is defined as the one that makes the top face commute. (It exists by the universal property of the categorical tensor product.) 
\end{lemma}

\begin{proof}
The left face commutes by definition of $\FF$. It is obvious from the definition that the front face commutes. The bottom face commutes by compatibility of $\Psm$ with tensoring coefficients. The top face commutes by definition. 

To see that the right face commutes, consider the diagram 
\begin{equation}\label{eq: right square factor}
\begin{tikzcd}
\Parity^0_{L^+G}(\Gr_G; \OO) \ar[r, "\TT^* \epsilon^*"] \ar[d] & \Parity^0_{(L^+H)}(\Gr_H; \cT_{\OO}) \ar[d] \ar[r, "L"] & \Parity^0_{(L^+H)}(\Gr_H; k) \ar[d, "\TT^* \epsilon^*"] \\
D_{c,L^+H}^b(\Gr_H; \OO) \ar[r, "\TT^* \epsilon^*"] & \Perf_{(L^+H)} (\Gr_H; \cT_{\OO}) \ar[r, "\FF"] & \Perf_{(L^+H)}(\Gr_H; \cT_k)  
\end{tikzcd}
\end{equation}
The right square of \eqref{eq: right square factor} is the right square of the cube \eqref{eq: integral cube}, which we want to show commutes. Since the middle column is a Verdier quotient of the left column, it suffices to show that the outer square of \eqref{eq: right square factor} commutes. Next note that the composite of the upper horizontal arrows in \eqref{eq: right square factor} is the modular reduction functor $\FF$ by definition \eqref{eq: lifting triangle}, so we can factor the outer square of \eqref{eq: right square factor} as the outer quadrilateral in the diagram below. 
\[
\begin{tikzcd}
\Parity^0_{L^+G}(\Gr_G; \OO) \ar[r, "\FF"] \ar[d, hook] &  \Parity^0_{L^+G}(\Gr_H; k) \ar[d, hook] \ar[ddr]  \\
D_{c, L^+H}^b(\Gr_H; \OO) \ar[r, "\FF"] \ar[drr] & D_{c, L^+H}^b(\Gr_H; k) \ar[dr, "\TT^* \epsilon^*"] \\
& &  \Perf_{(L^+H)}^b(\Gr_H; \cT_k) 
\end{tikzcd}
\]
Obviously the commutativity of the outer quadilateral follows from the commutativity of the inner rectangle, which is then immediate from the definition of $\FF$. 

It remains to show that the back face commutes. Consider juxtaposing the top and back faces of cube \eqref{eq: integral cube} to get: 
\begin{equation}\label{eq: back factor 1}
\begin{tikzcd}
\Parity^0_{L^+G}(\Gr_G; \OO) \ar[d, "\FF"]   \ar[r, "\Psm"] & \Parity_{(L^+H)}^0(\Gr_H; \cT_{\OO}) \ar[d, "L"]\\
\Parity^0_{L^+G}(\Gr_G; k) \ar[d] \ar[r, dashed] & \Parity^0_{L^+H}(\Gr_H; k) \ar[d, "\TT^* \epsilon^*"] \\
D_{c, L^+G}^b(\Gr_G;k) \ar[r, "\Psm"] & \Perf_{(L^+H)}(\Gr_H; \cT_k)
\end{tikzcd}
\end{equation}
We want to show that the lower square commutes. The composite vertical arrows on left and right columns are both the modular reduction functor $\FF$ from $\OO$-coefficients to $k$ coefficients, so the outer square commutes. The dashed arrow is defined as the $k$-linearization of $L \circ \Psm$, noting that $\Parity^0_{L^+G}(\Gr_G; k) \cong \Parity^0_{L^+G}(\Gr_G; \OO) \otimes_{\OO} k$. Therefore the right-then-down (resp. down-then-right) composite functor in the lower square is the $k$-linearization of the right-then-down (resp. down-then-right) composite functor in the upper square, so the commutativity of the lower square follows from that of the upper square, completing the proof.

\end{proof}

\begin{remark}
Let us try to make some vague remarks about the utility of Lemma \ref{lem: integral cube}. The unlabeled arrow in the top face is a priori somewhat mysterious, but the Lemma says that after projecting to the Tate category, it has a simple description in terms of $\Psm$. Later, we will take Tate cohomology with coefficients indexed by the type of parity sheaves constructed in this section. Note that Tate cohomology factors through the projection of these sheaves to the Tate category. Therefore, the computation of Tate cohomology is not so sensitive to the subtleties in the constructions of this section; the purpose of this section has more to do with the indexing of coefficient sheaves, in terms of the discussion of \S \ref{ssec: intro proofs}. 
\end{remark}

 \section{On global base change}\label{sec: global base change}
 
 In this section we will apply the preceding theory to \emph{moduli stacks of shtukas}, in the context of Lafforgue's construction of the global Langlands parametrization for function fields. In particular, we will prove Theorem \ref{thm: intro 1}, among other results. 
 
We briefly review the relevant parts of Lafforgue's construction in \S \ref{ssec: shtukas} and \S \ref{ssec: VLaff decomposition}. Then in \S \ref{subsec: fusion functor construction}, where we use a variant of Lafforgue's ideas to construct and analyze an action of the excursion action on Tate cohomology of moduli spaces of shtukas. In the situation of base change, equivariant localization mediates between the Tate cohomology of shtukas for $G$ and for $H$, allowing us to relate certain excursion operators for the two groups. This is then used in \S \ref{ssec: base change} to establish the existence of base change for mod $p$ automorphic forms; this relation will also be the crucial input for our local results in the next section. 
  
\subsection{Moduli of shtukas}\label{ssec: shtukas}

We will use the theory of moduli stacks of shtukas, due to Drinfeld and generalized by Varshavsky. Here we very briefly recall the relevant definitions in order to set notation. More comprehensive references include \cite{Var04} and \cite{Laff18}. 

\subsubsection{Shtukas} Fix a smooth projective curve $X$ over a finite field $\F_{\ell}$ of characteristic $\neq p$. Let $G$ be a smooth algebraic group scheme over $X$. We assume that $G$ is generically reductive, and let $\circX \inj X$ be the locus where $G$ is reductive. For each finite set $I$, the stack $\Sht_{G,I}$ has the following functor of points on $\F_{\ell}$-schemes $S$: 
\[
\Sht_{G,I} \co S \mapsto \left\{   \begin{array}{@{}c@{}} 
(x_i)_{i \in I}  \in X^I(S)  \\
\Cal{E}  = \text{\'{e}tale $G$-torsor over $X \times S$}  \\
\varphi \colon \Cal{E}|_{X \times S - \bigcup_{i \in I} \Gamma_{x_i}} \xrightarrow{\sim} \ft \Cal{E}|_{X \times S - \bigcup_{i \in I} \Gamma_{x_i}} 
 \end{array} \right\},
\]
where $\tau$ is the Frobenius $\Frob_{\ell}$ on the $S$ factor in $X \times  S$, and  ${}^{\tau} \Cal{E}$ is the pullback of $\Cal{E}$ under the map $1 \times \tau \co X \times S \rightarrow X \times S$. 

Geometrically, $\Sht_{G,I}$ has a Schubert stratification whose strata are Deligne-Mumford stacks locally of finite type. We regard it as an ind-(locally finite type) Deligne-Mumford stack.

There is a map 
\[
\pi_I \co \Sht_{G,I} \rightarrow X^I
\]
projecting a tuple $(\{x_i\}_{i \in I}, \Cal{E}, \varphi_i)$ to $\{x_i\}_{i \in I}$. Let $\circSht_{G,I} := \Sht_{G,I} \times_{X^I} (\circX)^I$. 

\subsubsection{Hecke stack} The Hecke stack $\Hk_{G,I}$ classifies
\[
\Hk_{G,I} \co S \mapsto \left\{   \begin{array}{@{}c@{}} 
(x_i)_{i \in I}  \in X^I(S)  \ \\
\Cal{E} , \Cal{E}'  = \text{\'{e}tale $G$-torsors over $X \times S$}  \\
\varphi \colon \Cal{E}|_{X \times S - \bigcup \Gamma_{x_i}} \xrightarrow{\sim}  \Cal{E}'|_{X \times S - \bigcup \Gamma_{x_i}} 
 \end{array} \right\}.
 \]
Recall that $G \rightarrow X$ is reductive over $\circX$. Let $\circHk_{G,I} := \Hk_{G,I} \times_{X^I} (\circX)^I$. The Geometric Satake equivalence provides a functor $\Rep_k((\ld G)^I) \rightarrow D(\circHk_{G,I};k)$, which we normalize as in \cite[Theorem 0.9]{Laff18}. 

\subsubsection{Satake sheaves}\label{sssec: satake functor} There is a map $\Sht_{G,I} \rightarrow \Hk_{G,I}$ sending $(\{x_i\}_{i \in I}, \Cal{E}, \varphi)$ to $(\{x_i\}_{i \in I}, \Cal{E}, \ft \Cal{E}, \varphi)$.  Composing with the $*$-pullback through $\Sht_{G,I} \rightarrow \Hk_{G,I}$ induces a functor
\[
\Sat^{\geom} \co \Rep_k(\wh{G}^I)^{B\Gal(F^s/F), \geom}\rightarrow D^b(\circSht_{G,I};k).
\]
Finally, we may identify $\Rep_k((\ld G^{\alg})^I) \xrightarrow{\sim} \Rep_k(\wh{G}^I)^{B\Gal(F^s/F), \geom}$ as in \S \ref{sssec: algebraic reps}, giving a functor (cf. \cite[Theorem 0.11]{Laff18})
\[
\Sat \co \Rep_k((\ld G^{\alg})^I) \rightarrow D^b(\circSht_{G,I};k).
\]
The Schubert stratification is defined by the support of the sheaves in the image of $\Sat$, with the closure relations corresponding to the Bruhat order. (In particular, $\Sat$ lands in the derived category of sheaves constructible with respect to the Schubert stratification on $\circSht_{G,I}$.)

\subsubsection{Level structures} For $D \subset X$ a closed finite subscheme, there are level covers $\Sht_{G,D,I} \rightarrow \Sht_{G,I}|_{(X \setminus D)^I}$ which parametrize the additional datum of a trivialization of $\Cal{E}$ over $S \times D$ compatible with $\tau$ and $\varphi$. Note that by definition, the ``legs'' $\{x_i\}_{i \in I} \in (X \setminus D)(S)^I$ avoid $D$.

\subsubsection{Iterated shtukas}

Let $I_1, \ldots, I_r$ be a partition of $I$. We define $\Sht_{G,D,I}^{(I_1, \ldots, I_r)}$ (sometimes called a moduli stack of \emph{iterated shtukas}) to be the stack 
 \[
\Sht_{G,D,I}^{(I_1, \ldots, I_r)} \co S \mapsto \left\{   \begin{array}{@{}c@{}} 
(x_i)_{i \in I} \in (X - D)^I(S)  \\
\Cal{E}_0,\ldots, \Cal{E}_r  = \text{\'{e}tale $G$-torsors over $X \times S$}  \\
\varphi_j \colon \Cal{E}_{j-1}|_{X \times S - \bigcup_{i \in I_j} \Gamma_{x_i}} \xrightarrow{\sim}  \Cal{E}_{j}|_{X \times S - \bigcup_{i \in I_j} \Gamma_{x_i}} \quad j = 1, \ldots, r \\ 
\varphi \co \Cal{E}_r \xrightarrow{\sim} \ft \Cal{E}_0  \\
\upsilon = \text{level structure over $D \times S$}
 \end{array} \right\}.
\]
Here by ``level structure'' we mean a trivialization of the restriction of each $\cE_i$ over $D \times S$, compatible with the $\varphi_j$ and $\varphi$. There is a map $\nu \co \Sht_{G,D,I}^{(I_1, \ldots, I_r)} \rightarrow \Sht_{G,D,I}$. A key property of this morphism is that it is \emph{stratified small} (with respect to the Schubert stratification), which is a consequence of the same property of the convolution morphism for Beilinson-Drinfeld Grassmannians.

Define 
\[
\circSht_{G,D,I} = \Sht_{G,D,I} \times_{(X \setminus D)^I} (\circX \setminus D)^I
\]
and 
\[
\circSht_{G,D,I}^{(I_1, \ldots, I_r)} = \Sht_{G,D,I}^{(I_1, \ldots, I_r)} \times \times_{(X \setminus D)^I} (\circX \setminus D)^I.
\]

 \subsubsection{Partial Frobenius}\label{sssec: partial frob}

There is a partial Frobenius $F_{I_1} \co  \Sht_{G,D,I}^{(I_1,  I_2, \ldots, I_r)}  \rightarrow \Sht_{G,D,I}^{(I_2, \ldots, I_r, I_1)} $ sending 
\begin{align*}
x_i  & \mapsto \begin{cases} \ft x_i & i \in I_1\\ x_i &\text{otherwise}\end{cases} \\
(\Cal{E}_0, \ldots, \Cal{E}_r)  & \mapsto (\Cal{E}_1, \ldots, \Cal{E}_r, \ft \Cal{E}_0) \\
(\varphi_1, \ldots, \varphi_r) &\mapsto (\varphi_2, \ldots, \varphi_r, \ft \varphi_1).
\end{align*} 
It lies over the partial Frobenius $\Frob_{I_1}$ on $X^I$ (applying $\Frob_\ell$ to the coordinates indexed by $i \in I_1$), so that the diagram below is commutative (and cartesian up to radiciel maps):
\begin{equation}\label{diag: partial Frob}
\begin{tikzcd}
\Sht_{G,D,I}^{(I_1, \ldots, I_r)}  \ar[r, "F_{I_1}"] \ar[d, "\pi_I"] & \Sht_{G,D,I}^{(I_2, \ldots, I_r, I_1)} \ar[d, "\pi_I"]  \\
(X-D)^I \ar[r, "{\Frob_{I_1}}"] & (X-D)^I
\end{tikzcd}
\end{equation}
 
\subsubsection{Base change setup}

We now consider the following ``base change setup''. Let $F$ be the function field of $X$ and $H_F$ a reductive group over $F$. We choose a parahoric extension of $H_F$ to a smooth affine group scheme $H$ over $X$.

Let $E/F$ be a cyclic extension of $F$ having degree $p$, so $E$ corresponds to the function field of a smooth projective curve $X'$. Define $G := \Res_{X'/X} (H_{X'})$, which is an algebraic group scheme over $X$ with generic fiber $G_F \cong \Res_{E/F}(H_E)$. The group scheme $G \rightarrow X$ comes with an induced action of $\langle \sigma \rangle = \Aut(X'/X)$.

\subsection{Review of V. Lafforgue's global Langlands correspondence}\label{ssec: VLaff decomposition}

Write $\Gamma := \Weil(F, \ol{F}) = \Weil(F^s/F)$. In \cite[\S 13]{Laff18}, Lafforgue constructs an action of $\Exc(\Gamma, \ld G^{\alg})$ on the space of cusp forms for $G$ with coefficients in $k$. This has been improved by Cong Xue, who extended the action to all compactly supported functions (\cite[\S 7]{Xue20} for split $G$ and \cite[\S 6]{Xue21} for all $G$).

We summarize the construction of the excursion action, as we shall make use of some of its internal aspects, and we also need to explain why it can be used to construct some excursion actions on Tate cohomology.

\subsubsection{Constructing actions of the excursion algebra}\label{sssec: admissible functors}

We will explain an abstract setup that gives rise to actions of the excursion algebra. 

\begin{defn}
Let $A$ be a (not necessarily commutative) ring. A family of functors $H_I \co \Rep_k((\ld G)^I) \rightarrow \Mod_A(\Gamma^I)$, where $I$ runs over (possibly empty) finite sets, is \emph{admissible} if it satisfies the two conditions below. 
\begin{enumerate}
\item \emph{(Compatibility with fusion)} For all $\zeta \co I \rightarrow J$, there is a natural isomorphism $\chi_{\zeta}$ between the functors $H_I \circ \Res_{\zeta}$ and $\Res_{\zeta} \circ H_J$ in the diagram: 
\begin{equation}\label{eq: functor fuse}
\begin{tikzcd}
\Rep_k((\ld G)^I) \ar[r, "H_I"] \ar[dr, Rightarrow, "\chi_{\zeta}"] \ar[d, "\Res_\zeta"'] &    \Mod_A(\Gamma^I)\ar[d, "\Res_{\zeta}"] \\
\Rep_k((\ld G)^J)   \ar[r, "H_J"'] &    \Mod_A(\Gamma^J)  
\end{tikzcd}
\end{equation}
\item \emph{(Compatibility with composition)} For $I' \xrightarrow{\zeta'} I \xrightarrow{\zeta} J $, we have $\chi_{\zeta \circ \zeta'} = \chi_{\zeta} \circ \chi_{\zeta'}$. 
\end{enumerate}
\end{defn}

\begin{const}\label{const: excursion construction}
Let $\bbm{1}$ denote the trivial representation of $\ld {G}$. Given an admissible family of functors  $H_I \co \Rep_k((\ld G)^I) \rightarrow \Mod_A(\Gamma^I)$, we get an $A$-linear action of $\Exc(\Gamma, \ld{G})$ on $H_{\{0\}}(\bbm{1})$ as follows.

For a tuple $(I, W, x, \xi, (\gamma_i)_{i \in I})$ we define an endomorphism, which gives the image of $S_{I, W, x, \xi, (\gamma_i)_{i \in I}}$ in $\End_A(H_{\{0\}}(\bbm{1}))$, by the following composition: 
\[
\begin{tikzcd}
 H_{\{0\}}(\bbm{1})  \ar[r, "H_{\{0\}}(x)"]  &  H_{\{0\}}( W^{\zeta})  \ar[r, "\sim"', "\chi_{\zeta}"] & H_{I}(W)  \ar[r, "(\gamma_i)_{i \in I}"]  & H_I(W)   \ar[r, "\sim"', "\chi_{\zeta}^{-1}"]  &   H_{\{0\}}(W^{\zeta})  \ar[r, "H_{\{0\}}(\xi)"]  & H_{\{0\}}(\bbm{1}).
\end{tikzcd}
\] 
(Here we again conflate the vector $x \in W^{\Delta(\wh{G})}$ with a $\wh{G}$-equivariant map $\bbm{1} \rightarrow W|_{\Delta(\wh{G})}$.) From the assumptions of admissibility it is straightforward to check the relations in \S \ref{sssec: relations 2}.
\end{const}

\begin{remark}\label{rem: trivial value}
Note that it follows from admissibility that the $A$-module underlying $H_I(\bbm{1})$ for any $I$ is identified with $H_{\emptyset}(\bbm{1})$ by $\chi_{\emptyset \rightarrow I}$. Proposition \ref{prop: galois rep} then attaches a Galois representation to each generalized eigenvector for the $\Exc(\Gamma, \ld G)$-action on $H_{\emptyset}(\bbm{1})$. (Of course, such an eigenvector is not guaranteed to exist a priori.)
\end{remark}

\subsubsection{Excursion action on the cohomology of shtukas}

Let $\Cal{H}_G$ be the Hecke algebra acting on $\Sht_{G, D,I}$; it is the tensor product of local Hecke algebras with the level structure dictated by $D$. For any finite set $I$, we have a map 
\[
\pi_I \co \Sht_{G,D, I} \rightarrow (X-D)^I
\]
remembering the points of the curve indexed by $I$ (which avoid $D$ by definition). Let $\eta^I$ denote the generic point of $X^I$ and $\ol{\eta^I}$ the spectrum of an algebraic closure, viewed as a geometric generic point of $X^I$. When $I$ is a singleton, we will just abbreviate these by $\eta$ and $\ol{\eta}$.

We will define a family of functors indexed by finite sets $I$: 
\begin{equation}\label{eq: shtuka admissible family}
H_I^j \co \Rep_k((\ld {G}^{\alg})^I) \rightarrow \Mod_{\Cal{H}_G}(\Gamma^I)
\end{equation}
 sending $V \in \Rep_k((\ld G^{\alg})^I)$ to 
\begin{equation}\label{eq: H_I group}
H_c^j(\Sht_{G,D, I}|_{\ol{\eta^{I}}}; \Sat(V)).
\end{equation}
Note that a priori $H_I^j(V)$ has an action of $\pi_1(\eta^{I}, \ol{\eta^{I}})$, which maps\footnote{The map is non-canonical: it depends on a choice of specialization as in \cite[Remark 8.18]{Laff18}.} to $\Gamma^{I}$ but neither injectively nor surjectively. 

\subsubsection{} We explain why the action of $\pi_1(\eta^I, \ol{\eta^I})$ extends canonically to an action of $\Gamma^I$. Assume $I$ is non-empty, since otherwise there is nothing to prove. The Satake functor of \S \ref{sssec: satake functor} generalizes to a functor 
\[
\Sat^{(I_1, \ldots, I_r)} \co \Rep_k((\ld G)^I) \rightarrow D^b(\circSht_{G,D, I}^{(I_1, \ldots, I_r)};k),
\]
such that the map
\[
\nu \co \circSht_{G,D, I}^{(I_1, \ldots, I_r)} \rightarrow \circSht_{G,D,I}
\]
has the property that $R\nu_! \Sat^{(I_1, \ldots, I_r)}(V) \cong \Sat(V)$. Furthermore, there are natural isomorphisms 
\[
F_{I_1}^* \Sat^{(I_1, I_2, \ldots, I_r)}(V) \cong \Sat^{(I_2, \ldots, I_r, I_1)}(V),
\]
where $F_{I_1}$ is the partial Frobenius from \S \ref{sssec: partial frob}.

Write $I = \{1, \ldots, n\}$. Thanks to the above properties and \eqref{diag: partial Frob}, the partial Frobenius maps on $\Sht_{G,D, I}^{(\{1\}, \ldots, \{n\})}$ then induce maps 
\[
\Frob_{\{1\}}^* H_I^j(V) \xrightarrow{\sim} H_I^j(V). 
\]
That equips $H_I^j(V)$ with the action of the larger group $\FWeil(\eta^I, \ol{\eta^I}) $ that we now recall, summarizing \cite[Remarque 8.18]{Laff18}. Let $F^I$ denote the function field of $X^I$, so $\eta^I  = \Spec F^I$, and $\ol{F^I}$ an algebraic closure, so we may take $\ol{\eta^I} = \Spec \ol{F^I}$. Write $(F^I)^{\perf}$ for the perfect closure of $F^I$, and $\Frob_{\{i\}}$ for the ``partial Frobenius'' automorphism of $(F^I)^{\perf}$ induced by $\Frob_\ell$ on the $i$th factor. We define 
\[
\FWeil (\eta^I, \ol{\eta^I}) := \{ \gamma \in \Aut_{\ol{\F}_q} (\ol{F^I}) \co \exists (n_i)_{i \in I} \in \Z^I \text{ such that } \gamma|_{(F^I)^{\perf}} = \prod_{i \in I} (\Frob_{\{i\}})^{n_i}\}.
\]
Writing $\pi_1^{\geom}(\eta^I, \ol{\eta^I}) := \ker (\pi_1(\eta^I, \ol{\eta^I}) \xrightarrow{\deg} \wh{\Z})$, this fits into an extension
\[
0 \rightarrow \pi_1^{\geom}(\eta^I, \ol{\eta^I}) \rightarrow \FWeil (\eta^I, \ol{\eta^I}) \rightarrow \Z^I \rightarrow 0.
\]

Fixing a specialization morphism $\ol{\eta^I} \rightsquigarrow \Delta(\eta^{\{1\}})$ induces a surjection 
\[
\FWeil (\eta^I, \ol{\eta^I})  \surj \Weil(\eta, \ol{\eta})^I.
\]
A form of Drinfeld's Lemma \cite[Lemma 7.4.2]{Xue20} is used to show that the action of $\FWeil (\eta^I, \ol{\eta^I})$ on $H_I^j(V)$ factors through $\Weil(F^s/F)^I$.

\begin{example}
Let us unravel 
\begin{equation}\label{eq: initial space}
H_{\{0\}}^0(\bbm{1})  = H^0_c(\Sht_{G,D, \{1\}}|_{\ol{\eta^{\{1\}}}}; \Sat(\bbm{1})).
\end{equation}
By Remark \ref{rem: trivial value} the underlying Hecke module of $H_{\{0\}}(\bbm{1})$ is isomorphic to $H_{\emptyset	}(\bbm{1})$. According to \cite[Remarque 12.2]{Laff18}, this is the space of compactly supported $k$-valued functions on the discrete groupoid 
\begin{equation}\label{eq: Bun_G rational points}
\Bun_{G,D}(\F_\ell) = \coprod_{\alpha \in \ker^1(F,G)}  \left( G_{\alpha}(F) \bs G_{\alpha}(\A_F) / \prod_v K_v \right),
\end{equation}
where $G_{\alpha}$ is the pure inner form of $G$ corresponding to $\alpha$, $K_v = G(\Cal{O}_v)$ for $v \notin D$, and $K_v  = \ker(G(\Cal{O}_v) \rightarrow G_D)$. 

The excursion action preserves the decomposition \eqref{eq: Bun_G rational points}, and so gives an action of $\Exc(\Gamma, \ld G)$ on each $H_c^0(\Sht_{G, D, \emptyset}; \bbm{1})_{\alpha} := C_c^{\infty}(G_{\alpha}(F) \bs G_{\alpha}(\A_F) / \prod_v K_v;k)$. 
\end{example}


The family of functors $H_I^j$ is admissible; this is an immediate consequence of the fact that $\Sat$ is already compatible with composition and fusion. Hence Construction \ref{const: excursion construction} applies to define an action of $\Exc(\Gamma, \ld G)$ on $C_c^{\infty}(\Bun_{G,D}(\F_\ell); k)$. Elements of the image of $\Exc(\Gamma, \ld G)$ in $\End(C_c^{\infty}(\Bun_{G,D}(\F_\ell); k))$ are called ``excursion operators''.  


\subsubsection{Xue's generalization}\label{sssec: Xue}

Lafforgue defined an $\Exc(\Gamma, \ld G)$-action on the finite-dimensional subspace of cuspidal functions $C_{\mrm{cusp}}^{\infty}(\Bun_{G,D}(\F_\ell); k) \subset C_c^{\infty}(\Bun_{G,D}(\F_\ell); k)$. This decomposes $C_{\mrm{cusp}}^{\infty}(\Bun_{G,D}(\F_\ell); k)$ into a direct sum of generalized eigenspaces under the action of $\Exc(\Gamma, \ld G)$. Using Proposition \ref{prop: galois rep}, this decomposition corresponds to a parametrization by Langlands parameters.

Thanks to Xue's extension of the action to $\Exc(\Gamma, \ld G) \acts  C_c^{\infty}(\Bun_{G,D}(\F_\ell); k)$, it is meaningful to speak of Langlands parameters arising from $C_c^{\infty}(\Bun_{G,D}(\F_\ell); k)$. However, since the excursion action does not stabilize any finite-dimensional subspaces of $C_c^{\infty}(\Bun_{G,D}(\F_\ell); k)$ unless they are contained in the space of cusp forms, we must broaden what it means to have an $L$-parameter ``come from'' an automorphic function.

\begin{defn} 
We say that an $L$-parameter $\rho \in H^1(\Gal(F^s/F), \wh{G}(k))$ \emph{arises from $C_c^{\infty}(\Bun_{G,D}(\F_\ell); k)$} if it arises via Proposition \ref{prop: galois rep} from the $\Exc(\Gamma, \ld G)$-action on some irreducible Hecke-subquotient of $C_c^{\infty}(\Bun_{G,D}(\F_\ell); k)$; equivalently, if the corresponding maximal ideal $\mf{m}_{\rho} \subset \Exc(\Gamma, \ld G)$ is in the support of $C_c^{\infty}(\Bun_{G,D}(\F_\ell); k)$ as an $\Exc(\Gamma, \ld G)$-module.

We say  will be called \emph{automorphic} if it arises via Proposition \ref{prop: galois rep} from $C_c^{\infty}(\Bun_{G,D}(\F_\ell); k)$ for some $D$; equivalently if the corresponding maximal ideal $\mf{m}_{\rho} \subset \Exc(\Gamma, \ld G)$ is in the support of $C_c^{\infty}(G_{\alpha}(F) \bs G_{\alpha}(\A_F);k)$ for some $\alpha$. 
\end{defn}



\subsection{Excursion action on the Tate cohomology of shtukas}\label{subsec: fusion functor construction}

For a category $\msf{C}$ with $\sigma$-action, recall that we let $\msf{C}^{B \sigma}$ denote the category of $\sigma$-equivariant objects in $\msf{C}$. This comes equipped with a forgetful functor to $\msf{C}$.

\subsubsection{Tate cohomology of moduli of shtukas} 

If $\sigma$ acts on $G$, it induces an action $V \mapsto {}^{\sigma}V$ on $\Rep(\ld G)$. 

Given a $\sigma$-equivariant representation $V \in \Rep_k( (\ld G^{\alg}) ^I)^{\eq}$, we can form $R\Gamma_c(\Sht_{G,D,I}|_{\ol{\eta^I}}; \Sat(V))$ as above. The $\sigma$-equivariant structure on $V$ equips this with a $\sigma$-equivariant structure; more formally, because $\Sat$ and $\pi_I \co \Sht_{G,D,I} \rightarrow (X \setminus D)^I$ are $\sigma$-equivariant, $R\Gamma_c(\Sht_{G,D,I}|_{\ol{\eta^I}}; \Sat(-))$ lifts to a functor $ \Rep_k( (\ld G^{\alg}) ^I)^{\eq} \rightarrow D^b((X \setminus D)^I; k)^{\eq}$. Hence we can form $T^j(R\Gamma_c (\Sht_{G,D,I}|_{\ol{\eta^I}}; \Sat(V)))$, the Tate cohomology (\S \ref{ssec: tate cohomology}) of $R\Gamma_c(\Sht_{G,D,I}|_{\ol{\eta^I}}; \Sat(V))$. 
To ease notation, we will abbreviate 
\begin{equation}\label{eq: tate cohomology of shtukas}
T^j (\Sht_{G,D,I}; V) := T^j(R\Gamma_c (\Sht_{G,D,I}|_{\ol{\eta^I}}; \Sat(V))).
\end{equation}
Let us explain in what category we regard \eqref{eq: tate cohomology of shtukas}. Since $R\Gamma_c(\Sht_{G,D,I}|_{\ol{\eta^I}}; \Sat(V))$ has commuting actions of $\FWeil(\eta^I, \ol{\eta^I})$ and the Hecke algebra $\Cal{H}_G$ (the former commuting with the $\sigma$-action), its Tate cohomology has commuting actions of $\FWeil(\eta^I, \ol{\eta^I})$ and of $T^0(\Cal{H}_G)$, where Tate cohomology is formed with respect to the $\sigma$-action. A priori we regard \eqref{eq: tate cohomology of shtukas} as a $T^0( \Cal{H}_G)[\FWeil(\eta^I, \ol{\eta^I})]$-module. However, in Appendix \ref{app: B}, we will prove: 

\begin{prop}\label{prop: FWeil on shtukas}
For any $G,D,I,V$, the $\FWeil(\eta^I, \ol{\eta^I})$-action on $T^j (\Sht_{G,D,I}; V)$ factors through $\FWeil(\eta^I, \ol{\eta^I}) \surj \Weil(\eta, \ol{\eta})^I$. 
\end{prop}

It will then be natural to regard \eqref{eq: tate cohomology of shtukas} as a $T^0 (\Cal{H}_G)[\Weil(\eta, \ol{\eta})^I]$-module.

Using Lemma \ref{lem: tate cohomology trivial action} we deduce the following simple but important identity: if $\sigma$ acts trivially on $\Sht_H$ and $\Cal{F}$, then 
\begin{equation}\label{eq: comparing cohomology and Tate cohomology}
T^* (\Sht_{H,D,I}; \Cal{F}) \cong H_c^*(\Sht_{H,D,I}|_{\ol{\eta^I}}; \Cal{F}) \otimes T^*(k).
\end{equation}



\subsubsection{Excursion action}

Since $\sigma$ acts on $G$, it acts on $\Exc(\Gamma, \ld G^{\alg})$ by transport of structure. Concretely, we have
\begin{equation}\label{eq: sigma-action on S}
\sigma \cdot S_{V, x, \xi, (\gamma_i)_{i \in I}} = S_{\sigma(V), \sigma(x), \sigma(\xi), (\gamma_i)_{i \in I}}.
\end{equation}
In general given a $k[\sigma]$-algebra $A$ and an $A$-module $M$, there is a natural $T^0(A) = A^\sigma/(N \cdot A)$-module structure on $T^*(M)$. This equips $T^*(\Sht_{G, D, \emptyset}; \bbm{1})$ with a natural $\Exc(\Gamma, \ld G)^{\sigma}$-action. If all the data $(V, x, \xi)$ is $\sigma$-equivariant, then the action of $S_{V, x, \xi, (\gamma_i)_{i \in I}} \in \Exc(\Gamma, \ld G)^{\sigma}$ can be described more concretely as follows: it is given by composition
\begin{equation}
\begin{tikzcd}
T^*(\Sht_{G,D,\emptyset}; \bbm{1}) \ar[r, "x"] & T^*(\Sht_{G,D,\emptyset}; V ) \ar[r, "{(\gamma_i)_{i \in I}}"] & T^*(\Sht_{G,D,\emptyset}; V )  \ar[r, "\xi"] & 
T^*(\Sht_{G,D,\emptyset}; \bbm{1}) .
\end{tikzcd}
\end{equation}
Here we used Proposition \ref{prop: FWeil on shtukas} to define the middle arrow. 

\subsection{Preparations for equivariant localization}
\subsubsection{Analysis of fixed points} We study the $\sigma$-fixed points of $\Sht_{G,D,I}$, in anticipation of applying the theory of \S \ref{sec: generalities} to it. 

According to \cite[Proposition 2.16]{Var04} (stated there for split $G$, but valid for all $G$ by the same argument), $\Sht_{G,D,I}$ is exhausted by quasi-compact open substacks $\Sht_{G,D,I}^{\leq \mu}$ as $\mu$ runs over dominant coweights and the Harder-Narasimhan truncation $\Sht_{G,D,I}^{\leq \mu}$ is defined as in \cite[(1.3)]{Laff18}. The open substack is determined by the Cartesian square 
\[
\begin{tikzcd}
\Sht_{G,D,I}^{\leq \mu} \ar[r, hook] \ar[d] & \Sht_{G,D,I} \ar[d]  \\
\Bun_G^{\leq \mu} \ar[r, hook] & \Bun_G
\end{tikzcd}
\]
Furthermore, for fixed $\mu$ the Deligne-Mumford stack $\Sht_{G,D,I}^{\leq \mu}$ can be presented as a quotient of a quasi-projective scheme by a finite group; for any closed point $x_0 \in X$, the quasi-projective scheme can be taken to be $\Sht_{G, D + n x_0, I}^{\leq \mu}$ for sufficiently large $n$ relative to $\mu$, and the group is then the automorphisms of the level structure. The same applies for the variants $\Sht_{G,D,I}^{(I_1, \ldots, I_r)}$.

We fix the following notation below. Let $\mu$ be a coweight of $H$ and let $\wt{\mu}$ be the induced coweight of $G$. Then we have a Cartesian square
\[
\begin{tikzcd}
\Bun_H^{\leq \mu} \ar[r] \ar[d] & \Bun_H  \ar[d]\\
\Bun_G^{\leq \wt{\mu}} \ar[r] & \Bun_G 
\end{tikzcd}
\]
which induces the Cartesian square 
\begin{equation}\label{eq: Cartesian HN}
\begin{tikzcd}
\Sht_{H,D,I}^{(I_1, \ldots, I_r), \leq \mu} \ar[r] \ar[d] & \Sht_{H,D,I}^{(I_1, \ldots, I_r)}  \ar[d]\\
\Sht_{G,D,I}^{(I_1, \ldots, I_r), \leq \wt{\mu}} \ar[r] & \Sht_{G,D,I}^{(I_1, \ldots, I_r)}
\end{tikzcd}
\end{equation}

\begin{lemma}\label{lemma: shtuka fixed point} If $n$ is sufficiently large so that $\Sht_{H,D + nx_0,I}^{(I_1, \ldots, I_r), \leq \mu}$ and $\Sht_{G,D + nx_0,I}^{(I_1, \ldots, I_r), \leq \wt{\mu}}$ are representable by schemes, then the diagonal map $H \rightarrow G$ induces an isomorphism 
\[
\Sht_{H,D+nx_0,I}^{(I_1, \ldots, I_r), \leq \mu} \xrightarrow{\sim} (\Sht_{G,D+nx_0,I}^{(I_1, \ldots, I_r), \leq \wt{\mu}})^{\sigma}.
\]
\end{lemma}

\begin{proof}
For notational convenience we just treat the case of non-iterated shtukas, $\Sht_{G,D,I}$; the general case is essentially the same but with cumbersome extra notation. 

There is an obvious map in one direction, $\Sht_{H,D+nx_0,I}^{ \leq \mu} \rightarrow (\Sht_{G,D+nx_0,I}^{\leq \wt{\mu}})^{\sigma}$. We will construct the inverse. 

Notate the $S$-points of $\Sht_{G,D+nx_0,I}^{\leq \wt{\mu}}$ as the set $\{(\{x_i\}_{i \in I}, \Cal{E}, \varphi, \upsilon)\}$. For any $S$, there is an equivalence of categories between $\Res_{X'/X}(H)$-torsors on $X_S$ and $H$-torsors on $X'_S$, which we denote $\Cal{E} \mapsto \Cal{E}'$. The datum of a $\sigma$-fixed point of $\Bun_{G,D}$ translates under the above equivalence to the datum of an $H$-torsor $\Cal{E}'$ on $X'_S$ together with an isomorphism $h \co \cE' \xrightarrow{\sim} \sigma^* \cE'$. We claim that, since the point $(\{x_i\}_{i \in I}, \Cal{E}, \varphi, \upsilon)$ has no non-trivial automorphisms, such an isomorphism automatically satisfies the cocycle condition, hence is equivalent to a descent datum from $\cE'$ to an $H$-torsor over $X_S$. Furthermore, the map $\varphi$ and level structure $\upsilon$ will similarly descend uniquely. 

Let $\Nm(h) := (\sigma^{p-1}h) \circ  \ldots \circ (\sigma h) \circ h \co \cE' \xrightarrow{\sim}   \cE'$. The claim amounts to checking that $\Nm(h)$ is the identity automorphism of $\cE'$. By definition, it corresponds to some automorphism of $\cE$ compatible with $\varphi$ and the level structure $\upsilon$. But by assumption, this datum had no non-trivial automorphisms, so $\Nm(h)$ can only be the identity automorphism. 

This constructs a map $\Sht_{H,D+nx_0,I} \leftarrow (\Sht_{G,D+nx_0,I}^{\leq \wt{\mu}})^{\sigma}$ which is manifestly a one-sided inverse; we conclude by using \eqref{eq: Cartesian HN} to see that it lands in $\Sht_{H,D+nx_0,I}^{\leq \mu}$. 

\end{proof}

\subsubsection{Cohomology at infinite level}\label{sssec: infinite level} We will use Lemma \ref{lemma: shtuka fixed point} to apply Smith theory. However, the excursion action does not stabilize the piece of cohomology coming from bounding the HN polygon, so we need to let $\mu$ and $n$ both go ``off to infinity''. 

\begin{defn}Fix a closed point $x_0 \in X$ and consider the system of Deligne-Mumford stacks, $\{\Sht_{H,D+nx_0,I}^{\leq \mu}\}$ as $n$ and $\mu$ vary. For $V \in \Rep((\ld H^{\alg})^I)$, we \emph{define}
\[
R\Gamma_c(\Sht_{H, D + \infty x_0, I}|_{\ol{\eta^I}}; \Sat(V)) = \colim_{n,\mu}  R\Gamma_c (\Sht_{H, D + n x_0, I}^{\leq \mu}|_{\ol{\eta^I}}; \Sat(V))
\]
where the maps in the $\mu$ variable are the covariant maps induced by open embeddings, while the maps in the $n$ variable are the contravariant maps induced by pullback. Note that the colimit is \emph{filtered} because both indexing posets are filtered. 
\end{defn}

\begin{remark}\label{rem: scheme cofinal}
As explained above, for any fixed $\mu$, and all sufficiently large $n$ depending on $\mu$, $\Sht_{H, D + n x_0, I}^{\leq \mu}$ is representable by a scheme. Hence, the subposet of indices $(n, \mu)$ for which $\Sht_{H, D + n x_0, I}^{\leq \mu}$ is representable by a scheme is cofinal, so $R\Gamma_c(\Sht_{H, D + \infty x_0, I}|_{\ol{\eta^I}}; \Sat(V))$ is naturally isomorphic to the colimit taken along this subposet.
\end{remark}

\begin{defn}
Fix a closed point $x_0 \in X$ and $V \in \Rep_k( (\ld H^{\alg}) ^I)^{\eq}$. We define 
\[
T^j(\Sht_{H,D + \infty x_0,I};V) := T^j (R\Gamma_c(\Sht_{H, D + \infty x_0, I}|_{\ol{\eta^I}}; \Sat(V))) .
\]
We note that $R\Gamma_c(\Sht_{H, D + \infty x_0, I}|_{\ol{\eta^I}}; \Sat(V)) $ is bounded, since the dimension of the support of $\Sat(V)$ on each $\Sht_{H, D + n x_0, I}^{\leq \mu}|_{\ol{\eta^I}}$ is uniformly bounded for all $n,\mu$. 
\end{defn}

Furthermore, note that for any cofinal subposet of HN polygons $\mu$ for $H$, the induced HN polygons $\wt{\mu}$ form a cofinal poset for $G$.  We make the analogous definitions $R\Gamma_c(\Sht_{G, D + \infty x_0, I}|_{\ol{\eta^I}}; \Sat(V))$ and $T^j(\Sht_{G,D + \infty x_0,I};V)$ for $G$.

\subsection{Equivariant localization for excursion operators}

We define $\Nm \co \Rep_k((\ld{G})^I) \rightarrow \Rep_k((\ld{G})^I)^{\eq}$ to be the functor taking a representation $V$ to $V \otimes_k {}^{\sigma} V \otimes_k \ldots \otimes_k {}^{\sigma^{p-1}}V$, with the $\sigma$-equivariant structure
\[
{}^{\sigma} \Nm(V) = {}^{\sigma} V \otimes_k {}^{\sigma^2} V \otimes_k \ldots \otimes_k {}^{\sigma^{p-1}}V \otimes_k V  \xrightarrow{\sim}  V \otimes_k {}^{\sigma} V \otimes_k \ldots \otimes_k  {}^{\sigma^{p-1}}V  = \Nm(V)
\]
given by the commutativity constraint for tensor products. It corresponds under Geometric Satake to Definition \ref{defn: Nm}. Given $h \co  V \rightarrow V' \in \Rep_k((\ld{G})^I)$, we set 
\[
\Nm(h) := h \otimes {}^{\sigma} h \otimes \ldots \otimes {}^{\sigma^{p-1}} h  \co \Nm(V) \rightarrow \Nm(V').
\]Note that $\Nm$ is \emph{not} an additive functor, nor is it even $k$-linear. We linearize it by defining $\Nm^{(p^{-1})} :=  \Nm \circ \Frob_p^{-1}$, where (as in \S \ref{sssec: base change functor}) $ \Frob_p^{-1}$ is the identity on objects and on morphisms it is $(-) \otimes_{k, \Frob_p^{-1}} k$. Then $\Nm^{(p^{-1})} \co \Rep_k((\ld{G})^I) \rightarrow \Rep_k((\ld{G})^I)^{\eq}$ is $k$-linear, although still not additive.

For $V \in \Rep_k((\ld{G})^I)$, we denote by $N \cdot V$ the $\sigma$-equivariant representation $V \oplus {}^{\sigma} V \oplus \ldots  \oplus {}^{\sigma^{p-1}}V$, with $\sigma$-equivariant structure 
\[
{}^{\sigma} (N \cdot V) = {}^{\sigma}V \oplus {}^{\sigma^2} V \oplus \ldots  \oplus {}^{\sigma^{p-1}}V  \oplus V  \xrightarrow{\sim} V \oplus {}^{\sigma} V \oplus \ldots \oplus   {}^{\sigma^{p-1}}V = (N \cdot V)	
\]
given by the commutativity constraint for direct sums. For $h \co  V \rightarrow V' \in \Rep_k((\ld{G})^I)$, we denote by $N \cdot h \co N \cdot V \rightarrow N \cdot V'$ the $\sigma$-equivariant map $h \oplus {}^{\sigma} h \oplus \ldots \oplus {}^{\sigma^{p-1}} h$. Let $\Delta_p \co \bbm{1} \rightarrow \bbm{1}^{\oplus p}$ denote the diagonal map and $\nabla_p \co \bbm{1}^{\oplus p} \rightarrow \bbm{1}$ denote the sum map. 

Our goal in this subsection is to prove the theorem below. 

\begin{thm}\label{thm: global actions coincide} Fixed a closed point $x_0$ on $X$ and let $D$ be any closed finite subscheme of $X$. 

(i) For each $j \in \{0,1\}$, with respect to the isomorphism 
\[
T^j(\Sht_{G, D + \infty x_0, \emptyset}; \bbm{1}) \cong H_c^j(\Sht_{H, D + \infty x_0, \emptyset};\bbm{1})
\]
induced by Lemma \ref{lemma: shtuka fixed point}, the action of $S_{I, \Nm^{(p^{-1})}(V), \Nm^{(p^{-1})}(x), \Nm^{(p^{-1})}(\xi), (\gamma_i)_{i \in I}} \in \Exc(\Gamma, \ld G)^{\sigma}$ on $T^j (\Sht_{G, D + \infty x_0, \emptyset}; \bbm{1})$ is identified with the action of $S_{I, \Res_{\BC}(V), x, \xi, (\gamma_i)_{i \in I}} \in \Exc(\Gamma, \ld H)$ on $T^j(\Sht_{H, D + \infty x_0, \emptyset};\bbm{1}) = H_c^0(\Sht_{H, D + \infty x_0, \emptyset};\bbm{1})$. 

(ii) For each $j \in \{0,1\}$, the action of $S_{I, N \cdot V, (N \cdot x) \circ \Delta_p, \nabla_p \circ (N \cdot \xi), (\gamma_i)_{i \in I}}\in \Exc(\Gamma, \ld G)^{\sigma}$ on $T^j (\Sht_{G, D + \infty x_0, \emptyset}; \bbm{1})$ is $0$. 
\end{thm}

This will be established in several steps. The heart of the matter is the following equivariant localization isomorphism. 

\begin{prop}\label{prop: equiv localization for shtukas} 
Fixed a closed point $x_0$ on $X$ and let $D$ be any closed finite subscheme of $X$. Let $G$ be a reductive group over $F$. We equip $\Rep_k(\wh{G}^I) $ with the $\pi_1(\eta, \ol{\eta})^I$-action coming from the geometric action of $\pi_1(\eta, \ol{\eta})$ on $\wh{G}$ (\S \ref{sssec: geometric L-group}). 

(i) For each $j \in \{0,1\}$, there is a natural isomorphism of $\pi_1(\eta^I, \ol{\eta}^I)$-equivariant functors $\Rep_k(\wh{G}^I)  \rightarrow \mrm{Vect}_{/k}$,
\begin{equation}\label{eq: rep equiv}
T^j(\Sht_{G, D+\infty x_0, I}; \Nm^{(p^{-1})} (V)) \cong T^j(\Sht_{H,D+\infty x_0,I}; \Res_{\BC}(V)), \quad V \in \Rep_k(\wh{G}^I).
\end{equation}
(ii) For each $j \in \{0,1\}$, there is a natural isomorphism of functors $\Rep_k((\ld{G})^I) \rightarrow \Mod_k(\Weil(\eta, \ol{\eta})^I)$,
\begin{equation}\label{eq: equiv localization for shtukas}
T^j(\Sht_{G,D + \infty x_0,I}; \Nm^{(p^{-1})}(V)) \cong T^j(\Sht_{H,D+ \infty x_0,I}; \Res_{\BC}(V)), \quad V \in \Rep_k((\ld G)^I),
\end{equation}
which is compatible with fusion and composition. 
\end{prop}

\begin{proof}(i) We will first construct a natural isomorphism  of $\pi_1(\eta^I, \ol{\eta}^I)$-equivariant functors $\Tilt_k(\wh{G}^I)  \rightarrow \Shv(\ol{\eta}^I; \cT_k)$,
\begin{equation}\label{eq: tilt equiv}
\TT^* R\Gamma_c(\Sht_{G, D+\infty x_0, I}; \Nm^{(p^{-1})} (V)) \cong \TT^* \epsilon^* R\Gamma_c(\Sht_{H,D+\infty x_0,I}; \Res_{\BC}(V)), \quad V \in \Tilt_k(\wh{G}^I).
\end{equation}
By Lemma \ref{lem: integral cube} (specifically, the commutativity of the back face), the diagram 
\[
\begin{tikzcd}
\Tilt_k(\wh{G}) \ar[r, "{\Sat}"]  &  \Parity^0_{L^+G}(\Gr_G;k) \ar[r, "{\Nm}"]  &   \Parity^0_{L^+G \rtimes \sigma}(\Gr_G;k) \ar[d]  \ar[r] & \Parity^0_{L^+H}(\Gr_H; k) \ar[d, "\TT^* \epsilon^*"] \\
& &  D_{L^+G \rtimes \sigma, c}^b(\Gr_G; k) \ar[r, "{\Psm}"]  & \Perf_{(L^+H)}(\Gr_H; \cT_k)  
\end{tikzcd}
\]
commutes. The composite functor that follows along the left and then bottom is $\Psm \circ \Nm \circ \Sat $, while Theorem \ref{thm: BC on Gr} identifies the composite functor that follows along the top and right with $\TT^* \epsilon^* \Res_{\BC^{(p)}}$. We therefore have a $\pi_1(\eta^I, \ol{\eta^I})$-equivariant natural isomorphism between these two composite functors. (The only non-tautological aspect of the equivariance is handled in Lemma \ref{lem: bc equivariance}.) Linearizing this natural isomorphism, we get a $\pi_1(\eta^I, \ol{\eta^I})$-equivariant natural isomorphism of functors $\Tilt_k(\wh{G}^I) \rightarrow  \Perf_{(L^+H)}(\Gr_H; \cT_k) $,
\begin{equation}\label{eq: revise 3 eq 0.5}
 \Psm \circ  \Nm \circ \Sat_G \circ \Frob_p^{-1} \cong \Sat_H \circ \Res_{\BC}. 
\end{equation}

For any coweight $\mu$ for $H$, inducing the coweight $\wt{\mu}$ for $G$, and $n$ sufficiently large so that $\Sht_{G, D+nx_0,I}^{\leq \wt{\mu}}$ and $\Sht_{H,D+nx_0,I}^{\leq \mu}$ are schemes, Lemma \ref{lemma: shtuka fixed point} identifies the $\sigma$-fixed points of $\Sht_{G, D+\infty x_0, I}$ with $\Sht_{H, D + \infty x_0, I}$. Therefore we may apply \eqref{eq: equivariant localization} to obtain an isomorphism in $\Shv(\ol{\eta^I}; \cT_k)$,
\begin{align}\label{eq: revise 3 eq 1}
\TT^* R\Gamma_c  (\Sht_{G, D+nx_0,I}^{\leq \wt{\mu}}; \Nm^{(p^{-1})}(V)) & := \TT^* R\Gamma_c (\Sht_{G,D+nx_0,I}^{\leq \wt{\mu}}; \Sat(\Nm^{(p^{-1})}(V))) \nonumber \\
&   \cong  R\Gamma_c(\Sht_{H,D+nx_0,I}^{\leq \mu};  \Psm( \Nm(\Sat( \Frob_p^{-1} V))))
\end{align}
which is natural in $V \in \Tilt_k(\wh{G}^I)$ and $\pi_1(\eta^I, \ol{\eta^I})$-equivariant. The commutative diagram 
\[
\begin{tikzcd}[ampersand replacement=\&]
\Sht_{H,I} \ar[d] \ar[r] \& \Hk_{H,I} \ar[d] \\
\Sht_{G,I} \ar[r] \& \Hk_{G,I}
\end{tikzcd}
\]
induces a natural isomorphism between the two restriction functors $D^b_c(\Hk_{G,I}|_{\ol{\eta^I}};k) \rightarrow D^b_c(\Sht_{H,I}|_{\ol{\eta^I}};k)$, one by $*$-pullback through the right then top maps and the other by $*$-pullback through the bottom then left maps. Furthermore, the pullbacks are $\pi_1(\eta^I, \ol{\eta^I})$-equivariant, as is the natural isomorphism between them. The same discussion applies with any level structure and HN truncation. Therefore, there is no risk of confusion in the expression $\Psm(\Nm(\Sat(V))$ whether we first regard $\Nm(\Sat(V))$ as a sheaf on $\Sht_{G, D+nx_0,I}^{\leq \wt{\mu}}$ and then apply $\Psm$, or first apply $\Psm$ and then pull back to $\Sht_{H,D+nx_0,I}^{\leq \mu}$; the two are naturally identified. Now, \eqref{eq: revise 3 eq 0.5} induces an isomorphism in $\Shv(\ol{\eta^I}; \cT_k)$,
\begin{equation}\label{eq: revise 3 eq 2}
R\Gamma_c(\Sht_{H,D+nx_0,I}^{\leq \mu};  \Psm( \Nm(\Sat(\Frob_p^{-1} V)))) \cong  \TT^* \epsilon^* R\Gamma_c(\Sht_{H,D+nx_0,I}^{\leq \mu};	 \Sat(\Res_{\BC}(V))),
\end{equation}
which is natural in $V \in \Tilt_k(\wh{G}^I)$ and $\pi_1(\eta^I, \ol{\eta^I})$-equivariant.

We conclude \eqref{eq: tilt equiv} by taking the colimit of these isomorphisms \eqref{eq: revise 3 eq 1} and \eqref{eq: revise 3 eq 2} along such $n$ and $\mu$, using that they form a cofinal poset by Remark \ref{rem: scheme cofinal}.

Next we bootstrap from $\Tilt_k(\wh{G}^I)$ to $\Rep_k(\wh{G}^I)$. For this, we use Proposition \ref{prop: tiltings generate}, which allows for any $V \in \Rep_k(\wh{G}^I)$ to produce a resolution of $V$ by a complex $V^{\bu} = (\ldots  \rightarrow V^{-1} \rightarrow V^0 \rightarrow \ldots)$ of tilting modules which is well-defined up to homotopy. Then Lemma \ref{lem: exactness} gives a natural isomorphism in $\Perf(\Sht_{H,D+nx_0,I}^{\leq \mu}; \cT_k)$,
\[
\Psm(\Sat(\Nm^{(p^{-1})}V^{\bu})) =( \ldots \rightarrow  \Psm(\Sat(\Nm^{(p^{-1})} V^{-1}))  \rightarrow \Psm(\Sat(\Nm^{(p^{-1})} V^0)) \rightarrow \ldots)  
\]
So, using the earlier observations and \eqref{eq: equivariant localization}, we have $\pi_1(\eta^I, \ol{\eta^I})$-equivariant natural isomorphisms in $\Shv(\ol{\eta^I}; \cT_k)$, 
\begin{align*}
\TT^* R\Gamma_c(\Sht_{G,D+nx_0,I}^{\leq \mu};  \Nm^{(p^{-1})} V)  & \cong R\Gamma_c(\Sht_{H,D+nx_0,I}^{\leq \mu}; \Psm  \Nm^{(p^{-1})} (V))   \\
& 
\cong   R\Gamma_c(\Sht_{H,D+nx_0,I}^{\leq \mu}; \Psm  \Nm^{(p^{-1})} (V^{\bu})).
\end{align*}
Then using \eqref{eq: tilt equiv}, we have 
\begin{align*}
 R\Gamma_c(\Sht_{H,D+nx_0,I}^{\leq \mu}; \Psm  \Nm^{(p^{-1})} (V^{\bu}))  &  \cong \TT^* \epsilon^*  R\Gamma_c(\Sht_{H,D+nx_0,I}^{\leq \mu}; \Res_{\BC}(V^{\bu})) \\
 & \cong \TT^* \epsilon^* R\Gamma_c(\Sht_{H,D+nx_0,I}^{\leq \mu}; \Res_{\BC}(V)).
\end{align*}
Now take the $j$th Tate cohomology group to obtain \eqref{eq: rep equiv}.

(ii) Since $\FWeil(\eta^I, \ol{\eta^I}) \surj \Weil(\eta, \ol{\eta})^I$, it suffices to show a natural isomorphism as $\FWeil(\eta^I, \ol{\eta}^I)$-modules. Then, since the $\FWeil(\eta^I, \ol{\eta^I})$-actions on $T^*(\Sht_{G,D+ \infty x_0, I}; \Nm^{(p^{-1})}(V)) $ and on $T^*(\Sht_{H, D+ \infty x_0, I}; \Res_{\BC}(V)) $ are determined by their respective $\pi_1(\eta^I, \ol{\eta^I})$-actions plus partial Frobenius morphisms, we can and will focus on these two equivariance structures separately.

Applying $(-)^{B \pi_1(\eta^I, \ol{\eta^I})}$ (see the references on equivariantization in \S \ref{sssec: equivariantization}) to the natural isomorphism in (i) gives a natural isomorphism 
\[
T^j(\Sht_{G,D + \infty x_0,I}; \Nm^{(p^{-1})}(V)) \cong T^j(\Sht_{H,D+ \infty x_0,I}; \Res_{\BC}(V))
\]
of functors $V \in \Rep_k((\ld G)^I) \rightarrow \Mod_k(\pi_1(\eta^I, \ol{\eta^I}))$.

Finally, we check the compatibility with partial Frobenius. We want to show that the diagram 
\begin{equation}\label{eq: commutative diagram}
\begin{tikzcd}
F_{\{1\}}^* T^j(\Sht_{G,D+ \infty x_0,I}; \Nm^{(p^{-1})}(V)) \ar[r, "{\sim}"] \ar[d, "\sim"] &  T^j(\Sht_{G,D+ \infty x_0,I}; \Nm^{(p^{-1})}	(V)) \ar[d,"\sim"] \\
F_{\{1\}}^* T^j(\Sht_{H,D+ \infty x_0,I}; \Res_{\BC}(V)) \ar[r, "\sim"]  &  T^j(\Sht_{H,D+ \infty x_0,I}; \Res_{\BC}(V))
\end{tikzcd}
\end{equation}
commutes, where the vertical isomorphisms (as $\pi_1(\eta^I, \ol{\eta^I})$-modules) have just been established. By Lemma \ref{lemma: shtuka fixed point}, there is a cofinal system of $n,\mu,\mu'$ such that applying $\sigma$-fixed points to the diagram 
\[
F_{\{1\}} \co \Sht_{G,D+n x_0,I}^{(\{1\},  \{2\}, \ldots, \{r \}), \leq \wt{\mu}} \rightarrow \Sht_{G,D + nx_0,I}^{(\{2\}, \ldots, \{r\}, \{1\}), \leq \wt{\mu}'} 
\]
yields the diagram 
 \[
 F_{\{1\}} \co \Sht_{H,D+nx_0,I}^{(\{1\}, \{2\}, \ldots, \{r \}), \leq \mu} \rightarrow \Sht_{H,D+nx_0,I}^{(\{2\}, \ldots, \{r\}, \{1\}), \leq \mu'}.
 \]
 (The need for $\mu'$ arises because $F_{\{1\}}$ does not preserve HN polygons.) This implies that the natural isomorphims \eqref{eq: revise 3 eq 1} and \eqref{eq: revise 3 eq 2} are compatible with the maps $F_{\{1\}}^*$. Taking the (filtered) colimit along such $n,\mu,\mu'$ completes the proof.

\end{proof}

\begin{proof}[Proof of Theorem \ref{thm: global actions coincide}] (i) Proposition \ref{prop: equiv localization for shtukas}(ii) gives a chain of compatible identifications 
\[
\adjustbox{scale=0.75,center}{\begin{tikzcd}[column sep = huge]
T^j(\Sht_{G,D + \infty x_0,I}; \bbm{1}) \ar[r, "\Nm^{(p^{-1})}(x)"] \ar[d, "\sim"] & T^j(\Sht_{G,D+ \infty x_0,I}; \Nm^{(p^{-1})}(V) ) \ar[d, "\sim"] \ar[r, "{(\gamma_i)_{i \in I}}"] & T^j(\Sht_{G,D+ \infty x_0,I}; \Nm^{(p^{-1})}(V) ) \ar[d, "\sim"]  \ar[r, "\Nm^{(p^{-1})}(\xi)"] & T^j(\Sht_{G,D+ \infty x_0,I}; \bbm{1}) \ar[d, "\sim"] \\
T^j(\Sht_{H,D+ \infty x_0,I}; \bbm{1}) \ar[r, "x"] & T^j(\Sht_{H,D+ \infty x_0,I}; \Res_{\BC}(V) ) \ar[r, "{(\gamma_i)_{i \in I}}"] & T^j(\Sht_{H,D+ \infty x_0,I}; \Res_{\BC}(V) )  \ar[r, "\xi"] & T^j(\Sht_{H,D+ \infty x_0,I}; \bbm{1})
\end{tikzcd}}
\]
The operator $S_{I, \Nm^{(p^{-1})}(V), \Nm^{(p^{-1})}(x), \Nm^{(p^{-1})}(\xi), (\gamma_i)_{i \in I}}$ on $T^j(\Sht_{G,D + \infty x_0,I}; \bbm{1})$ is obtained by tracing along the upper row, while the operator $S_{I, \Res_{\BC}(V), x, \xi, (\gamma_i)_{i \in I}}$ on $T^j(\Sht_H; \bbm{1})$ is obtained by tracing along the lower row. Hence they coincide under the vertical identifications.

(ii) By Lemma \ref{lemma: shtuka fixed point} and \eqref{eq: equivariant localization} we have a chain of compatible identifications
 \[
\adjustbox{scale=0.75,center}{\begin{tikzcd}[column sep = huge]
T^j(\Sht_{G,D + \infty x_0,I}; \bbm{1}) \ar[r, "(N \cdot x) \circ \Delta_p"] \ar[d, "\sim"] & T^j(\Sht_{G,D + \infty x_0,I}; N \cdot V) \ar[d, "\sim"] \ar[r, "{(\gamma_i)_{i \in I}}"] & T^j(\Sht_{G,D + \infty x_0,I}; N \cdot V) \ar[d, "\sim"]  \ar[r, "\nabla_p \circ (N \cdot \xi)"] & T^j(\Sht_{G,D + \infty x_0,I}; \bbm{1}) \ar[d, "\sim"] \\
T^j(\Sht_{H,D + \infty x_0,I}; \bbm{1}) \ar[r, "(N \cdot x) \circ \Delta_p"]  & T^j(\Sht_{H,D + \infty x_0,I}; \Psm( N \cdot V) )  \ar[r, "{(\gamma_i)_{i \in I}}"] & T^j(\Sht_{H,D + \infty x_0,I}; \Psm(N \cdot V))  \ar[r, "\nabla_p \circ  (N \cdot \xi)"] & T^j(\Sht_{H,D + \infty x_0,I}; \bbm{1})
\end{tikzcd}
}
\]
The operator $S_{I, N \cdot V, (N \cdot x) \circ \Delta_p, \nabla_p \circ (N \cdot \xi), (\gamma_i)_{i \in I}}$ on $T^j(\Sht_{G,D + \infty x_0,I}; \bbm{1}) $ is obtained by tracing along the upper row. But the stalks and costalks of $N \cdot \Sat(V)|_{\Gr_H}$ are all complexes of induced $k[\sigma]$-modules, and in fact they are perfect complexes over $k[\sigma]$. Hence $\Psm( N \cdot V)$ is equivalent to $0$ in the Tate category of $\Sht_{H,D,I}$ for all $D$, so $T^j(\Sht_{H,D + \infty x_0,I}; \Psm( N \cdot V) )  = 0$. Therefore the endomorphism in question factors through the zero map, hence is itself zero. 
\end{proof}

 \subsection{Applications to base change for automorphic forms}\label{ssec: base change}
 
 In \S \ref{ssec: VLaff decomposition} we described Lafforgue's action of $\Exc(\Gamma, \ld G)$ on $H_{\emptyset}(\bbm{1})$. By \eqref{eq: Bun_G rational points}, we have
\[
 H_{\emptyset}(\bbm{1}) = \bigoplus_{\alpha \in \ker^1(F, G)} C_c^{\infty}(G_{\alpha}(F) \bs G_{\alpha}(\A_F) / \prod_v K_v; k). 
 \]
Here $\ker^1(F,G) := \ker(H^1(F,G) \rightarrow \prod_v H^1(F_v, G)) $ is the isomorphism class of the generic fiber of the $G$-torsor. More generally, this defines a decomposition
\begin{equation}\label{eq: decomp inner forms}
\Sht_{G,D,I} = \coprod_{\alpha \in \ker^1(F,G)} (\Sht_{G,D,I})_{\alpha}
\end{equation}
according to the isomorphism class of the generic fiber of $\Cal{E}$. 

In the base change situation, the ``diagonal embedding'' map $\phi \co H \rightarrow G$ induces a map $\phi_* \co \ker^1(F, H) \rightarrow \ker^1(F, G)$, compatible with the map $\Bun_{H,D}(\F_{\ell}) \rightarrow \Bun_{G,D}(\F_{\ell})$. Theorem \ref{thm: intro 1} is evidently implied by the theorem below, whose proof occupies this subsection.

\begin{thm}\label{thm: global main}
Fix any closed point $x_0 \in X$ and any closed finite subscheme on $X$. Let $[\rho] \in H^1(\Gamma,  \wh{H}(k))$ be an $L$-parameter arising from the action of  $\Exc(\Gamma, \ld{H})$ on $H^0_c(\Sht_{H,D +\infty x_0,I}; \Sat(\bbm{1}))_{\alpha}$ in the sense of \S \ref{sssec: Xue}. Then the image of $[\rho]$ in $H^1(\Gamma, \wh{G}(k))$ arises in the action of $\Exc(\Gamma, \ld G)$ on $H^0_c(\Sht_{G,D + \infty x_0,I}; \Sat(\bbm{1}))_{\phi(\alpha)}$ in the sense of \S \ref{sssec: Xue}.
\end{thm}

We establish some preliminaries in preparation for the proof.

\begin{defn}[The Tate diagonal]\label{def: tate diagonal}
For a commutative algebra $A$ in characteristic $p$ with $\sigma$-action, we denote by $N \cdot A$ the subset consisting of elements of the form $(1+\sigma + \ldots + \sigma^{p-1}) a$ for $a \in A$. One easily checks that $N \cdot A$ is an ideal in $A^{\sigma}$. 

We denote by $\Nm \co A \rightarrow A^{\sigma}$ the set map sending $a \mapsto a \cdot  \sigma(a) \cdot \ldots \cdot \sigma^{p-1}(a)$. It is multiplicative but not additive. It is an exercise to verify that the composition of $\Nm$ with the quotient $A^{\sigma} \surj A^{\sigma}/N \cdot A$ is an algebra homomorphism, which we call the \emph{Tate diagonal homomorphism} $\Delta^p \co A \rightarrow T^0(A)$.
\end{defn}

\begin{lemma}\label{lem: unique extension of char}
Let $A$ be a commutative ring over $\F_p$ with a $\sigma$-action. Let $\kappa$ be any perfect field over $\F_p$. Let $A' \subset A^{\sigma}$ be a subring containing $\Nm (A)$ and $N \cdot A$. (Since $N \cdot A$ is an ideal in $A^{\sigma}$, it is also an ideal in any such $A'$.) Any character $\chi \co A' \rightarrow \kappa$ factoring through $A'/N\cdot A$ extends uniquely to a character $\wt{\chi} \co A \rightarrow \kappa$, which is expclitily given by 
\begin{equation}\label{eq: ext of hom}
\wt{\chi}(a) = \chi(\Nm(a))^{1/p}.
\end{equation}
\end{lemma}

\begin{proof} The same proof as that of \cite[\S 3.4]{TV} works, but since our situation is a little more general we reproduce it. One easily checks that the given formula \eqref{eq: ext of hom} defines a valid extension (it is a ring homomorphism since $\kappa$ is in characteristic $p$, and it clearly extends $\chi$). 

Next we check that it is the unique extension. Note that $\sigma$ acts on characters $A \rightarrow \kappa$ by pre-composition; we denote this action by $\wt{\chi} \mapsto \sigma \cdot \wt{\chi} := \wt{\chi} \circ \sigma^{-1}$. Clearly \eqref{eq: ext of hom} is the unique $\sigma$-fixed extension, so we will show that any extension $\wt{\chi}'$ must be $\sigma$-fixed. Indeed, since any extension $\wt{\chi}'$ is trivial on $N \cdot A$ by the assumption that $\chi$ factors through $A'/N \cdot A$, we have 
\[
\sum_{i=0}^{p-1} \sigma^i \cdot \wt{\chi}' = 0.
\]
By linear independence of characters \cite[\href{https://stacks.math.columbia.edu/tag/0CKK}{Tag 0CKK}]{stacks-project} we must have $\sigma^i \cdot \wt{\chi}' = \wt{\chi}'$ for all $i$, i.e. $\wt{\chi}'$ is $\sigma$-fixed. 
\end{proof}

\begin{lemma}\label{lem: norm extension} Inside $\Exc(\Gamma, \ld G)$ we have 
\[
\Nm (S_{I, V, x, \xi, (\gamma_i)_{i \in I}}) = S_{I, \Nm (V), \Nm(x), \Nm(\xi), (\gamma_i)_{i \in I}}
\]
and 
\[
N \cdot S_{I, V, x, \xi, (\gamma_i)_{i \in I}} = S_{I, N \cdot V, (N \cdot x) \circ \Delta_p, \nabla_p \circ (N \cdot \xi), (\gamma_i)_{i \in I}}.
\]
\end{lemma}

\begin{proof}
The first equality follows from repeated application of the relations \eqref{eq: relation a-1.5}, \eqref{eq: relation a-2} and the explicit description of the $\sigma$-action in \eqref{eq: sigma-action on S}. The second equality follows from repeated application of relations  \eqref{eq: relation a-1.5}, \eqref{eq: relation a-2b} and the explicit description of the $\sigma$-action in \eqref{eq: sigma-action on S}. 
\end{proof}

\begin{const}[Frobenius twist of algebras]\label{const: frob twist algebra} This discussion is parallel to Construction \ref{const: frob twist C}. Let $\Frob$ be the absolute Frobenius of $k$. Given a $k$-algebra $A$, we denote by $A^{(p)} := A \otimes_{k, \Frob} k$ its Frobenius twist. The map $A \rightarrow A^{(p)}$ sending $a \mapsto a \otimes 1$ is a $k$-semilinear isomorphism. It is characterized by the universal property that any $\Frob$-semilinear homomorphism $f \co A \rightarrow B$ (i.e., $f(\lambda a) = \lambda^p f(a)$ for $\lambda \in k$) factors uniquely through a $k$-linear homomorphism 
\[
\begin{tikzcd}
A \ar[dr, "f"] \ar[d]  \\
A^{(p)} \ar[r, dashed] & B
\end{tikzcd}
\]

If $A$ is equipped with an $\F_p$-structure $\varphi_0 \co A \cong A_0 \otimes_{\F_p} k$, then there is a $k$-linear isomorphism $A \cong A^{(p)}$, characterized by the property that it sends $A_0$ to $A_0 \otimes 1$ via $\Id \otimes 1$. We denote by $\Frob_{\varphi_0} \co A \xrightarrow{\sim} A$ the $\Frob$-semilinear isomorphism which is is the identity on $A_0$; it is the composition of the $k$-linear isomorphism $ A^{(p)} \cong A$ above and the $\Frob$-semilinear isomorphism $A \xrightarrow{\sim} A^{(p)}$ above. The $k$-linear homomorphism $A \rightarrow B$ obtained by precomposing $f$ with the inverse of $\Frob_{\varphi_0}$ will be call the \emph{linearization} of $f$ (with respect to $\varphi_0$). 
\end{const}

\begin{example}Note that $\Exc(\Gamma, \ld G)$ has an $\F_p$-structure coming from the fact that $\ld G$ is defined over $\F_p$. 


\end{example}





\begin{defn}\label{def: prime subalg}
Let $\Exc(\Gamma, \ld G)' \subset \Exc(\Gamma, \ld G)$ be the $k$-subalgebra generated by $N \cdot \Exc(\Gamma, \ld G)$ and all elements of the form $\Nm (S_{I, V, x, \xi, (\gamma_i)_{i \in I}}) = S_{I, \Nm (V), \Nm(x), \Nm(\xi), (\gamma_i)_{i \in I}} $ (the equality by Lemma \ref{lem: norm extension}). 
\end{defn}

\begin{proof}[Proof of Theorem \ref{thm: global main}]
A trivial case of Proposition \ref{prop: equiv localization for shtukas} gives an isomorphism 
\begin{equation}\label{eq: 1 trivial equiv loc isom}
T^0(\Sht_{G, D+ \infty x_0, \emptyset}; \bbm{1}) \cong H_c^0(\Sht_{H, D + \infty x_0, \emptyset}; \bbm{1}).
\end{equation}
The right side has an action of $\Exc(\Gamma, \ld H)$ through the quotient map $\Exc(\Gamma, \ld H) \surj \Exc(\Gamma, \ld H)$, and the left side has an action of $\Exc(\Gamma, \ld G)^{\sigma}$ through the map $\Exc(\Gamma, \ld G)^{\sigma} \rightarrow \Exc(\Gamma, \ld G)$. Theorem \ref{thm: global actions coincide} says that the isomorphism \eqref{eq: 1 trivial equiv loc isom} is equivariant for the action of the subalgebra $\Exc(\Gamma, \ld G)' \subset \Exc(\Gamma, \ld G)^{\sigma}$, which acts on the right side via the map 
\[
S_{I, \Nm^{(p^{-1})}(V), \Nm^{(p^{-1})}(x), \Nm^{(p^{-1})}(\xi), (\gamma_i)_{i \in I}}  \rightarrow S_{I, \Res_{\BC}(V), x, \xi, (\gamma_i)_{i \in I}}.
\]


The $L$-parameter $[\rho] \in H^1(\Gamma, \wh{H}(k))$ corresponds to a character $\chi_{\rho}  \co \Exc(\Gamma, \ld H) \rightarrow k$ under Proposition \ref{prop: galois rep}. Therefore $\chi_{\rho}$ induces a character $\chi_{\rho}' \co \Exc(\Gamma, \ld G)' \rightarrow k$ sending 
\begin{align*}
S_{I, \Nm^{(p^{-1})}(V), \Nm^{(p^{-1})}(x), \Nm^{(p^{-1})}(\xi), (\gamma_i)_{i \in I}} & \mapsto  \chi_{\rho}(S_{I, \Res_{\BC}(V), x, \xi, (\gamma_i)_{i \in I}}), \\
N \cdot S & \mapsto 0 \text{ for any $S \in \Exc(\Gamma, \ld G)$}.
\end{align*}
Let $\mf{m}_{\rho}' = \ker(\chi_{\rho}') \subset \Exc(\Gamma, \ld G)'$ be the corresponding maximal ideal. The assumption that $\chi_{\rho}'$ arises in the action of $\Exc(\Gamma, \ld G)'$ on $H^0_c(\Sht_{H,D +\infty x_0,\emptyset}; \bbm{1})_{\alpha}$ means that $H^0_c(\Sht_{H,D +\infty x_0,\emptyset}; \bbm{1})_{\alpha}$ is supported at $\mf{m}_{\rho}'$. For any $\F_p$-module $M_0$, there is semilinear action of $\Aut(k)$ on $M := M_0 \otimes_{\F_p} k$ through the second factor. This applies in particular to $\Exc(\Gamma, \ld G)'$ and $H^0_c(\Sht_{H,D +\infty x_0,\emptyset}; \bbm{1})_{\alpha}$ since they are defined over $\F_p \subset k$. Since the action of $\Exc(\Gamma, \ld G)'$ on $H^0_c(\Sht_{H,D +\infty x_0,\emptyset}; \bbm{1})_{\alpha}$ is also defined over $\F_p$, the image of $\mf{m}_{\rho}'$ under the automorphism of $\Exc(\Gamma, \ld G)'$ induced by $\Frob_p \in \Aut(k)$ also appears in the support of $H^0_c(\Sht_{H,D +\infty x_0,\emptyset}; \bbm{1})_{\alpha}$. We denote this maximal ideal by $\mf{n}_{\rho}'$; its associated character $\eta_\rho$ is characterized by the property that $\eta_\rho$ kills $N \cdot S$ for any $S \in \Exc(\Gamma, \ld G)$, and for $(V, x, \xi)$ defined over $\F_p$ it sends 
\begin{align*}
S_{I, \Nm(V), \Nm(x), \Nm(\xi), (\gamma_i)_{i \in I}} \mapsto \chi_{\rho}(S_{I, \Res_{\BC}(V), x, \xi, (\gamma_i)_{i \in I}})^p
\end{align*}
(we omitted Frobenius twists because they have no effect on maps defined over $\F_p$). 

At the start of the proof, we identified $H^0_c(\Sht_{H,D +\infty x_0,\emptyset}; \bbm{1})_{\alpha}$ with $T^0(\Sht_{G, D+ \infty x_0, \emptyset}; \bbm{1})$ as modules over $\Exc(\Gamma, \ld G)'$. The latter is a subquotient of $H_c^0(\Sht_{G, D+\infty x_0, \emptyset}; \bbm{1})$ viewed as a $\Exc(\Gamma, \ld G)'$-module via the composition $\Exc(\Gamma, \ld G)' \inj \Exc(\Gamma, \ld G)$. Since $H_c^0(\Sht_{G, D+\infty x_0, \emptyset}; \bbm{1})$ is supported at $\mf{n}_\rho' \subset \Exc(\Gamma, \ld G)'$, it is also supported at some maximal ideal $\mf{n}_{\rho}$ of $\Exc(\Gamma, \ld G)$ lying over $\mf{n}_\rho'$. Lemma \ref{lem: unique extension of char} implies that there is a \emph{unique} maximal ideal of $\Exc(\Gamma, \ld G)$ lying over $\mf{n}_\rho'$, which by Lemma \ref{lem: norm extension} corresponds to the character sending
\begin{align*}
S_{I, V, x, \xi, (\gamma_i)_{i \in I}} \mapsto & \eta_{\rho} (S_{I, \Nm(V), \Nm(x), \Nm(\xi), (\gamma_i)_{i \in I}})^{1/p} = \chi_{\rho}(S_{I, \Res_{\BC}(V), x, \xi, (\gamma_i)_{i \in I}}).
\end{align*}
This is precisely $\chi_{\rho} \circ \phi_{\BC}^*$, so its kernel is $\mf{n}_{\rho}$. We conclude that $H_c^0(\Sht_{G, D+\infty x_0, \emptyset}; \bbm{1})$ must be supported at $\ker(\chi_{\rho} \circ \phi_{\BC}^*) \subset \Exc(\Gamma, \ld G)$, as desired.

\end{proof}

\section{On local base change}\label{sec: local base change}

 In this section we will prove the main local results mentioned in the Introduction. We begin by reviewing the relevant aspects of the Genestier-Lafforgue correspondence in \S \ref{ssec: Genestier-Lafforgue}. Its key property is local-global compatibility, which will allow us to leverage the global results proved in the preceding section. 
 
 After that we embark on the construction of the map $\mf{Z}_{\TV}$ from Theorem \ref{thm: intro bernstein center}. Its definition does not require any geometry, and works equally well over local fields of characteristic zero (and residue characteristic different from $p$), but requires some technical preliminaries on Hecke algebras, which we establish in \S \ref{sssec: double coset fixed point}. Then we review the Brauer homomorphism in \S \ref{ssec: brauer homomorphism}, which is needed to finally construct $\mf{Z}_{\TV}$ and prove Theorem \ref{thm: intro bernstein center}. We then use it (and intermediate results established along the way) to prove Theorem \ref{thm: intro 2} in \S \ref{ssec: TV conjectures} and Theorem \ref{thm: local existence} in \S \ref{ssec: local base change}.

\subsection{Review of the Genestier-Lafforgue correspondence}\label{ssec: Genestier-Lafforgue}

Let $F_v$ be a local function field with ring of integers $\Cal{O}_v$ and residue characteristic $\ell \neq p$. Let $W_v$ be the Weil group of $F_v$. Let $G$ be a reductive group over $F_v$ and denote $G_v := G(F_v)$.  In \cite[Th\'{e}or\`{e}m 8.1]{GL18}, Genestier-Lafforgue construct a map 
\begin{align*}
\left\{ \begin{array}{@{}c@{}}  \text{irreducible admissible representations} \\  \text{$\pi$ of $G_v$ over $k$}\end{array} \right\}/\sim  & \longrightarrow \left\{ \begin{array}{@{}c@{}}  \text{semi-simple $L$-parameters} \\ 
\rho_{\pi} \co W_v \rightarrow \ld G(k)  \end{array} \right\}/\sim,
\end{align*}
which is characterized by local-global compatibility with Lafforgue's Global Langlands correspondence.

We briefly summarize the aspects of the Genestier-Lafforgue correspondence that we will need. 

\subsubsection{The Bernstein center}
We begin by recalling the formalism of the Bernstein center \cite{Bern84}. Let $K \subset G_v$ be an open compact subgroup with prime-to-$p$ order. The \emph{Hecke algebra of $G$ with respect to $K$ with coefficients in $\Lambda$} is
\[
\Cal{H}(G,K; \Lambda) := C_c(K \bs G_v / K; \Lambda),
\]
the compactly supported functions on $K \bs G_v / K$ valued in $\Lambda$. This forms an algebra under convolution, normalized so that the indicator function $\bbm{1}_K$ is the unit. We let $\mf{Z}(G,K; \Lambda) :=  Z(\Cal{H}(G,K; \Lambda))$ be the center of $\Cal{H}(G,K; \Lambda)$. 

If $K \subset K'$ have prime-to-$p$ pro-order (e.g., this will be true as long as they are sufficiently small), then convolution with $\bbm{1}_{K'}$ gives a homomorphism $\mf{Z}(G,K; \Lambda) \rightarrow \mf{Z}(G,K'; \Lambda)$. The \emph{Bernstein center of $G$ (with coefficients in $\Lambda$)} is 
\[
\mf{Z}(G; \Lambda) := \varprojlim_{K} \mf{Z}(G,K; \Lambda),
\]
where the transition maps to $\mf{Z}(G,K; \Lambda)$ are as above, and the inverse limit is taken over $K$ with prime-to-$p$ pro-order. 

If $\Lambda = k$, we abbreviate $\Cal{H}(G,K) := \Cal{H}(G,K; k)$, $\mf{Z}(G,K) := \mf{Z}(G,K; k)$, and $\mf{Z}(G) := \mf{Z}(G; k)$.

The ring $\mf{Z}(G)$ has another interpretation as the ring of endomorphisms of the identity functor of the category of smooth $k$-representations of $G_v$. Explicitly, smoothness of $\Pi$ implies that $\Pi  = \bigcup_{\text{open compact }K \subset G_v} \Pi^{K}$, and $\mf{Z}(G,K)$ acts on $\Pi^{K}$ as an $\Cal{H}(G,K)$-module; this assembles into action of $\mf{Z}(G)$ on $\Pi$. In particular, any irreducible admissible representation $\Pi$ of $G_v$ over $k$ induces a character of $\mf{Z}(G)$. 

\subsubsection{Action of the excursion algebra} 
The main result of \cite{GL18} is the construction of a homomorphism
\begin{equation}\label{eq: Z_G}
Z_G  \co \Exc(W_v, \ld G) \rightarrow \mf{Z}(G).
\end{equation}

Let $x \in \cB(G/F_v)$ be a point of the Bruhat-Tits building of $G$, and use it to extend $G$ to a parahoric group scheme over $\Cal{O}_v$. (Some reminders on Bruhat-Tits theory will appear in \S \ref{sssec: BT theory}.) For $r \geq 0$, let $K_{r} := G(F_v)_{x,r}$; this is an open compact subgroup of $F_v$. We write $Z_{G,r} \co \Exc(W_v, \ld {G}) \rightarrow \mf{Z}(G, K_{r})$ for the composition of $Z_G$ with the tautological projection to $\Cal{H}(G, K_{r})$.

We will shortly give a characterization of \eqref{eq: Z_G}. First let us indicate how \eqref{eq: Z_G} defines the correspondence $\Pi \mapsto \rho_{\Pi}$. An irreducible admissible $\Pi$ induces a character of $\mf{Z}(G)$, as discussed above. Composing with $Z_G$ then gives a character of $\Exc(W_v, \ld {G}) $, which by Proposition \ref{prop: galois rep} gives a semisimple Langlands parameter $\rho_{\Pi} \in H^1(W_v, \wh{G}(k))$.

\subsubsection{Local-global compatibility}  Choose a smooth projective curve $X$ over $\F_{\ell}$ and a point $v \in X$ so that $X_v = \Spec \Cal{O}_v$, and $G$ extends to a reductive group over the generic point of $X$. Then choose a further extension of $G$ to a parahoric group scheme over all of $X$, such that $G/\cO_v$ is the parahoric group scheme corresponding to the chosen point $x \in \cB(G/F_v)$. 

Choose an embedding of the local Weil group $W_v$ into the global Weil group of $X$, which we denote $\Weil(\eta, \ol{\eta})$. The map \eqref{eq: Z_G} is characterized by local-global compatibility with the global excursion action. The idea is that for $(\gamma_i)_{i \in I} \subset W_v^I \subset \Weil(\eta, \ol{\eta})^I$, the action of the the excursion operator $S_{I, f, (\gamma_i)_{i \in I}}$ on $H_c^0(\Sht_{G, D, \emptyset}; \bbm{1})$ is local at $v$, i.e. it acts through a Hecke operator for $G_v$. Moreover, it commutes with other Hecke operators because all excursion operators commute with all Hecke operators, hence it must actually be in the center of the relevant Hecke algebra. This idea is affirmed by the Proposition below.

\begin{prop}[Genestier-Lafforgue Proposition 1.3]\label{prop: GL 1.3}  
Let $r \geq 0$ be an integer and $D := rv + D^v$ a closed finite subscheme of $X$, with $D^v$ supported away from $v$. For $(\gamma_i)_{i \in I} \subset W_v^I$, the operator $S_{I, f, (\gamma_i)_{i \in I}}$ acts on $H_c^0(\Sht_{G, D, \emptyset}; \bbm{1})$ as convolution by $Z_{G,r}(S_{I, f, (\gamma_i)_{i \in I}}) \in \mf{Z}(G, K_{r})$. 
\end{prop}

\begin{remark}\label{rem: end injective}
By \cite[Lemme 1.4]{GL18}, for large enough (depending on $r$) $D^v$ the action of $Z_{G,r}(S_{I, f, (\gamma_i)_{i \in I}})$ on $H_c^0(\Sht_{G, D, \emptyset}; \bbm{1})$ is faithful. Therefore, Proposition \ref{prop: GL 1.3} certainly characterizes the map \eqref{eq: Z_G}. What is not clear is that the resulting $Z_{G,r}(S_{I, f, (\gamma_i)_{i \in I}})$ is independent of choices (of the global curve, or the integral model of the affine group scheme). In \cite{GL18} this is established by giving a purely local construction of \eqref{eq: Z_G} in terms of ``restricted shtukas'', but for our purposes it will be enough to accept Proposition \ref{prop: GL 1.3} as a black box. 
\end{remark}

\subsection{Preliminary results on Hecke algebras}\label{sssec: double coset fixed point} We next establish some technical lemmas which aid to study the properties of the Brauer homomorphism. The only result that will be needed in later subsections is Corollary \ref{cor: surjectivity}.


\subsubsection{Assumptions}\label{sssec: HA assumptions} In this subsection, we allow $F_v$ to be any local field (including one of characteristic zero) of residue characteristic $\ell \neq p$. Let $E_v/F_v$ be a finite Galois assumption such that $\Gal(E_v/F_v)$ has order coprime to $\ell$. We let $H$ be any (connected) reductive group over $F_v$ and $G := \Res_{E_v/F_v}(H_{E_v})$. We abbreviate $H_v = H(F_v)$ and $G_v = G(F_v) = H(E_v)$. 

\subsubsection{Reminders on Bruhat-Tits theory}\label{sssec: BT theory} First we recall some relevant facts from Bruhat-Tits theory, originally developed in \cite{BT72} but explained in the form used below in \cite{KP}.

Let $\cB(H/F_v)$ be the Bruhat-Tits building of $H/F_v$ and $x \in \cB(H/F_v)$. Associated to $x$ there is the parahoric group $H(F_v)_x^0 := H(F_v)_{x,0}$, along with its decreasing filtration $H(F_v)_{x,r}$ for $r \geq 0$. The subgroup $H(F_v)_{x,0+} := \bigcup_{s > 0} H(F_v)_{x,s}$ is pro-$\ell$. 

We record some descent properties: if $E_v/F_v$ is unramified then $H(E_v)_{x,r}^{\Gal(E_v/F_v)} = H(F_v)_{x,r}$, and if $E_v/F_v$ is tamely ramified then $H(E_v)_{x,r}^{\Gal(E_v/F_v)} \supset H(F_v)_{x,r}$ as explained in \cite[\S 12]{KP}. 

We shall make use of the \emph{Cartan decomposition}. We follow the description in \cite[\S 5.2]{KP}. Let $S$ be a maximal $F_v$-split torus of $H$, and $Z = Z_H(S)$. Referring to \cite[(5.2.1)]{KP} for undefined notation, we have a subset 
\begin{equation}\label{eq: cartan fund domain}
\cZ := \{ z \in Z(F_v) \co \alpha(\omega_Z(z)) \geq 0 \text{ for all } \alpha \in \Phi^+\} \subset Z(F_v).
\end{equation}
According to \cite[Theorem 5.2.1]{KP}, for a \emph{special vertex} $x$ in the apartment of $S$, we have 
\begin{enumerate}
\item $H(F_v) = H(F_v)_{x,0} \cdot \cZ \cdot H(F_v)_{x,0}$, and 
\item for $z, z' \in \cZ$, $H(F_v)_{x,0} \,  z \, H(F_v)_{x,0}   = H(F_v)_x^0 \,  z'  \,  H(F_v)_x^0$ if and only if $z'(z^{-1}) \in Z(F_v)^0$,
\end{enumerate}
where for a Levi subgroup $M \subset G$, the subgroup $M(F_v)^0 \subset M(F_v)$ is the kernel of the Kottwitz homomorphism $\kappa_M \co M(F_v) \rightarrow \pi_1(M)_I^{\Gal(\breve{F}_v/F_v)}$; it may be defined alternatively as in \cite[Definition 2.6.23]{KP}. We have $Z(F_v)^0 \subset H(F_v)_{x,0}$.

\subsubsection{Maps of double coset spaces}\label{sssec: double cosets} Following the notation of \S \ref{sssec: HA assumptions}, let $S_G$ be a maximal $F_v$-split torus of $G$, and define $Z_E := Z_G(S_G)$, and $\cZ_E \subset Z_E(F_v)$ as in \eqref{eq: cartan fund domain}. Let $S_F$ be a maximal $F_v$-split torus of $H$ contained in $S_G$, $Z_F = Z_H(S_H)$, and $\cZ_F \subset Z_F(F_v)$ as in \eqref{eq: cartan fund domain}. Let $x \in \cB(G/F_v)$ be a special vertex in the apartment of $S_G$. We will abbreviate $K_r := G(F_v)_{x,r}$ and $U_r = K_r^{\Gal(E_v/F_v)} \supset H(F_v)_{x,r}$. The goal of this subsection is to prove the following. 

\begin{prop}\label{prop: inj}
If $r \geq 0$, then the map $U_r \bs H_v  / U_r \rightarrow K_r \bs G_v  / K_r$ is injective. 
\end{prop}

\begin{proof}
We first handle the case $r=0$, which will seen to be a consequence of Cartan decomposition. Since $U_0 \supset H(F_v)_{x,0}$, the Cartan decomposition implies that double cosets $U_0 \bs G_v / U_0$ are represented by $z \in \cZ_F$. If $z_1, z_2 \in \cZ_F$ are such that $K_0 z_1 K_0 = K_0 z_2 K_0$, then $z_2z_1^{-1} \in Z_E(F_v)^0 \subset K_0$. On the other hand, clearly $z_2z_1^{-1} \in H(F_v)$, so we conclude that $z_2z_1^{-1} \in K_0 \cap H_v = U_0$. Therefore, $U_0 z_1 U_0 = U_0 z_2 U_0$, and the case $r=0$ is concluded. 

Now suppose $r>0$. Let $a,b \in U_r \bs H_v$ be two elements whose images in $K_r \bs G_v$ lie in the same orbit for the right translation of $K_r$. In other words, $a = bk$ for some $k \in K_r$. Since $a,b$ are fixed by $\sigma$, this implies that 
\[
a = b \sigma(k)
\]
and therefore $\sigma(k) k^{-1} \in \Stab_{K_r}(b) =: K^b_r$. Note that $\Gal(E_r/F_r)$ is of order prime-to-$\ell$ while $K_r$ is pro-$\ell$ thanks to the assumption $r>0$, so then $H^1(\Gal(E_r/F_r); K^b_r) = 0$. This means that there exists $y \in K^b_r$ such that $\sigma(k) k^{-1} = \sigma(y) y^{-1}$. Then $y^{-1} k$ is fixed by $\sigma$, so $y^{-1} k \in H_v \cap K_r = U_r$. But then 
\[
 a = bk = (by^{-1}) k  = b (y^{-1} k),
 \]
which shows that $a$ and $b$ lie in the same orbit for the right translation of $U_r$ on $U_r \bs H_v$. 
\end{proof}

In the following Corollary, we let $\Gal(E_v/F_v)$ act on $G_v = H(E_v)$ by the natural Galois action, which induces an action on $\cH(G, K_r; \Lambda)$. 

\begin{cor}\label{cor: surjectivity}
Suppose $r\geq 0$. Then the restriction map $\cH(G, K_r; \Lambda)^{\Gal(E_v/F_v)} \rightarrow \cH(H, U_r; \Lambda)$ is surjective. 
\end{cor}

\subsection{The Brauer homomorphism}\label{ssec: brauer homomorphism}

We introduce the notion of the Brauer homomorphism from \cite{TV}, whose utility for our purpose is to capture the relationship between $\Pi$ and its Tate cohomology from the perspective of Hecke algebras. 

\subsubsection{Assumptions}\label{sssec: assumptions Brauer} In this subsection we allow $F_v$ to be any local field (including one of characteristic zero) of residue characteristic $\ell \neq p$. We assume that $\Gal(E_v/F_v)$ is cyclic of order $p$, and we let $\sigma \in \Gal(E_v/F_v)$ be a generator. We let $H$ be any (connected) reductive group over $F_v$ and $G := \Res_{E_v/F_v}(H_{E_v})$. Subgroups $K_r \subset G_v$, $U_r \subset H_v$ are defined as in \S \ref{sssec: double cosets}.

\subsubsection{The (un-normalized) Brauer homomorphism}
Let $K \subset G_v$ be an open compact subgroup, and let $U := K^{\sigma} \subset H_v$. We say that $K \subset G_v$ is a \emph{plain subgroup} if $(G_v/K)^{\sigma} = H_v/U$. 

We can view $\Cal{H}(G, K)$ as the ring of $G_v$-invariant (for the diagonal action) functions on $(G_v / K) \times (G_v / K)$ under convolution. 

\begin{lemma}\label{lem: Br homomorphism}
If $K \subset G_v$ is a plain subgroup, then the restriction map 
\begin{align}\label{eq: Br}
\Cal{H}(G, K)^{\sigma}  &= \mrm{Fun}_{G_v,c}((G_v / K) \times (G_v / K), k)^{\sigma} \\
&  \xrightarrow{\text{restrict}} \mrm{Fun}_{H_v,c}((H_v/U) \times (H_v/U), k)  = \Cal{H}(H_v, U) \nonumber
\end{align}
is an algebra homomorphism.
\end{lemma}

\begin{proof}
What we must verify is that for $x,z \in H_v/U$, and $f,g \in\Cal{H}(G, K)^{\sigma}$, we have 
\begin{equation}\label{eq: rest alg}
\sum_{y \in G_v/K} f(x, y) g(y,z) = \sum_{y \in H_v/U} f(x,y) g(y,z).
\end{equation}
Since $f$ and $g$ are $\sigma$-invariant, we have 
\[
f(x,y) = f(\sigma x, \sigma y ) = f(x, \sigma y) \text{ and } g(y,z ) = g(\sigma y, \sigma z) = g(\sigma y, z).
\]
If $y \notin H_v/U$, then the plain-ness assumption implies that $y$ is not fixed by $\sigma$. Therefore the contribution from the orbit of $\sigma$ on $y$ to \eqref{eq: rest alg} is a multiple of $p$, which is $0$ in $k$. 
\end{proof}

The map of Lemma \ref{lem: Br homomorphism} was introduced in \cite[\S 4]{TV} and called the \emph{(un-normalized) Brauer homomorphism}. We denote it 
\[
\Br \co \Cal{H}(G, K)^{\sigma} \rightarrow \Cal{H}(H, U).
\]

\begin{lemma}\label{lem: deep lattices are plain}
If $K \subset G(F_v)_{x,0+}$ for any $x \in \cB(G/F_v)$, then $K$ is plain.
\end{lemma}

\begin{proof}
By the long exact sequence for group cohomology, the plain-ness of $K \subset G_v$ is equivalent to condition that the map on non-abelian cohomology $H^1( \langle \sigma \rangle; K) \rightarrow H^1(\langle \sigma \rangle; G_v)$ has trivial fiber over the trivial class. But since $G(F_v)_{x,0+}$ is pro-$\ell$, all its subgroups are acyclic for $H^1( \langle \sigma \rangle, -)$ as $\sigma$ has order $p$. Therefore $H^1( \langle \sigma \rangle, K)$ vanishes for all such $K \subset G(F_v)_{x,0+}$.
\end{proof}

\begin{lemma}[Relation to the Brauer homomorphism]\label{lem: brauer homomorphism}
Assume $K \subset G_v$ is plain. Suppose $\Pi$ is a $\sigma$-fixed representation of $G_v$. Then the map of Tate cohomology groups $T^*(\Pi^{K}) \rightarrow	 T^*(\Pi)$ lands in the $U$-invariants, and for any $h \in \Cal{H}(G, K)^{\sigma}$ we have the commutative diagram below.
\[
\begin{tikzcd}
T^*(\Pi^{K}) \ar[d, "T^0h"] \ar[r] & T^*(\Pi)^{U} \ar[d, "\Br(h)"]  \\
T^*(\Pi^{K})\ar[r] & T^*(\Pi)^{U}  
\end{tikzcd}
\]
(Here $T^0 h$ is the element of $T^0 (\Cal{H}(G, K))$ represented by $h$.) 
\end{lemma}

\begin{proof} This is \cite[\S 6.2]{TV}; it follows from a direct computation similar to the proof of Lemma \ref{lem: Br homomorphism}.
\end{proof}

\subsubsection{Treumann-Venkatesh homomorphism}
If we take $K = K_r$ as in Corollary \ref{cor: surjectivity}, then the Brauer homomorphism $
\Br \co \Cal{H}(G, K_r)^{\sigma} \rightarrow \Cal{H}(H, U_r)$ is a surjective algebra homomorphism, hence induces a map on centers 
\begin{equation}\label{eq: un-normalized center}
Z(\Br) \co Z(\Cal{H}(G, K_r)^{\sigma}) \rightarrow \mf{Z}(H, U_r).
\end{equation}
It is evident from the definition that $Z(\Br)$ through the quotient $Z(\Cal{H}(G, K_r)^{\sigma})  / N \cdot \mf{Z}(G, K_r)$.

Since $\mf{Z}(G, K_r)$ is commutative, it has a Tate diagonal homomorphism (Definition \ref{def: tate diagonal}) $\mf{Z}(G, K_r) \xrightarrow{\Delta^p} T^0(\mf{Z}(G,K_r))$. Since $Z(\Cal{H}(G, K_r))^{\sigma}  \subset Z(\Cal{H}(G, K_r)^{\sigma}) $, we may compose with $Z(\Br)$ to obtain a map  
\begin{equation}\label{eq: Z(Br) tilde}
\mf{Z}(G, K_r) \xrightarrow{\Delta^p} T^0(\mf{Z}(G,K_r)) \xrightarrow{Z(\Br)} \mf{Z}(H, U_r).
\end{equation}
Note however that it is not $k$-linear, since $\Delta^p$ is $\Frob$-semilinear over $k$. Then there is a (unique) homomorphism $Z'$ fitting into the commutative diagram 
\begin{equation}\label{eq: linearize Brauer}
\begin{tikzcd}
\mf{Z}(G, K_r) \ar[d, "{\sim}"] \ar[dr, "Z(\Br) \circ \Delta^p"] \\
 \mf{Z}(G, K_r)^{(p)} \ar[r, "Z'"', dashed] &   \mf{Z}(H, U_r)
 \end{tikzcd}
\end{equation}
We have $\mf{Z}(G, K_r) \cong \mf{Z}(G, K_r; \F_p) \otimes_{\F_p} k$, which as explained in Construction \ref{const: frob twist algebra} induces a $k$-linear isomorphism 
\begin{equation}\label{eq: Hecke F_p structure}
\mf{Z}(G, K_r)^{(p)} \cong \mf{Z}(G, K_r).
\end{equation}

\begin{defn}\label{def: TV homomorphism}
We define \emph{Treumann-Venkatesh homomorphism} $\mf{Z}_{\TV,r}$ to be the homomorphism $ \mf{Z}(G, K_r) \rightarrow \mf{Z}(H, U_r)$ obtained by linearization of $Z'$ in the sense of Construction \ref{const: frob twist algebra} with respect to the $\F_p$ structure \eqref{eq: Hecke F_p structure}. 
\end{defn}

\begin{remark}
Definition \ref{def: TV homomorphism} is not considered in \cite{TV}, but it is inspired by the definition of the \emph{normalized Brauer homomorphism} in \cite[\S 4.3]{TV}. 
\end{remark}

\subsection{The base change homomorphism for Bernstein centers}\label{ssec: BC homomorphism for BC}

For now, assumptions are as in \S \ref{sssec: assumptions Brauer}. Suppose $K_r$ has prime-to-$p$ pro-order. For $s>r$, so that $K_s \subset K_r$, we have a map $e_G^{s \rightarrow r} \co \mf{Z}(G, K_s) \rightarrow \mf{Z}(G, K_r)$ given by convolution with $\bbm{1}_{K_r}$. (Technically $e_G^{s \rightarrow r}$ also depends on the point $x \in \cB(G/F_v)$ used to define the $K_r$, but we suppress this from our notation.) Similarly, we have $e_H^{s \rightarrow r} \co \mf{Z}(H, U_s) \rightarrow \mf{Z}(H, U_r)$ given by convolution with $\bbm{1}_{U_r}$. 

\begin{lemma}\label{lem: TV compatible}
The diagram 
\[
\begin{tikzcd}
\mf{Z}(G, K_s) \ar[r, "\mf{Z}_{\TV,s}"]  \ar[d, "e_G^{s \rightarrow r}"] & \mf{Z}(H, U_s) \ar[d, "e_H^{s \rightarrow r}"]  \\
\mf{Z}(G, K_r) \ar[r, "\mf{Z}_{\TV,r}"] & \mf{Z}(H, U_r) 
\end{tikzcd}
\]
commutes. 
\end{lemma}

\begin{proof}
This follows by direct computation, using Lemma \ref{lem: Br homomorphism} and Lemma \ref{lem: brauer homomorphism}. 
\end{proof}

\begin{defn}[Base change homomorphism for Bernstein centers]\label{def: BC for BC} We define the map $\mf{Z}_{\TV} \co \mf{Z}(G) \rightarrow \mf{Z}(H)$ as
\[
 \varprojlim_r \mf{Z}_{\TV,r} \co \varprojlim_r \mf{Z}(G, K_r) \rightarrow \varprojlim_r \mf{Z}(H, U_r) .
\]
\end{defn}

Definition \ref{def: BC for BC} is well-defined over local fields of any residue characteristic $\ell \neq p$, but in this paper we will only prove properties of it for local function fields. Hence, \textbf{for the rest of the paper, we assume that $F_v$ is a local field of \emph{positive} characteristic}	. The rest of this subsection shall be devoted to the proof of Theorem \ref{thm: intro bernstein center}.

\subsubsection{}\label{sssec: local step1} The maps
\[
\Exc(W_v, {}^L G)  \xrightarrow{Z_{G,r}}  \mf{Z}(G, K_{r}) \rightarrow \End_{\Cal{H}_G}(H_c^0(\Sht_{G,D,\emptyset}; \bbm{1}))
\]
induce upon applying Tate cohomology, 
\[
T^0\Exc(W_v, {}^L G)   \xrightarrow{T^0 Z_{G,r}} T^0 \mf{Z}(G, K_{r}) \rightarrow \End_{T^0\Cal{H}_G}(T^0 (H_c^0(\Sht_{G,D,\emptyset}; \bbm{1}))).
\]
Fix a closed point $x_0 $ on $X$ distinct from $v$. For each integer $r$, we will impose level structure along $D := r v + \infty x_0$, interpreted as in \S \ref{sssec: infinite level}. By Remark \ref{rem: end injective}, the map $\mf{Z}(G, K_{r}) \rightarrow \End_{\Cal{H}_G}(H_c^0(\Sht_{G,D,\emptyset}; \bbm{1}))$ is injective.

\subsubsection{}\label{sssec: local step2} Theorem \ref{thm: global actions coincide} implies that under the identification $T^0 (\Sht_{G, D,\emptyset}; \bbm{1}) \cong T^0 (\Sht_{H,D,\emptyset}; \bbm{1})$, we have 
\[
\left(\begin{array}{@{}c@{}} 
\text{the action of} \\ \text{$S_{I, \Nm^{(p^{-1})}(V), \Nm^{(p^{-1})}(x), \Nm^{(p^{-1})}(\xi), (\gamma_i)_{i \in I}}$} \\  \text{on $T^0 (\Sht_{G, D,\emptyset}; \bbm{1})$} 
\end{array} \right)
 = 
\left( \begin{array}{@{}c@{}} \text{the action of $S_{I, \Res_{\BC}(V), x, \xi, (\gamma_i)_{i \in I}}$}\\ \text{on $T^0 (\Sht_{H,D,\emptyset}; \bbm{1})$}
 \end{array}\right).
 \]

\subsubsection{}\label{sssec: local step3} For any set $S$, we let $k[S]$ denote the $k$-vector space of $k$-valued functions on $S$. 

Now suppose $\wt{S}$ is a set with an action of $G_v \rtimes \langle \sigma \rangle$, on which an open compact subgroup $K \subset G_v$ acts freely. Then for $S := \wt{S}/K$, there is a natural action of $\Cal{H}(G,K)$ on $k[S]$ since we may view $\Cal{H}(G, K) = \Hom_{G_v}(k[G_v/K], k[G_v/K])$ and $k[S] = \Hom_{G_v}(k[G_v/K], k[\wt{S}])$. This induces an action of $T^0 (\Cal{H}(G, K))$ on $T^0(k[S]) \cong k[S^{\sigma}]$, and then by inflation an action of $\Cal{H}(G,K)^{\sigma}$ on $k[S^{\sigma}]$. 

By the same mechanism, for $U := K^\sigma$ there is an induced action of $\Cal{H}(H, U)$ on $k[\wt{S}^{\sigma}/K^{\sigma}] = k[\wt{S}^{\sigma}/U]$. We define $Z_{H,r} \co \Exc(W_v, \ld {H}) \rightarrow \mf{Z}(H, U_{r})$ similarly to $Z_{G,r}$. 

\begin{lemma}\label{lem: Br action}
Assume $K \subset G_v$ is a plain subgroup. Then $k[\wt{S}^{\sigma}/U]$ is a $\Cal{H}(G,K)^{\sigma}$-direct summand of $k[S^{\sigma}]$, and for all $h \in \Cal{H}(G,K)^{\sigma}$ we have
\[
\left(\begin{array}{@{}c@{}} 
\text{the action of $h$ on $k[\wt{S}^{\sigma}/U]$}
\end{array} \right)
 = 
\left( \begin{array}{@{}c@{}} \text{the action of $\Br(h) \in \Cal{H}(H,U)$ on $k[\wt{S}^{\sigma}/U]$}
 \end{array}\right).
 \] 
\end{lemma}

\begin{proof}
See \cite[equation (4.2.2)]{TV}. 
\end{proof}

From \S \ref{sssec: local step1} we have the diagram 
\begin{equation}\label{diag: tate actions}
\begin{tikzcd}
T^0\Exc(W_v, {}^L G)  \ar[r, "Z_{G,r}"] & T^0 \mf{Z}(G, K_{r}) \ar[r] \ar[d, "Z(\Br)"]  & \End_{T^0 \Cal{H}_G} (T^0 (\Sht_{G, D, \emptyset}; \bbm{1}))  \ar[d] \\
\Exc(W_v, \ld H) \ar[r, "Z_{H,r}"] &  \mf{Z}(H, U_{r}) \ar[r, hook] & \End_{\Cal{H}_H} (T^0(\Sht_{H, D, \emptyset}; \bbm{1}))
\end{tikzcd}
\end{equation}
Here the right vertical map is the identity map on endomorphisms, with respect to the identification $T^0 (\Sht_{G, D, \emptyset}; \bbm{1}) \cong 
T^0(\Sht_{H, D, \emptyset}; \bbm{1})$. Note that the surjectivity statement of Corollary \ref{cor: surjectivity}, plus Lemma \ref{lem: Br homomorphism} and Lemma \ref{lem: brauer homomorphism} giving compatibility of the respective actions with the Brauer homomorphism, are what guarantees that an endomorphism commuting with the $T^0(\cH_G)$-action also commutes with the $\cH_H$-action.

\begin{cor}\label{cor: actions agree} For all $r \geq 1$, the action of $z \in T^0\mf{Z}(G, K_{r}) $ on $T^0(\Sht_{G, D, \emptyset}; \bbm{1}) $ in \eqref{diag: tate actions} agrees with the action of $Z(\Br)(z)$ on $T^0(\Sht_{H, D, \emptyset}; \bbm{1})$ in \eqref{diag: tate actions} under the identification $T^0(\Sht_{G, D, \emptyset}; \bbm{1}) \cong T^0(\Sht_{H, D, \emptyset}; \bbm{1})$ from \S \ref{ssec: general localization}. In other words, the square in diagram \eqref{diag: tate actions} commutes.
\end{cor}

\begin{proof}
Apply Lemma \ref{lem: Br action} with $S := \Sht_{G, D, \emptyset}$ and $\wt{S} := \Sht_{G, \infty v + \infty x_0, \emptyset} := \varprojlim_{j \geq 0} \Sht_{G, (r+j) v + j x_0, \emptyset}$. Then $k[S]$ is identified with the functions on $\Sht_{G, rv+ \infty x_0, \emptyset}$, and Lemma \ref{lemma: shtuka fixed point} plus \S \ref{ssec: general localization} identify $k[\wt{S}^{\sigma}/K^{\sigma} ]$ with the functions on $\Sht_{H, rv+\infty x_0, \emptyset}$. 

As compactly supported functions are dual to functions, the assertions for compactly supported functions then follow by duality. 
\end{proof}

\begin{cor}\label{cor: local operators identified} 
For all $r \geq 1$, for all $\{V, x, \xi, (\gamma_i)_{i \in I}\}$ as in \S \ref{ssec: excursion presentation 2}, $\Br$ sends 
\begin{align*}
 Z_{G,r}(S_{I, \Nm^{(p^{-1})}(V), \Nm^{(p^{-1})}(x), \Nm^{(p^{-1})}(\xi), (\gamma_i)_{i \in I}}) & \in \mf{Z}(G, K_{r}) \subset \Cal{H}(G, K_{r}) \\
 \hspace{1cm} \xrightarrow{\Br} Z_{H,r}(S_{I, \Res_{\BC}(V), x, \xi, (\gamma_i)_{i \in I}}) & \in \mf{Z}(H,U_{r}) \subset \Cal{H}(H, U_{r}).
\end{align*}
\end{cor}

\begin{proof}
The equality from \S \ref{sssec: local step2} shows that 
\begin{equation}\label{eq: excursion and brauer}
\left(\begin{array}{@{}c@{}} 
\text{the image of }\\
\underbrace{S_{I, \Nm^{(p^{-1})}(V), \Nm^{(p^{-1})}(x), \Nm^{(p^{-1})}(\xi), (\gamma_i)_{i \in I}}}_{\in T^0\Exc(W_v, \ld G)}  \\ \text{in $\End_{\Cal{H}_H} (T^0 (\Sht_{H, rv+\infty x_0, \emptyset}; \bbm{1}))$} \\
\text{via \eqref{diag: tate actions}}
\end{array} \right)
 = 
\left(\begin{array}{@{}c@{}} 
\text{the image of } \\
\underbrace{S_{I, \Res_{\BC}(V), x, \xi, (\gamma_i)_{i \in I}}}_{\in \Exc(W_v, \ld H)}  \\ \text{in $\End_{\Cal{H}_H} (T^0 (\Sht_{H, rv+\infty x_0, \emptyset}; \bbm{1}))$}  \\
\text{via \eqref{diag: tate actions}}
\end{array} \right).
 \end{equation}
On the other hand, Corollary \ref{cor: actions agree} shows that the left hand side of \eqref{eq: excursion and brauer} agrees with the image of \\ $\Br(Z_{G,r}(S_{I, \Nm^{(p^{-1})}(V), \Nm^{(p^{-1})}(x), \Nm^{(p^{-1})}(\xi), (\gamma_i)_{i \in I}}))$ via \eqref{diag: tate actions}, for all $r \geq 1$. We conclude using that the map $\mf{Z}(H, U_{r}) \hookrightarrow \End_{\Cal{H}_H} (T^*(\Sht_{H, D, \emptyset}; \bbm{1}))$ in \eqref{diag: tate actions} is injective, which follows from Remark \ref{rem: end injective}.
\end{proof}

\subsubsection{} Recall that in Definition \ref{defn: phi_{BC}} we have defined a map $\phi_{\BC} \co \ld H \rightarrow \ld G$ over $k$.

\begin{cor} The following diagram commutes:
\begin{equation}\label{eq: base change diagram}
\begin{tikzcd}
\Exc(W_v, \ld G)  \ar[r, "\phi_{\BC}^*"] \ar[d, "Z_G"]  & \Exc(W_{v}, \ld H)  \ar[d, "Z_H"]\\
\mf{Z}(G)\ar[r, "\mf{Z}_{\TV}"] & \mf{Z}(H)
\end{tikzcd}
\end{equation}
\end{cor}

\begin{proof}The commutativity of the diagram 
\[
\begin{tikzcd}
\Exc(W_v, \ld G) \ar[r, "\Delta^p"] \ar[d, "Z_G"] & T^0 \Exc(W_v, \ld G) \ar[d, "T^0 (Z_G)"]   \\
\mf{Z}(G)  \ar[r, "\Delta^p"] & T^0 \mf{Z}(G)
\end{tikzcd}
\]
implies that $Z(\Br) \circ \Delta^p \circ Z_G  = Z(\Br) \circ T^0(Z_G) \circ \Delta^p $. By definition $\mf{Z}_{\TV} \circ Z_G$ is the linearization of $Z(\Br) \circ \Delta^p \circ Z_G$, so it is also the linearization of $Z(\Br) \circ T^0(Z_G) \circ \Delta^p$.

By Lemma \ref{lem: norm extension}, the Tate diagonal $\Delta^p \co \Exc(W_v, \ld G) \rightarrow T^0(\Exc(W_v, \ld G))$ sends 
\begin{equation}\label{eq: 6 eq revise 1}
S_{I,V, x, \xi, (\gamma_i)_{i \in I}} \xrightarrow{\Delta^p} S_{I, \Nm(V), \Nm(x), \Nm(\xi), (\gamma_i)_{i \in I}}.
\end{equation}
It linearization therefore sends 
\[
S_{I, V, x, \xi, (\gamma_i)_{i \in I}} \mapsto S_{I, \Nm^{(p^{-1})}(V), \Nm^{(p^{-1})}(x), \Nm^{(p^{-1})}(\xi), (\gamma_i)_{i \in I}}
\]
since this is a $k$-algebra homomorphism that agrees with \eqref{eq: 6 eq revise 1} when $(V, x, \xi)$ are defined over $\F_p$. Applying Corollary \ref{cor: local operators identified} with $r \rightarrow \infty$, we have 
\begin{align*}
Z(\Br) \circ Z_{G}(S_{I, \Nm^{(p^{-1})}(V), \Nm^{(p^{-1})}(x), \Nm^{(p^{-1})}(\xi), (\gamma_i)_{i \in I}}) &=  Z_H(S_{I, \Res_{\BC}(V), x, \xi, (\gamma_i)_{i \in I}}) \\
& = Z_H( \phi_{\BC}^* (S_{I, V, x,  \xi , (\gamma_i)_{i \in I}})).
\end{align*}
Therefore the linearization of $Z(\Br) \circ T^0(Z_G) \circ \Delta^p$ agrees with $Z_H \circ \phi_{\BC}^*$. 
\end{proof}


\begin{proof}[Completion of the proof of Theorem \ref{thm: intro bernstein center}]
Let $\pi$ be an irreducible representation of $H_v$ and $\chi_{\pi} \co \mf{Z}(H) \rightarrow k$ the induced character. By the definition of the Genestier-Lafforgue parametrization, $\rho_{\pi}$ corresponds to $\chi_{\pi} \circ Z_H$ via Proposition \ref{prop: galois rep}. Then \eqref{eq: base change diagram} implies that $\chi_{\pi} \circ   \mf{Z}_{\TV} \circ Z_G  =  \chi_{\pi} \circ Z_H \circ \phi_{\BC}^*$ is associated to the $L$-parameter $\phi_{\BC} \circ \rho_{\pi}$. 
\end{proof}

\subsection{The Treumann-Venkatesh Conjecture}\label{ssec: TV conjectures} In this subsection we will prove Theorem \ref{thm: intro 2}. We begin by formulating the Treumann-Venkatesh Conjecture precisely in this setting. (The original phrasing of \cite{TV} is in terms of a hypothetical Local Langlands correspondence which was not defined at the time for general groups.) 

\subsubsection{Assumptions} In this subsection the assumptions are as in \S \ref{ssec: brauer homomorphism}, and we furthermore assume $F_v$ is a local function field. We note, however, that the formulation of all the statements in \S \ref{def: BC for BC} makes sense for any local field $F_v$ of residue field $\ell \neq p$, with a suitable replacement for the Genestier-Lafforgue correspondence, and that all our arguments in this subsection apply if those statements are true for $F_v$.

\subsubsection{Formulation of the Conjecture} Let $\Pi$ be an irreducible admissible representation of $G_v$ over $k$. Let $\Pi^{\sigma}$ be the representation of $G_v$ obtained by composing $\Pi$ with $\sigma \co G_v \rightarrow G_v$. We say that $\Pi$ is \emph{$\sigma$-fixed} if $\Pi \approx \Pi^{\sigma}$ as $G_v$-representations. 

\begin{lemma}[{\cite[Proposition 6.1]{TV}}]\label{lem: 6.1} If $\Pi$ is $\sigma$-fixed, then the $G_v$-action on $\Pi$ extends uniquely to an action of $G_v \rtimes \langle \sigma \rangle$. 
\end{lemma}


Using Lemma \ref{lem: 6.1} we can form the Tate cohomology groups $T^0(\Pi)$ and $T^1(\Pi)$ with respect to the $\sigma$-action, which are then representations of $H_v$. Treumann-Venkatesh conjecture that they are in fact admissible representations of $H_v$, but we do not prove or use this. 

\begin{defn}[Linkage] An irreducible admissible representation $\pi$ of $H_v$ is \emph{linked} with an irreducible admissible representation $\Pi$ of $G_v$ if $\pi^{(p)}$ appears in $T^0(\Pi)$ or $T^1(\Pi)$, where $\pi^{(p)}$ is the Frobenius twist
\[
\pi^{(p)} := \pi \otimes_{k, \Frob} k.
\]
\end{defn}

\begin{conj}[{\cite[Conjecture 6.3]{TV}}]\label{conj: TV local}
If $\pi$ is linked to $\Pi$, then $\pi$ base changes to $\Pi$.
\end{conj}

\begin{example}
The need for the Frobenius twist can be seen in a simple example. Suppose $G = H^p$ and $\sigma$ acts by cyclic permutation. Then $G^\sigma$ is the diagonal copy of $H$. In this case a representation $\pi$ of $H_v$ should transfer to $\pi^{\boxtimes p}$ of $G_v$. And indeed, 
\[
T^0(\pi^{\boxtimes p}) = \frac{\ker (1-\sigma \mid \pi^{\boxtimes p})}{N \cdot  \pi^{\boxtimes p}} \cong \pi^{(p)}.
\]
\end{example}

\begin{remark} Conjecture \ref{conj: TV local} is highly non-trivial even for groups such as $\GL_n$ where the full Local Langlands correspondence, hence in particular existence of cyclic base change, is already known. In fact, the main result of \cite{Ron16} is a special case of the conjecture, for depth-zero supercuspidal representations of $\GL_n$ compactly induced from cuspidal Deligne-Lusztig representations. Despite the very explicit nature of the Local Langlands Correspondence for such representations, the proof in \emph{loc. cit.} involves rather hefty calculations, which were not amenable to generalization. 

Our proof of Conjecture \ref{conj: TV local} (when $p$ is odd and good for $\wh{G}$) is conceptual and applies to all representations, without using any explicit models such as models for supercuspidal representations as compact inductions. Furthermore, the unramified and tamely ramified base change are handled completely differently in \cite{Ron16}, whereas our proof will be completely uniform in the field extension, the reductive group, and the irreducible representation. 
\end{remark}

\begin{thm}\label{thm: main local 2}
Assume $p$ is an odd good prime for $\wh{G}$. Let $\Pi$ be an irreducible admissible representation of $G_v$ and let 
\[
\chi_{\Pi^{(p)}} \co \Exc(W_v, \ld G) \rightarrow k
\]
be the associated character of $\Pi^{(p)}$. Form $T^*(\Pi) := T^*(\langle \sigma \rangle, \Pi)$, viewed as a smooth $H_v$-representation. Then for any irreducible character $\chi \co \Exc(W_v, \ld H) \rightarrow k$ appearing in the action on $T^*(\Pi)$ via $Z_H \co \Exc(W_v, \ld H) \rightarrow \mf{Z}(H)$, the composite character 
\[
\Exc(W_v , \ld G) \xrightarrow{\phi_{\BC}^*} \Exc(W_v, \ld H) \xrightarrow{\chi} k
\]
agrees with $\chi_{\Pi^{(p)}}$. 
\end{thm}

It is clear that Theorem \ref{thm: main local 2} implies Theorem \ref{thm: intro 2}.

\begin{proof}
Let $\Pi$ be a representation of $G_v$. Then $\mf{Z}(G)$ acts $G_v$-equivariantly on $\Pi$, inducing an $H_v$-equivariant action of $\mf{Z}(G)^{\sigma}$ on $T^*(\Pi)$. In particular, as $Z_G$ maps $\Exc(W_v, \ld G)^\sigma \subset \Exc(W_v, \ld G)$ into $\mf{Z}(G)^{\sigma}$, we get an $H_v$-equivariant action of $\Exc(W_v, \ld G)^{\sigma}$ on $T^*(\Pi)$. 

By Lemma \ref{lem: deep lattices are plain}, $K_{r}$ is plain as soon as $r \geq 1$. Taking the (filtered) colimit over $r$ in Lemma \ref{lem: brauer homomorphism} with $K = K_r$, we find that for all $S \in \Exc(W_v, \ld G)^{\sigma}$, we have
 \[
\left(\begin{array}{@{}c@{}} 
\text{the action on $T^*(\Pi)$ of} \\ Z_{G}(S) 
\end{array} \right)
 = 
\left( \begin{array}{@{}c@{}} \text{the action on $T^*(\Pi)$ of}\\ \Br(Z_{G}(S)) 
 \end{array}\right).
 \]
In other words, the diagram below commutes:
\begin{equation}\label{eq: some diagram commutes}
\begin{tikzcd}
\mf{Z}(G)^{\sigma} \ar[r]  \ar[d, "Z(\Br)"]& \End_{H_v}(T^* \Pi) \ar[d, "\sim"] \\
 \mf{Z}(H) \ar[r] & \End_{H_v}(T^*\Pi)
\end{tikzcd}
\end{equation}

On the other hand, taking the inverse limit over $r$ in Corollary \ref{cor: local operators identified} yields that
\begin{equation}\label{eq: some operators identified}
\Br(Z_{G}(S_{I, \Nm^{(p^{-1})}(V), \Nm^{(p^{-1})}(x), \Nm^{(p^{-1})}(\xi), (\gamma_i)_{i \in I}}) ) = Z_{H}(S_{I, \Res_{\BC}(V), x, \xi, (\gamma_i)_{i \in I}})
 \end{equation}
 for all $V \in \Rep_k((\ld G)^I)$. 

Combining \eqref{eq: some diagram commutes} and \eqref{eq: some operators identified} shows that
\begin{equation}\label{eq: some crucial equality}
\left(\begin{array}{@{}c@{}} 
\text{the action on $T^*(\Pi)$ of} \\ Z_{G}(S_{I, \Nm^{(p^{-1})}(V), \Nm^{(p^{-1})}(x), \Nm^{(p^{-1})}(\xi), (\gamma_i)_{i \in I}}) 
\end{array} \right)
 = 
\left( \begin{array}{@{}c@{}} \text{the action on $T^*(\Pi)$ of}\\ Z_{H}(S_{I, \Res_{\BC}(V), x, \xi, (\gamma_i)_{i \in I}})
 \end{array}\right)
 \end{equation}
for all $V \in \Rep_k((\ld G)^I)$.






From now on, assume $\Pi$ is an irreducible (smooth) representation of $G_v$. Then $\End_{G_v} (\Pi) \cong k$ (by Schur's Lemma applied to the Hecke action on the invariants of $\Pi$ for every compact open subgroup of $G_v$). The $L$-parameter attached to $\Pi$ corresponds under Proposition \ref{prop: galois rep} to the character
\[
\chi_{\Pi} \co \Exc(W_v, \ld G) \surj \Exc(W_v, \ld G) \xrightarrow{Z_G} \mf{Z}(G) \rightarrow \End_{G_v}(\Pi) \cong k. 
\]
This induces 
\[
T^0 \chi_{\Pi} \co T^0  \Exc(W_v, \ld G) \rightarrow T^0 \Exc(W_v, \ld G) \xrightarrow{T^0 Z_G} T^0 \mf{Z}(G) \rightarrow T^0 \End_{G_v}(\Pi) \cong k.
\]
The action of $T^0 \Exc(W_v, \ld G)$ on $T^*(\Pi)$ is through $T^0 \chi_{\Pi}$ composed with the natural map $\iota \co  T^0 \End_{G_v}(\Pi) \rightarrow \End_{H_v}(T^* \Pi)$.

We also consider the homomorphism 
\[
\chi_{T^* \Pi} \co \Exc(W_v, \ld H) \xrightarrow{Z_H} \mf{Z}(H) \rightarrow \End_{H_v}(T^* \Pi).
\]
Then \eqref{eq: some crucial equality} implies that for all $V \in \Rep_k((\ld G)^I)$, we have 
\begin{equation}\label{eq: 2nd crucial equality}
\iota \circ T^0 \chi_{\Pi} (S_{I, \Nm^{(p^{-1})}(V), \Nm^{(p^{-1})}(x), \Nm^{(p^{-1})}(\xi), (\gamma_i)_{i \in I}} ) = \chi_{T^* \Pi}(S_{I, \Res_{\BC}(V), x, \xi, (\gamma_i)_{i \in I}} ).
\end{equation}
Note that the fact that the right hand side of \eqref{eq: 2nd crucial equality} lies in $k$ is already non-obvious. In particular, \eqref{eq: 2nd crucial equality} implies that for any irreducible subquotient $\pi$ of $T^* \Pi$, we have
\begin{align}\label{eq: pi vs T^0Pi}
\chi_{\pi}(S_{I, \Res_{\BC}(V), x, \xi, (\gamma_i)_{i \in I}}) &=  \chi_{T^* \Pi}(S_{I, \Res_{\BC}(V), x, \xi, (\gamma_i)_{i \in I}} ) \nonumber \\
& =   (T^0 \chi_{\Pi}) (S_{I, \Nm^{(p^{-1})}(V), \Nm^{(p^{-1})}(x), \Nm^{(p^{-1})}(\xi), (\gamma_i)_{i \in I}}  )  \nonumber\\
& =  \chi_{\Pi}(S_{I, \Nm^{(p^{-1})}(V), \Nm^{(p^{-1})}(x), \Nm^{(p^{-1})}(\xi), (\gamma_i)_{i \in I}} ) .
\end{align}

The character $\chi_{\Pi^{(p)}}$ giving the $k$-linearized action of $\Exc(W_v, \ld G)$ on $\Pi^{(p)} := \Pi \otimes_{k, \Frob_p} k$ satisfies 
\begin{equation}\label{eq: frob twist action}
\chi_{\Pi^{(p)}}(S_{I, \Nm(V), \Nm(x), \Nm(\xi), (\gamma_i)_{i \in I}}) = \chi_{\Pi}(S_{I, \Nm^{(p^{-1})}(V), \Nm^{(p^{-1})}(x), \Nm^{(p^{-1})}(\xi), (\gamma_i)_{i \in I}})^p.
\end{equation}

\begin{defn}\label{defn: local tiltexc'}Similarly to Definition \ref{def: prime subalg}, let $\Exc(W_v, \ld G)' \subset \Exc(W_v, \ld G)$ be the $k$-subalgebra generated by $N \cdot \Exc(W_v, \ld G)$ and all elements of the form $\Nm (S_{I, V, x, \xi, (\gamma_i)_{i \in I}}) = S_{I, \Nm (V), \Nm(x), \Nm(\xi), (\gamma_i)_{i \in I}}$. 
\end{defn}

Then the combination of \eqref{eq: pi vs T^0Pi} and \eqref{eq: frob twist action} tells us that the action of $\Exc(W_v, \ld G)'$ on $\Pi^{(p)}$ via $\Exc(W_v, \ld G)' \rightarrow \Exc(W_v, \ld G) \xrightarrow{Z_G} \mf{Z}(G)$ is given by the character $\chi'_{\Pi^{(p)}}$ that sends
\begin{equation}\label{eq: 6 revise eq 2}
S_{I, \Nm(V), \Nm(x), \Nm(\xi), (\gamma_i)_{i \in I}} \mapsto \chi_{\pi}(S_{I, \Res_{\BC}(V), x, \xi, (\gamma_i)_{i \in I}})^p.
\end{equation}
and (using Theorem \ref{thm: global actions coincide}(ii)) $N \cdot S  \mapsto 0$ for any $S \in \Exc(\Gamma, \ld G)$. By Lemma \ref{lem: unique extension of char}, the unique extension of this character to $\Exc(W_v, \ld G)$ is (using Lemma \ref{lem: norm extension})
\begin{align*}
S_{I, V, x, \xi, (\gamma_i)_{i \in I}} & \mapsto  \chi_{\Pi^{(p)}}'(S_{I, \Nm(V), \Nm(x), \Nm(\xi), (\gamma_i)_{i \in I}})^{1/p} \\
[\eqref{eq: 6 revise eq 2} \implies] &= \chi_{\pi}(S_{I, \Res_{\BC}(V), x, \xi, (\gamma_i)_{i \in I}}) \\
&= \chi_{\pi} \circ \phi_{\BC}^* (S_{I, V, x, \xi, (\gamma_i)_{i \in I}}).
\end{align*}
On the other hand, since $\chi_{\Pi^{(p)}}$ is tautologically an extension of $\chi_{\Pi^{(p)}}'$ to $\Exc(W_v, \ld G)$, it must be the case that $\chi_{\Pi^{(p)}} = \chi_{\pi} \circ \phi_{\BC}^*$ as characters of $\Exc(W_v, \ld G)$ for any irreducible subquotient $\pi$ of $T^*(\Pi)$. 

\end{proof}

\subsection{Local mod $p$ cyclic base change}\label{ssec: local base change} In this subsection, we will prove Theorem \ref{thm: local existence}. Assumptions are as in \S \ref{ssec: TV conjectures}. We note, however, that the formulation of all the statements in \S \ref{def: BC for BC} makes sense for any local field $F_v$ of residue field $\ell \neq p$, with a suitable replacement for the Genestier-Lafforgue correspondence, and that all our arguments in this subsection apply if those statements are true for $F_v$.

 \subsubsection{Formulation of local base change}\label{ssec: local formulation} We begin by formulating a precise notion of local base change.

\begin{defn}\label{def: base change}
Let $\pi$ be an irreducible admissible representation of $H_v$ over $k$, and $\Pi$ be an irreducible admissible representation of $G_v$ over $k$. We say that \emph{$\pi$ base changes to $\Pi$} if $\rho_{\Pi} \cong \phi_{\BC} \circ \rho_{\pi} \in H^1(W_v, \wh{G}(k))$. 
\end{defn}

This definition is an approximation to the notion of base change for $L$-packets. An $L$-packet for $H_v$ should be said to base change to an $L$-packet for $G_v$ if the corresponding $L$-parameters are related by $\phi_{\BC}$. A more refined version of Definition \ref{def: base change} would declare $\pi$ to base change to $\Pi$ if the $L$-packet of $\pi$ base changes to the $L$-packet of $\Pi$, but we lack a definition of $L$-packets for general groups and representations; therefore, we use the fibers of the Genestier-Lafforgue correspondence as a substitute for $L$-packets. 

\subsubsection{Finiteness conditions on Hecke algebras} 




We will use the following recent result of Dat-Helm-Kurinczuk-Moss. We are keeping the running assumption that $p$ differs from the residue characteristic of $F_v$. 

\begin{thm}[\cite{DHKM}]\label{thm: finiteness} For every $x \in \cB(G/F_v)$ and every $r \geq 0$, and $K_r := G(F_v)_{x,r}$, the Hecke algebra $\cH(G, K_r)$ is finite over its center $\mf{Z}(G, K_r)$, which is itself a finitely generated algebra over $k$. 
\end{thm}

\begin{remark}
The paper \cite{DHKM} proves a much stronger result, where coefficients are allowed to be an arbitrary $\Z_p$-algebra. The analogous result with coefficients in a characteristic zero field is an old result of Bernstein. The case where $p$ is banal for $G$ was known to experts to follow in a similar manner from work of Vign\'{e}ras, although it is not explicitly written down in the literature. 
\end{remark}

\subsubsection{Existence of local base change} Fix $x \in \cB(H/F_v)$, and let $K_r := G(F_v)_{x,r}$ and $U_r := K_r^{\sigma}$. We prove the following theorem, which in particular implies Theorem \ref{thm: local existence}.

\begin{thm}\label{thm: mod ell base change}
Suppose $p$ is an odd good prime for $\wh{G}$. Let $\pi$ be an irreducible representation of $H_v$ over $k$, having non-zero $U_r$-fixed vectors, with $L$-parameter $\rho_{\pi} \in H^1(W_v, \wh{H}(k))$. Then there is an irreducible representation $\Pi$ of $G_v$ over $k$, having non-zero $K_r$-fixed vectors, such that $\rho_{\Pi} \cong \phi_{\BC} \circ \rho_{\pi}$. 
\end{thm}

\begin{proof}
If $r = 0$ then the result is classical, so we assume $r>0$. Then $U_r, K_r$ have prime-to-$p$ pro-order so the theory of the Bernstein center applies. Recall that the functor $\Pi \mapsto \Pi^{K_r}$ induces a bijection between irreducible admissible $G_v$-representations with non-zero $K_r$-invariants and irreducible $\cH(G_v, K_r)$-modules. It therefore suffices to construct an irreducible representation of $\cH(G_v, K_r)$ whose induced character of $\Exc(W_v, \ld G)$ is $\chi_{\pi} \circ \phi_{\BC}^*$, where $\chi_{\pi} \co \Exc(W_v, \ld H) \rightarrow k$ is the character of $\Exc(W_v, \ld H)$ corresponding to $\pi$.  

By hypothesis, we have a non-zero algebra homomorphism $\cH(H, U_r)  \rightarrow \End(\pi^{U_r})$, which has the property that the composite homomorphism 
\[
\Exc(W_v, \ld H) \xrightarrow{Z_{H,r}} \mf{Z}(H, U_r) \rightarrow \cH( H,U_r ) \rightarrow \End_k(\pi^{U_r})
\]
has kernel the maximal ideal $\mf{m}_{\pi}  = \ker ( \chi_{\pi}) \subset \Exc(W_v, \ld H)$. The Brauer homomorphism $\Br \co \cH(G, K_r)^{\sigma} \rightarrow  \cH(H, U_r)$ fits into a commutative diagram\footnote{Since we are not assuming here that $x$ is a special vertex, we cannot invoke Corollary \ref{cor: surjectivity} to say that $\Br$ is surjective, so it may not induce a map of centers.}
\[
\begin{tikzcd}
 \Exc(W_v, \ld H) \ar[r, "Z_{H,r}"] & \mf{Z}(H, U_r) \ar[r] & \cH(H, U_r)  \ar[r] & \End_k(\pi^{U_r}) \\
 \Exc(W_v, \ld G)' \ar[r, "Z_{G,r}"]   \ar[d, hook] & \mf{Z}(G, K_r)^{\sigma}\ar[r, hook]  \ar[d, hook] & \cH(G, K_r)^{\sigma}  \ar[u, "\Br"] \ar[d, hook]  \\
 \Exc(W_v, \ld G) \ar[r, "Z_{G,r}"] & \mf{Z}(G, K_r) \ar[r, hook] &  \cH( G, K_r)
\end{tikzcd}
\]
where $\Exc(W_v, \ld G)' \subset \Exc(W_v, \ld G)^\sigma$ is as in Definition \ref{defn: local tiltexc'}. Let $\mf{m}_{\pi}' \subset \Exc(W_v, \ld{G})' $ be the kernel of the map $\Exc(W_v, \ld {G})'  \rightarrow  \End_k(\pi^{U_r}) $ obtained by tracing through the diagram above. We claim that $\mf{m}_{\pi}'$ is a maximal ideal. First of all, we observe that the map $\Exc(W_v, \ld G)' \rightarrow \End_k(\pi^{U_r})$ lands in the subring of scalars $k \subset \End_k(\pi^{U_r})$, since by Corollary \ref{cor: local operators identified} the action of $\Exc(W_v, \ld {G})'$ on $\End_k(\pi^{U_r})$ factors through the action of $\Exc(W_v, \ld H)$, which is through $\chi_{\pi}$. On the other hand, since all maps in the diagram are maps of $k$-algebras, $\Exc(W_v, \ld G)' $ must surject onto the full ring of scalars $	k \subset \End_k(\pi^{U_r})$.

Note that $\Br \co \cH(G, K_r)^{\sigma}  \rightarrow \cH(H, U_r)  $ vanishes on $N \cdot \cH( G, K_r) \subset \cH(G, K_r)^{\sigma} $. By the commutativity of the bottom part of the diagram, the composition from $\Exc(W_v, \ld G)' $ to $\cH(H, U_r)$ therefore vanishes on $N \cdot \Exc(W_v, \ld G)  \subset \Exc(W_v, \ld G)' $. Therefore we may apply Lemma \ref{lem: unique extension of char} to see that the homomorphism $\chi_{\pi}' \co \Exc(W_v, \ld{G} )' \rightarrow k$ corresponding to $\mf{m}_{\pi}'$ has a \emph{unique} extension to a character $\Exc(W_v, \ld{G}) \rightarrow k$. Since Corollary \ref{cor: local operators identified} shows that $\chi_{\pi} \circ \phi_{\BC}^*$ is such an extension, its kernel must be the unique maximal ideal of $\Exc(W_v, \ld G)$ lying over $\mf{m}_{\pi}'$.

The preceding paragraph implies that the localization of $\cH(G,K_r)^{\sigma}$ at $\mf{m}_{\pi}'$ is non-zero, since the character $\chi_\pi'$ factors through this localization. Since the action of $\Exc(W_v, \ld G)'$ on $\cH(G, K_r)^\sigma$ factors through the action of $\mf{Z}(G, K_r)^\sigma$, there exists a maximal ideal $\mf{n}_{\pi}$ of $\mf{Z}(G, K_r)$ lying over $\mf{m}_{\pi}' $ at which $\cH(G, K_r)$ is supported. Since the pullback of $\mf{n}_{\pi}$ to $\Exc(W_v, \ld{G})$ contains $\mf{m}_{\pi}'$, it must equal $\ker(\chi_{\pi} \circ \phi_{\BC}^*)$ by the preceding paragraph. 

By Theorem \ref{thm: finiteness} and the Artin-Tate Lemma, $\mf{Z}(G, K_r)$ is finite over $\mf{Z}(G, K_r)^{\sigma}$ and then $\cH(G, K_r)$ is finite over $\mf{Z}(G, K_r)^{\sigma}$. So Nakayama's Lemma implies that the left $\cH(G, K_r)$-module quotient $\cH(G, K_r) /\cH(G, K_r) \mf{n}_{\pi}$ is finite-dimensional and non-zero. By design, the only maximal ideal in its support over $\Exc(W_v, \ld G)$ is $\ker(\chi_{\pi} \circ \phi_{\BC}^*)$, so there is an irreducible $\cH(G, K_r)$-subquotient $\Xi$ of $\cH(G, K_r) / \cH(G, K_r)  \mf{n}_{\pi}$ on which $\Exc(W_v, \ld G)$ acts through $ \chi_{\pi} \circ \phi_{\BC}^*$, as was to be showed. 


\end{proof}


\begin{remark}[Depth estimates]\label{rem: depth}
For applications it is useful to have control of the \emph{depth} of the base change. The proof of Theorem \ref{thm: mod ell base change} implies an estimate on the depth, which we now spell out. Recall from \cite{MP96} that the \emph{depth} of an irreducible representation $\Pi$ of $G_v$ is the minimal $r$ such that for some $x \in \cB(G/F_v)$, $\Pi^{G(F_v)_{x, r+}} \neq 0$. Let us emphasize that the definition of the Moy-Prasad filtration $\{G(F_v)_{x, r}\}$ is normalized so that $F_v^{\times}$ has value group $\Z$. 

Let $\pi$ be an irreducible representation of $H$ of depth $r$, and let $x \in \cB(H/F_v)$ such that $\pi^{H(F_v)_{x,r+}} \neq 0$. First assume that $E_v/F_v$ is \emph{unramified}. Then $\cB(H/F_v) = \cB(G/F_v)^{\Gal(E_v/F_v)}$, and we have $H(F_v)_{x,r+} = G(F_v)_{x,r+}^{\Gal(E_v/F_v)}$ \cite[\S 9]{KP}. The proof of Theorem \ref{thm: mod ell base change} shows that there exists a local base change $\Pi$ of $\pi$ such that $\Pi^{G(F_v)_{x,r+}} \neq 0$, so that $\depth(\Pi)\leq \depth(\pi)$. (The proof does not use Corollary \ref{cor: surjectivity} or the Treumann-Venkatesh homomorphism.) 

Next suppose $E_v/F_v$ is \emph{tamely ramified}. By \cite{Pra20} we still have $\cB(H/F_v) = \cB(G/F_v)^{\Gal(E_v/F_v)}$, and by \cite[Proposition 12.9.2]{KP} we have $(G(F_v)_{x, r^+})^{\Gal(E_v/F_v)} = H(F_v)_{x, r^+}$ for all $r \geq 0$. Hence in this case, the proof of Theorem \ref{thm: mod ell base change} shows that there exists a local base change $\Pi$ of $\pi$ such that $\Pi^{G(F_v)_{x,r+}} \neq 0$, so that $\depth(\Pi)\leq  \depth(\pi)$. Let us caution, however, that if we regard $x \in \cB(H/E_v)$ instead of $\cB(G/F_v)$ and $\Pi$ as a representation of $H(E_v)$ instead of $G(F_v)$, then it is natural to define the Moy-Prasad filtration $H(E_v)_{x,r}$ so that $E_v^{\times}$ has value group $\Z$, for which $H(E_v)_{x,pr}= G(F_v)_{x, r}$. Hence, in this normalization our estimate would instead be ``$\depth(\Pi) \leq p \cdot \depth(\pi)$''. 

In both cases, the inequalities we obtain are expected to be optimal \cite{AL10}. 
\end{remark}

\appendix

\section{The base change functor realizes Langlands functoriality\\ by Tony Feng and Gus Lonergan}\label{sec: categorical base change}

In this section we prove Theorem \ref{thm: BC on Gr}. 

First we recall some general properties of Smith theory for \emph{schemes}.  

\subsection{Recollections on Smith theory for schemes}\label{ssec: six-functor}

The Tate category for schemes enjoys a robust 6-functor formalism (observed in the topological case in \cite[\S 4.3]{Tr19}, and proved for schemes in \cite[\S 2,3]{RW}). Let us recall the statements for later use. Let $f \co Y \rightarrow S$ be a $\sigma$-equivariant morphism of varieties with admissible $\sigma$-action, over a field of characteristic $\ell \neq p$. Let $\Lambda$ be a $p$-adic ring of coefficients; we are most interested in $\Lambda \in \{ W(k), k\}$. 

\begin{itemize}
\item The pullback functor $f^* \co D^b_{c}(S^{\sigma};\Lambda[\sigma]) \rightarrow D^b_{c}(Y^{\sigma}; \Lambda[\sigma])$ descends to
\[
f^* \co \Perf(S^{\sigma}; \cT_{\Lambda}) \rightarrow \Perf(Y^{\sigma}; \cT_{\Lambda}).
\]
The proper pushforward $Rf_!  \co   D^b_{c}(Y^{\sigma}; \Lambda[\sigma]) \rightarrow D^b_{c}(S^{\sigma};\Lambda[\sigma])$ descends to
\[
Rf_! \co \Perf(Y^{\sigma}; \cT_{\Lambda}) \rightarrow \Perf(S^{\sigma}; \cT_{\Lambda}).
\]
\item As Verdier duality $\DD \co D^b_{c,\sigma}(Y^{\sigma}; \Lambda) \rightarrow D^b_{c,\sigma}(Y^{\sigma}; \Lambda)$ preserves $\Perf(Y^{\sigma}; \Lambda[\sigma])$, it descends to the Tate category to define 
\[
\DD \co \Perf(Y^{\sigma}; \cT_{\Lambda}) \rightarrow \Perf(Y^{\sigma}; \cT_{\Lambda}).
\]
Using this, we may define the operations 
\[
f^! := \DD \circ f^* \circ \DD \co \Perf(S^{\sigma}; \cT_{\Lambda}) \rightarrow \Perf(Y^{\sigma}; \cT_{\Lambda})
\]
and 
\[
Rf_* := \DD \circ  f_! \circ \DD \co \Perf(Y^{\sigma}; \cT_{\Lambda}) \rightarrow \Perf(S^{\sigma}; \cT_{\Lambda}).
\]
\end{itemize}



We now list some properties which could be remembered under the slogan, ``the Smith operation commutes with all operations'' (cf. \cite[\S 4.4]{Tr19}). 

\subsubsection{Compatibility with pullback}\label{sssec: compatibility with pullback} If $f$ satisfies the assumptions above, then the following diagrams commute: 
\[
\begin{tikzcd}
D^b_{c,\sigma}(Y;\Lambda) \ar[d, "\Psm"] &  D^b_{c,\sigma}(S;\Lambda)   \ar[d, "\Psm"]  \ar[l, "f^*"] \\
\Perf(Y^{\sigma}; \Cal{T}_{\Lambda}) & \Perf(S^{\sigma} ; \Cal{T}_{\Lambda}) \ar[l, "f^*"] 
\end{tikzcd} \hspace{2cm} 
\begin{tikzcd} 
D^b_{c,\sigma}(Y;\Lambda) \ar[d, "\Psm"] &  D^b_{c,\sigma}(S;\Lambda)   \ar[d, "\Psm"]  \ar[l, "f^!"] \\
\Perf(Y^{\sigma}; \Cal{T}_{\Lambda}) & \Perf(S^{\sigma}; \Cal{T}_{\Lambda}) \ar[l, "f^!"] 
\end{tikzcd}
\]
The proof for the first square is formal; for the second it follows immediately from the first plus \cite[Lemma 3.5]{RW}, whose proof is the same as that for Lemma \ref{lem: * vs !}.

\subsubsection{Compatibility with pushforward}\label{sssec: compatibility with pushforward} If $f$ satisfies the assumptions above, then the following diagrams commute: 
\[
\begin{tikzcd}
D^b_{c,\sigma}(Y;\Lambda) \ar[d, "\Psm"] \ar[r, "Rf_*"]  &  D^b_{c,\sigma}(S;\Lambda)   \ar[d, "\Psm"]  \\
\Perf(Y^{\sigma}; \Cal{T}_\Lambda)  \ar[r, "Rf_*"]  & \Perf(S^{\sigma}; \Cal{T}_\Lambda)
\end{tikzcd} \hspace{2cm}
\begin{tikzcd}
D^b_{c,\sigma}(Y;\Lambda) \ar[d, "\Psm"] \ar[r, "Rf_!"]  &  D^b_{c,\sigma}(S;\Lambda)   \ar[d, "\Psm"]  \\
\Perf(Y^{\sigma}; \Cal{T}_\Lambda)  \ar[r, "Rf_!"]  & \Perf(S^{\sigma}; \Cal{T}_\Lambda)
\end{tikzcd} 
\]
The proof for the second diagram is the same as that of Proposition \ref{prop: equivariant localization}. Then the commutativity of the first diagram follows by applying Verdier duality and using Lemma \ref{lem: * vs !}.

\subsection{Setup for the proof of Theorem \ref{thm: BC on Gr}}

We keep the setup of \S \ref{sssec: base change setup}: $H$ is any reductive group over a separably closed field $\F$ of characteristic $\neq p$, and $G = H^p$. We let $\sigma$ act on $G$ by cyclic rotation, sending the $i$th factor to the $(i+1)$st (mod $p$) factor. 



\subsection{Proof of additivity}\label{ssec: linearity}

We first prove that $\BC$ is additive, i.e., we exhibit a natural isomorphism $\BC(\Cal{F} \oplus \Cal{F}') \cong \BC(\Cal{F}) \oplus \BC(\Cal{F}')$. We have
\begin{align*}
\Nm(\Cal{F} \oplus \Cal{F}')  &=  (\Cal{F} \oplus \Cal{F}') * ({}^{\sigma}\Cal{F} \oplus {}^{\sigma}\Cal{F}') * \ldots * ({}^{\sigma^{p-1}}\Cal{F} \oplus {}^{\sigma^{p-1}}\Cal{F}') \\
&\cong \Nm(\Cal{F}) \oplus \Nm(\Cal{F}') \oplus (\text{direct sum of free $\sigma$-orbits}). 
\end{align*}
Therefore, the restrictions of $\Nm(\Cal{F} \oplus \Cal{F}') $ and $ \Nm(\Cal{F}) \oplus \Nm(\Cal{F}')$ to $Y^{\sigma}$ differ by a perfect complex of $\OO[\sigma]$-modules, and hence project to isomorphic objects in the Tate category $\Perf(Y^{\sigma}; \Cal{T}_{\OO})$. This shows that $\Psm \circ \Nm$ is additive. Since the lifting functor $L$ is also additive, diagram \eqref{eq: dashed arrow} shows that $\BC^{(p)} \circ \FF$ is additive. Since $\FF$ is essentially surjective as our assumptions imply that all parity sheaves exist, $\BC^{(p)}$ is additive. Finally, $\Frob_p^{-1}$ is an equivalence so also additive, so $\BC$ is additive. \qed

\subsection{Reduction to the case of a torus}\label{ssec: reduction to torus}

Fix a regular $\sigma$-equivariant cocharacter $\kappa \co \G_m \rightarrow H$. Its centralizer in $G$ is a maximal torus $T \subset G$, such that $T_H := T \cap H$ is a maximal torus of $H$. We recall the following statements relating the restriction functor $\Rep(\wh{H}) \rightarrow \Rep(\wh{T}_H)$ with the \emph{hyperbolic localization} functor under the Geometric Satake equivalence. 

The fixed points of $\G_m$ acting by left translation on $\Gr_G$ via $\kappa$ are the $t^{\nu}$ for $\nu \in X_*(T)$. The attracting locus to $t^{\nu}$ is the semi-infinite orbit $S_{\nu}$. These semi-infinite orbits form a stratification of $\Gr_G$. Let $i_{\nu} \co S_{\nu} \rightarrow \Gr_G$ be the inclusion (a locally closed embedding of ind-schemes) and $p_{\nu} \co S_{\nu} \rightarrow t^\nu$, viewed as a point of $\Gr_T$. The hyperbolic localization functor (for $G$) is the functor
\[
\CT_G = \bigoplus_{\nu \in X_*(T)} Rp_{\nu !} i_{\nu}^* \co D^b_c(\Gr_G)^{\mon} \rightarrow D^b_c(\Gr_T)^{\mon}
\]
where the superscript $\mon$ means monodromic for the $\G_m$-action via $\kappa$ (i.e., the full subcategory spanned by objects pulled back from the $\G_m$-equivariant derived category). 

Denote by $[\deg_G]$ the function $X_*(T) \xrightarrow{\langle 2 \rho_G, - \rangle} \Z$, and similarly for $H$. Set  
\[
\CT_G[\deg_G] :=\bigoplus_{\nu \in X_*(T)} Rp_{\nu !} i_{\nu}^* [\deg_G(\nu)]
\]
and $\CT_H[\deg_H]$ similarly. Then $\CT_G[\deg_G]$ and $\CT_H[\deg_H]$ are t-exact and under the Geometric Satake equivalence, they are intertwined with restriction along $\wh{T} \rightarrow \wh{G}$ and $\wh{T}_H \rightarrow \wh{H}$, respectively:
\begin{equation}\label{eq: CT restrict}
\begin{tikzcd}[column sep = huge, row sep = huge]
\Perv_{L^+G}(\Gr_G; \Lambda) \ar[r, "\sim", "\text{Geom. Sat.}"'] \ar[d, "{\CT_G[\deg_G]}"] & \Rep_{\Lambda}(\wh{G}) \ar[d, "\Res^{\wh{G}}_{\wh{T}}"] \\
\Perv_{L^+T} (\Gr_T; \Lambda) \ar[r, "\sim", "\text{Geom. Sat.}"'] & \Rep_{\Lambda}(\wh{T})
\end{tikzcd}
\hspace{1cm}
\begin{tikzcd}[column sep = huge, row sep = huge]
\Perv_{L^+G}(\Gr_H; \Lambda) \ar[r, "\sim", "\text{Geom. Sat.}"'] \ar[d, "{\CT_H[\deg_H]}"] & \Rep_{\Lambda}(\wh{H}) \ar[d, "\Res^{\wh{H}}_{\wh{T}_H}"] \\
\Perv_{L^+T} (\Gr_{T_H}; \Lambda) \ar[r, "\sim", "\text{Geom. Sat.}"'] & \Rep_{\Lambda}(\wh{T}_H)
\end{tikzcd}
\end{equation}
(Here $\CT$ is defined because equivariance implies monodromicity.) 

\begin{remark}
In the stated generality -- with the scheme-theoretic $\Gr_G$ in equal characteristic, and the coefficients being modular \'{e}tale sheaves -- the commutative diagram \eqref{eq: CT restrict} is perhaps not completely documented in the literature. It does appear for general coefficients on the complex affine Grassmannian \cite[p.66]{BR18} and the $B_{\mrm{dR}}^+$-affine Grassmannian (in arbitrary characteristic) \cite[p.233]{FS}. The proofs in either case are essentially the same -- the commutativity of the diagram is baked into the step of identifying the Tannakian group $\wh{G}$ -- and they carry over essentially verbatim to our setting. 
\end{remark}

The functor $\CT_H[\deg_H]$ induces 
\[
D^b_{c, \sigma}(\Gr_H)^{\mon} \rightarrow D^b_{c, \sigma}(\Gr_{T_H})^{\mon}
\]
Since $*/!$-restriction and $*/!$-pushforward all commute with $\Psm$ by \S \ref{ssec: six-functor}, the Constant Term functor commutes with $\Psm$ in the sense of the following commutative diagram 
\begin{equation}\label{eq: CT PSm}
\begin{tikzcd}
D^b_{c, L^+G \rtimes \sigma}(\Gr_G; \OO) \ar[r, "\Psm"] \ar[d, "{\CT_G[\deg_G]}"] & \Perf_{(L^+H)}(\Gr_H; \cT_{\OO}) \ar[d, "{\CT_H[\deg_H]}"] \\
D^b_{c, L^+G \rtimes \sigma}(\Gr_T; \OO) \ar[r, "\Psm"] & \Perf_{(L^+H)}(\Gr_{T_H}; \cT_{\OO})
\end{tikzcd}
\end{equation}
and the same holds with the shifts by $\deg_H$ and $\deg_G$, thanks to the parity calculations in \S \ref{sssec: base change functor}.

\begin{lemma}\label{lem: cube}
Consider the cube 
\begin{equation}\label{eq: cube}
\begin{tikzcd}
& & \Tilt_k(\wh{G}) \ar[rrr, "\Res_{\BC}"] \ar[ddd, "\Res^{\wh{G}}_{\wh{T}}"] & & & \Tilt_k(\wh{H}) \ar[ddd, "\Res^{\wh{H}}_{\wh{T}_H}"] \\
\\ 
\Parity^0_{L^+G}(\Gr_G;k) \ar[rrr, "\BC"] \ar[ddd, "{\CT_G[\deg_G]}"] \ar[uurr]  & & &  \Parity^0_{L^+H}(\Gr_H;k)  \ar[ddd, "{\CT_H[\deg_H]}"] \ar[uurr] \\
& & \Tilt_k(\wh{T}) \ar[rrr, "\Res_{\BC}"] &  & & \Tilt_k(\wh{T}_H) \\ 
\\
\Parity^0_{L^+T}(\Gr_T;k) \ar[rrr, "\BC"] \ar[uurr]  & & & \Parity^0_{L^+ T_H}(\Gr_{T_H};k) \ar[uurr] 
\end{tikzcd}
\end{equation}
where all diagonal arrows are the Geometric Satake equivalence (using Theorem \ref{thm: parity = tilting}). The back, front, left, and right faces commute. 
\end{lemma}

\begin{proof}
The back face obviously commutes. The left and right faces commute by \eqref{eq: CT restrict}. It remains to analyze the front square. 

Consider the diagram 
\begin{equation}\label{eq: cube diag 1}
\begin{tikzcd}
\Parity^0_{L^+G}(\Gr_G; \OO) \ar[r, "\Nm"] \ar[d, "{\CT_G[\deg_G]}"] & \Parity^0_{L^+G \rtimes \sigma}(\Gr_G; \OO) \ar[d, "{\CT_G[\deg_G]}"] \ar[r, "\Psm"] & \Parity_{(L^+H)}(\Gr_H; \cT_{\OO}) \ar[r, "L"] \ar[d, "{\CT_H[\deg_H]}"] & \Parity^0_{L^+H}(\Gr_H; k) \ar[d, "{\CT_H[\deg_H]}"]\\
\Parity^0_{L^+T}(\Gr_T; \OO) \ar[r, "\Nm"]  & \Parity^0_{L^+T \rtimes \sigma}(\Gr_T; \OO) \ar[r, "\Psm"]  & \Parity_{(L^+ T_H)}(\Gr_{T_H}; \cT_{\OO}) \ar[r, "L"] & \Parity^0_{L^+ T_H}(\Gr_{T_H}; k)
\end{tikzcd}
\end{equation}
The left square commutes because $\CT_G[\deg G]$ is symmetric monoidal. The middle square commutes because $\Psm$ is compatible with $*$-pullback and $!$-pushforward, as explained in \S \ref{ssec: six-functor}. We claim that the right square commutes. To see this, we consider the diagram 
\begin{equation}\label{eq: cube diag 2}
\begin{tikzcd}
\Parity^0_{L^+H \times \sigma}(\Gr_H; \OO) \ar[r, "\TT^* \epsilon^*"] \ar[d, "{\CT_H[\deg_H]}"]\ar[rr, bend left, "\FF"] & \Parity_{(L^+H)}(\Gr_H; \cT_{\OO}) \ar[r, "L"]  \ar[d, "{\CT_H[\deg_H]}"] & \Parity^0_{L^+ H}(\Gr_H; k)  \ar[d, "{\CT_H[\deg_H]}"] \\
\Parity^0_{L^+H \times \sigma}(\Gr_{T_H}; \OO) \ar[r, "\TT^* \epsilon^*"] \ar[rr, bend right, "\FF"]  & \Parity_{(L^+ T_H)}(\Gr_{T_H}; \cT_{\OO}) \ar[r, "L"] & \Parity^0_{L^+ T_H}(\Gr_{T_H}; k)
\end{tikzcd}
\end{equation}
The upper and lower caps commute by \eqref{eq: lifting triangle}. Then it is immediate from the definition of the modular reduction functor $\FF$ that the outer square commutes. In the left square, the vertical arrows are essentially surjective since all (Tate-)parity sheaves exist for all strata. The maps on morphisms are given by \eqref{eq: morphisms between Tate-parity}. Hence the outer commutative diagram induces the right one. 

We have now established that the outer rectangle in \eqref{eq: cube diag 1} commutes. Therefore, by \eqref{eq: dashed arrow}, the diagram 
\[
\begin{tikzcd}
\Parity^0_{L^+G}(\Gr_G; k) \ar[r, "\BC^{(p)}"] \ar[d, "{\CT_G[\deg_G]}"] & \Parity^0_{L^+H} (\Gr_H; k) \ar[d, "{\CT_H[\deg_H]}"] \\
\Parity^0_{L^+T}(\Gr_T;k) \ar[r, "\BC^{(p)}"] & \Parity^0_{L^+T_H}(\Gr_{T_H};k)
\end{tikzcd}
\]
commutes. Finally, applying the Frobenius linearization process of Definition \ref{defn: BC on Gr} completes the proof for commutativity of the front face of \eqref{eq: cube}. 

\end{proof}

Theorem \ref{thm: BC on Gr} is the statement that the top face commutes. The bottom face is the special case of Theorem \ref{thm: BC on Gr} for a torus, which we will check directly. We may reduce the general case to the torus case as follows. The restriction functor $\Rep(\wh{H}) \rightarrow \Rep(T_{\wh{H}})$ is faithful and injective on tilting objects (i.e. ``tilting modules are determined by their characters'') by \cite[p. 46]{Don93}. Hence, by Lemma \ref{lem: cube}, to check that the top face commutes it suffices to check that the bottom face commutes, i.e., to prove Theorem \ref{thm: BC on Gr} in the special case where $H$ is a \emph{torus}.

\subsection{Proof in the case of a torus}

Finally, we examine the case when $H$ is a torus. Since the theorem is compatible with products, we can even reduce to the case $H = \G_m$. For $H = \G_m$ the underlying reduced scheme of $\Gr_H$ is a disjoint union of points labeled by the integers. 

The irreducible algebraic representations of $\wh{H}$ are indexed by $n \in \Z$, with $V_n \in \Rep(\wh{H})$ corresponding to the constant sheaf supported on the component $\Gr_{H}^n$ labeled by $n$. The irreducible algebraic representations of $\wh{G}$ are then labeled by $p$-tuples of integers $(n_1, \ldots, n_p) \in \Z^p$. By the additivity of $\BC$ established in \S \ref{ssec: linearity} and the complete reducibility of algebraic representations of tori, we may assume that $\Cal{F}$ is irreducible, say $\Cal{F}  = \Cal{F}(n_1, \ldots, n_p)$ is the constant sheaf supported on $\Gr_{G}^{(n_1, \ldots, n_p)}$. Then the $\sigma$-equivariant sheaf $\Nm(\Cal{F})$ is the constant sheaf $\ul{k}$ supported on the component $\Gr_{G}^{ (n_1 + \ldots + n_p, \ldots, n_1 + \ldots + n_p)}$. Its restriction to the diagonal copy of $\Gr_H$ is the constant sheaf with value $k$ supported on $\Gr_H^{n_1+\ldots+n_p}$. This is already an indecomposable $k$-parity sheaf, which tautologically lifts its own image in the Tate category. Hence we have shown that 
\[
\ul{k}_{\Gr_H^{n_1+\ldots+n_p}} = \BC^{(p)}(V_{n_1, \ldots, n_p}).
\]
And indeed, this is precisely the sheaf which corresponds under geometric Satake to $\Res_{\BC}(V_{n_1} \boxtimes V_{n_2} \boxtimes \ldots \boxtimes V_{n_p}) \cong V_{n_1 + n_2 + \ldots + n_p} \in \Rep(\wh{H})$. This confirms the commutativity of the diagram 
\[
\begin{tikzcd}
 \Parity^0_{L^+G}(\Gr_G; k) \ar[d, "\sim"]  \ar[r, "\BC"] &  \Parity^0_{L^+H}(\Gr_H; k) \ar[d, "\sim"] \\
\mrm{Tilt}_k(\wh{G}) \ar[r, "\Res_{\BC}"] & \mrm{Tilt}_k(\wh{H})
\end{tikzcd}
\]
at the level of objects. Our final step is to verify the commutativity on morphisms. Since (as $H$ is a torus) the categories involved are all semi-simple, the commutativity at the level of morphisms reduces to examining a scalar endomorphism of the simple object $\Cal{F}$ above, which corresponds to the simple representation $V_{n_1, \ldots, n_p}$. The restriction functor $\Res_{\BC}$ is $k$-linear, so what we have to check is that $\BC$ sends multiplication by $\lambda$ on $\Cal{F}$ to multiplication by $\lambda$ on $\BC(\Cal{F})$. Now, multiplication by $\lambda$ on $\Cal{F}$ is sent under $\Nm$ to multiplication by $\lambda^p$ on $\Nm(\Cal{F})$, which restricts to multiplication by $\lambda^p$ on $\BC^{(p)}(\Cal{F})$. Then the inverse Frobenius twist $\Frob_p^{-1}$ sends it to multiplication by $\lambda$, so $\BC := \Frob_p^{-1} \circ \BC^{(p)}$ behaves as desired. 
\qed

\subsection{Proof of Lemma \ref{lem: exactness}}\label{ssec: proof of exactness}
In this subsection we prove Lemma \ref{lem: exactness}. We keep the notations from \S \ref{ssec: reduction to torus}. 

\begin{lemma}\label{lem: conservativity}
The functor $\CT \co \Perf_{(L^+H)}(\Gr_H; \cT_k) \rightarrow \Perf_{(L^+T_H)}(\Gr_{T_H}; \cT_k)$ is conservative. 
\end{lemma}

\begin{proof}
Suppose $\CT(\cK) = 0$ for some non-zero $\cK \in \Shv_{(L^+H)}(\Gr_H; \cT_k)$. Let $\Gr_H^\lambda = L^+H t^\lambda \subset \Gr_H$ be the maximal stratum (for the closure order) on which $\cK$ is supported. Then for $w_0$ the longest Weyl element, the semi-infinite orbit $S_{w_0(\lambda)}$ intersects $\Gr_H^\lambda$ in a single point $t^{w_0(\lambda)}$ \cite[Theorem 5.2]{BR18}. Hence the stalk at $t^{w_0(\lambda)}$ vanishes, which shows (by the assumed constructibility for $L^+H$-orbits) that $\cK$ vanishes on $\Gr_\lambda$, which contradicts the assumption on the support of $\cK$. 
\end{proof}

Let $i \co \Gr_H \inj \Gr_G$. Since apply $\Frob_p^{-1}$ preserves exact sequences, Lemma \ref{lem: exactness} is equivalent to: if $A \rightarrow B \rightarrow C$ is an exact sequence in $\Rep_k(\wh{G})$, then the map in $D^b_{c,L^+H}(\Gr_H; k[\sigma])$,
\[
\mrm{Cone}\left[ i^*  \Nm  \Sat (A) \rightarrow i^* \Nm \Sat (B) \right] \rightarrow i^* \Nm \Sat(C)
\]
projects to an isomorphism in $\Perf_{(L^+H)}(\Gr_H; \cT_k)$. By Lemma \ref{lem: conservativity}, this can be checked after applying $\CT[\deg_G]$. Using the commutative diagram analogous to \eqref{eq: CT PSm} but with $k$-coefficients, and that $\CT[\deg_G]$ is intertwined with restriction from $\wh{G}$ to $\wh{T}$ under Geometric Satake, we have a commutative diagram 
\[
\begin{tikzcd}
\Rep_k(\wh{G})  \ar[r, "\Sat"] \ar[d, "\Res"] &  \Perv_{L^+G}(\Gr_G; k) \ar[r, "\Nm"] \ar[d, "{\CT[\deg_G]}"] &  \Perv_{L^+G \rtimes \sigma}(\Gr_G; k) \ar[d, "{\CT[\deg_G]}"]  \ar[r, hook] & D^b_{L^+H ,c}(\Gr_H; k[\sigma]) \ar[r, "\TT^*"] \ar[d, "{\CT[\deg_G]}"]  &  \Perf_{(L^+H)}(\Gr_H; \cT_k) \ar[d, "{\CT[\deg_G]}"]   \\
\Rep_k(\wh{T})  \ar[r, "\Sat"] &  \Perv_{L^+T}(\Gr_T; k) \ar[r, "\Nm"] &  \Perv_{L^+T \rtimes \sigma}(\Gr_T; k) \ar[r, hook] &  D^b_{L^+T_H ,c}(\Gr_{T_H}; k[\sigma]) \ar[r, "\TT^*"] & \Perf_{(L^+T_H)}(\Gr_{T_H}; \cT_k)
\end{tikzcd}
\] 
The question of whether the composition of functors along the top then right is an isomorphism is equivalent to the question of whether the composition of functors along the left then bottom is an isomorphism. This reduces us to the case where $\wh{G} = \wh{T}$ is a torus. In this case, since all maps in $\Rep(\wh{T})$ have splittings, the exactness statement reduces to the additivity, which was verified in \S \ref{ssec: linearity}. 

\section{Applying Drinfeld's Lemma to Tate cohomology}\label{app: B}

Here we prove Proposition \ref{prop: FWeil on shtukas}, that the action of $\FWeil(\eta^I, \ol{\eta}_I)$ on $T^j(\Sht_{G,D,I}; V)$ factors through the quotient $\FWeil(\eta^I, \ol{\eta^I}) \surj \Weil(\eta, \ol{\eta})^I$. This statement is analogous to results in \cite{Laff18} and \cite{Xue21} for ordinary cohomology, the latter of which incorporates simplifications by \cite{XZ}, and our argument will follow the same broad lines. Here is an outline of the strategy: 
\begin{enumerate}
\item Prove an Eichler-Shimura relation, relating the action of partial Frobenii and Hecke operators. 
\item Using (1), express $T^j(\Sht_{G,D,I};  V)$ as the filtered colimit of submodules stable under the partial Frobenii, with each module being finite type over some finitely generated $k$-algebra (depending on the submodule; it will be taken to a suitable tensor product of local Hecke algebras). 
\item Apply Drinfeld's Lemma, which says roughly that any continuous $A$-linear $\FWeil(\eta^I, \ol{\eta}_I)$-action on a finite-type $A$-module automatically factors over $\Weil(\eta, \ol{\eta})^I$, to each of the submodules produced in (1). 
\end{enumerate}

\subsection{The Eichler-Shimura relation}
Regarding the first step, Xue proves: 

\begin{prop}[{\cite[Proposition 7.2.6]{Xue21}}]\label{prop: xue eichler-shimura}
Let $v \in \circX$ be a closed point, with degree $\deg v$. For any finite set $I = \wt{I} \sqcup \{0\}$ and any $V \in \Rep_k(\wh{G}^I)$, there exists $W \in \Rep_k(\wh{G})$ such that 
\[
\sum_{\alpha = 0}^{\dim W} (-1)^\alpha S_{\wedge^{\dim W - \alpha} W, v}(F_{\{0\}}^{\deg v})^\alpha = 0 \in \End_{D^c((\circX \setminus D)^{\wt{I}} \times v; k)}(R^j \pi_{I!} (\circSht_{G,D,I}|_{(\circX \setminus D)^{\wt{I}} \times v}; \Sat(V))).
\]
\end{prop}
Here for a representation $W$ of $\wh{G}$, the operator $S_{W,v}$ is defined in \cite[\S 6]{Laff18} by a process similar to one defining excursion operators. The only thing we need to know about $S_{W,v}$ is the following.

\begin{thm}[``S=T Theorem'']\label{B: s=t}
 $S_{W,v}$ agrees with the Hecke operator $T_{W,v}$ after restricting to $(X \setminus (D \cup v))^I$.
 \end{thm}
 
  This fact is proved in \cite[\S 6]{Laff18} in characteristic zero; a simpler proof is a consequence \cite[Theorem 6.0.1(2)]{XZ}, which already shows the equality at the level of cohomological correspondences on local shtukas. The argument of Xiao-Zhu is written with integral coefficients in \cite[Theorem 5.1 and Corollary 5.5]{Yu22}. Since it holds at the level of cohomological correspondences, it holds in particular for Tate cohomology. 

Xue's proof of Proposition \ref{prop: xue eichler-shimura} (which is a small generalization of an argument appearing in \cite[\S 6]{XZ}) 
works essentially verbatim for Tate cohomology, replacing her $\cH^{j, \cO_E}_{G, N, I, W}$ by $T^j \pi_{I!} (\Sht_{G,D,I}|_{(\circX \setminus D)^{\wt{I}} \times v}; \Sat(V)))$. It yields:

\begin{lemma}\label{B: eichler-shimura}
Let $v \in X$ be a closed point, with degree $\deg v$. For any finite set $I = \wt{I} \sqcup \{0\}$ and any $V \in \Rep_{k}(\wh{G}^I)$, there exists $W \in \Rep_k(\wh{G})$ such that 
\[
\sum_{\alpha = 0}^{\dim W} (-1)^\alpha S_{\wedge^{\dim W - \alpha} W, v}(F_{\{0\}}^{\deg v})^\alpha = 0  \in \End_{D((\circX \setminus D)^{\wt{I}} \times v; k)}(T^j R\pi_{I!} (\Sht_{G,D,I}|_{(\circX \setminus D)^{\wt{I}} \times v}; \Sat(V))).
\]
\end{lemma}

\subsection{The filtration}\label{sssec: filtration} We carry out Step (2) of the outline, following \cite[\S 1]{Xue21}. 

Harder-Narasimhan truncation presents $\Sht_{G,D,I}$ as a filtered colimit
\[
\Sht_{G,D,I} = \colim_{\mu} \Sht_{G,D,I}^{\leq \mu}
\]
where $\mu$ runs over dominant coweights of $G$. 

Since the support of $\Sat(V)$ on $\Sht_{G,D,I}^{\leq \mu}|_{\ol{\eta^I}}$ is of finite type over $\ol{\eta}^I$, $T^j(\Sht_{G,D,I}^{\leq \mu};V) := T^j(R\Gamma_c(\Sht_{G,D,I}^{\leq \mu}|_{\ol{\eta^I}}; \Sat(V)))$ is finite-dimensional over $k$. We have a filtered colimit 
\[
T^j(\Sht_{G,D,I};V) \cong \colim_{\mu} T^j(\Sht_{G,D,I}^{\leq \mu}; V).
\]
We will express $T^j(\Sht_{G,D,I}; V)$ as an increasing union of submodules $\mf{M}^{\mu}$ which are stable under $\FWeil(\eta^I, \ol{\eta^I})$, and such that for each $\mf{M}^{\mu}$ there is a finite set of points $v_i$ (depending on $\mf{M}^{\mu}$) so that $\mf{M}^{\mu}$ is stable under the partial Frobenii and finite type over $\otimes_{i \in I} \cH_{G, v_i}$. 

Write 
\[
\mf{T} := R\pi_{I!} (\Sht_{G,D,I}|_{(\circX \setminus D)^I}; \Sat(V)) \in D^b((\circX \setminus D)^I; k)
\]
and 
\[
\mf{T}^{\mu} :=  R\pi_{I!} (\Sht_{G,D,I}^{\leq \mu}|_{(\circX \setminus D)^I}; \Sat(V)) \in D^b_c((\circX \setminus D)^I; k).
\]
Note that proper base change gives an isomorphism $T^j(\Sht_{G,D,I}; V) \cong \mf{T}|_{\ol{\eta^I}}$. Since $ \mf{T}^{\mu}$ is constructible, there is an open dense subscheme $\Omega \subset (\circX \setminus D)^I$ such that $ \mf{T}^\mu$ is a local system over $\Omega$. Choose a closed point $v \in \Omega$ and let $v_i = \pr_i(v)$ for $i \in I$. Then $\times_{i \in I} v_i \in (\circX \setminus D)^I$ is a finite union of closed points containing $v$. Let $\mf{M}^{\mu}$ be the subspace of $T^j(\Sht_{G,D,I}; V)$ given by 
\begin{equation}\label{B: M^mu}
\mf{M}^\mu := \sum_{n_i \in \N^I}  \left( \otimes_{i \in I} \cH_{i,v_i} \right) \prod_{i \in I} F_{\{i\}}^{n_i} \left( \prod_{i \in I} (\Frob_{\{i\}}^{n_i})^* \Ima\left(\mf{T}^{\mu}|_{\ol{\eta^I}} \rightarrow \mf{T}|_{\ol{\eta^I}} \right) \right) .
\end{equation}
Then it is clear that $ T^j(\Sht_{G,D,I};V) = \bigcup_{\mu} \mf{M}^{\mu}$. 

We regard $\mf{M}^{\mu}$ as a module over the finite type $k$-algebra $A^{\mu} := \otimes_{i \in I} \cH_{i,v_i}$. The following Lemma and its proof are variants of \cite[Lemma 1.3.11]{Xue21}. 

\begin{lemma}\label{B: filtration}
The submodule $\mf{M}^{\mu} \subset \mf{T}|_{\ol{\eta^I}}  \cong T^j(\Sht_{G,D,I}; V)$ is stable under the partial Frobenii $F_{\{i\}}$ and of finite type over $A^\mu$. 
\end{lemma}

\begin{proof}The stability under partial Frobenii is clear by construction. Let $\ol{v}$ be a geometric point over $v$. We have a specialization map $\fsp \co \ol{\eta^I} \rightsquigarrow \ol{v}$. For any $n_i$, we have the partial Frobenius 
\[
F_{\{i\}}^{n_i \deg(v_i) } \co (\Frob_{\{i\}}^{n_i \deg(v_i) })^* \mf{T}^{\mu} \rightarrow \mf{T}.
\]
Altough partial Frobenius does not preserve the HN truncation, there exists $\kappa$ fitting into a commutative diagram 
\begin{equation}\label{B: specialization}
\begin{tikzcd}
(\Frob_{\{i\}}^{n_i \deg(v_i) })^*  \mf{T}^\mu|_{\ol{v}} \ar[d, "F_{\{i\}}^{n_i \deg(v_i) }"] \ar[r, "\mf{sp}^*"] & 
(\Frob_{\{i\}}^{n_i \deg(v_i) })^*  \mf{T}^{\mu}|_{\ol{\eta^I}} \ar[d, "F_{\{i\}}^{n_i \deg(v_i) }"] \\
\mf{T}^{\mu+\kappa}|_{\ol{v}} \ar[r, "\mf{sp}^*"] & \mf{T}^{\mu+\kappa}|_{\ol{\eta^I}}
\end{tikzcd}
\end{equation}
We have $\Frob_{\{i\}}^{n_i \deg(v_i) }(v) = v \in \Omega$. Then using Proposition \ref{B: eichler-shimura}, we may eliminate all powers of partial Frobenius with exponent $\geq \dim W$ in \eqref{B: M^mu} in terms of $S$-operators, because for $d \geq \dim W$ we have 
\[
F_{\{i\}}^{d \deg v_i} (\Frob_{\{i\}}^{d \deg(v_i)})^*  \Ima(\mf{T}^\mu|_{\ol{v}} \rightarrow \mf{T}|_{\ol{v}} ) \subset \sum_{\alpha=0}^{\dim W-w} S_?  F_{\{i\}}^{\alpha \deg(v_i) } (\Frob_{\{i\}}^{\alpha \deg(v_i) })^*  \Ima(\mf{T}^\mu|_{\ol{v}} \rightarrow \mf{T} |_{\ol{v}}) .
\]

Since the $S$-operators and the $F_{\{i\}}$ are morphisms of sheaves, they commute with the specialization map $\fsp^*$. The upper arrow in \eqref{B: specialization} an isomorphism because $\ol{v} \in \Omega$ lies in the lisse locus of $\mf{T}^\mu$ by construction. Therefore the Eichler-Shimura relation from Proposition \ref{B: eichler-shimura} is also satisfied in the right column. Now over $\ol{\eta^I}$ we can apply the same elimination argument and use Theorem \ref{B: s=t} to replace S-operators by Hecke operators, thus deducing that for $d \geq \dim W$, we have 
\[
F_{\{i\}}^{d \deg v_i} (\Frob_{\{i\}}^{d \deg(v_i)})^*  \Ima(\mf{T}^\mu|_{\ol{\eta^I}} \rightarrow \mf{T}|_{\ol{\eta^I}} ) \subset \sum_{\alpha=0}^{\dim W-1}  \left( \otimes_{i \in I} \cH_{i,v_i} \right)  F_{\{i\}}^{\alpha \deg(v_i) } (\Frob_{\{i\}}^{\alpha \deg(v_i) })^*  \Ima(\mf{T}^\mu|_{\ol{\eta^I}} \rightarrow \mf{T} |_{\ol{\eta^I}}) .
\]
Therefore, we actually have
\[
\mf{M}^\mu =  \sum_{0 \leq n_i < \dim W \deg v_i}  \left( \otimes_{i \in I} \cH_{i,v_i} \right) \prod_{i \in I} F_{\{i\}}^{n_i} \left( \prod_{i \in I} (\Frob_{\{i\}}^{n_i})^* \Ima\left(\mf{T}^{\mu} |_{\ol{\eta^I}}\rightarrow  \mf{T}|_{\ol{\eta^I}}) \right) \right) .
\]
Since $\mf{T}^{\mu}|_{\ol{\eta^I}}$ is finite-dimensional over $k$, $\mf{M}^{\mu}$ is finite-type over $A^{\mu}$.
\end{proof}

\subsection{Drinfeld's Lemma}

The following result of Xue is a generalization of the so-called ``Drinfeld's Lemma''. 

\begin{lemma}[{\cite[Lemma 7.4.2]{Xue20}}]\label{B: drinfeld lemma} Let $A$ be a finitely generated $k$-algebra. Let $M$ be an $A$-module of finite type. Then any continuous $A[\FWeil(\eta^I, \ol{\eta^I})]$-action on $M$ factors through $\Weil(\eta, \ol{\eta})^I$. 
\end{lemma}

\begin{proof}[Proof of Proposition \ref{prop: FWeil on shtukas}]

Applying Lemma \ref{B: drinfeld lemma} to each $\mf{M}^\mu$, we deduce that the $\FWeil(\eta^I, \ol{\eta^I})$-action on $\mf{M}^\mu$ factors through $\Weil(\eta, \ol{\eta})^I$. Then the same holds for $\colim_{\mu} \mf{M}^\mu = T^j(\Sht_{G,D,I}; V)$. 
\end{proof}

\bibliographystyle{amsalpha}
\bibliography{Bibliography}




\end{document}